\documentclass[reqno,10pt]{amsart}
\author{Chiung-ju Liu}
\title{Bando-Futaki Invariants on Hypersurfaces}
\email{cjliu@math.ncku.edu.tw}
\address{Department of Mathematics, National Cheng Kung
University, Tainan, Taiwan 701}
\date{December 08, 2007}
\usepackage{amsfonts}
\usepackage{amsmath,amsthm,amscd}
\usepackage{showkeys}

\thanks{The author was partially supported by NSF:DMS-0202508 and NSF:DMS-0347033 during her Ph.D. study.}
\begin{document}
\begin{abstract}
In this paper, the Bando-Futaki invariants on hypersurfaces are
derived in terms of the degree of the defining polynomials, the
dimension of the underlying projective space, and the given
holomorphic vector field. In addition, the holomorphic invariant
introduced by Tian and Chen (Ricci Flow on K\"ahler-Einstein
surfaces) is proven to be the Futaki invariant on compact K\"ahler
manifolds with positive first Chern class.
\end{abstract} \maketitle \tableofcontents\
\section{Introduction}
In 1983, Futaki introduced the well-known Futaki invariant
\cite{FA1}, which is an obstacle to the existence of
K\"ahler-Einstein metrics on a compact complex manifold with
positive first Chern class. Other generalizations of the Futaki
invariant were introduced later, all of which are obstructions to
certain geometric structures. The Calabi-Futaki invariant
\cite{CE} is an obstruction to the existence of K\"ahler metrics
of constant scalar curvature on a compact K\"ahler manifold. The
Bando-Futaki invariants, raised by Bando \cite{BS} in 1983, are
obstructions to the harmonicity of the higher order Chern forms.
The Bando-Futaki invariants vanish if and only if the short-term
solutions of the almost K\"ahler-Einstein exist (cf. ~Leung
\cite{LNC}). Tian and Zhu found a holomorphic invariant
\cite{T-Z}, which is an obstruction to the existence of
K\"ahler-Ricci soliton. Recently, Futaki \cite{FA3} generalized
the Bando-Futaki invariants and the Futaki-Morita invariants
\cite{F-M}. The new invariants give obstructions to asymptotic
Chow semi-stability when the invariant polynomials are Todd
polynomials.
\par
Efficient methods for computing the Futaki invariant and the
generalized Futaki invariants are essential to characterizing the
existence of certain geometric structures. Lu \cite{LZ1}
constructed a formula to evaluate the Futaki invariant on complete
intersections. The formula depends on the dimension of the
projective space, the degree of the defining polynomials, and the
given tangent holomorphic vector field. Concurrently, Yotov
\cite{YM} derived the same result with a different approach.  On
complete intersections, Phong and Sturm \cite{P-S2} formularized
the Futaki invariant and the Mabuchi energy functional using the
Deligne pairing. Their methods may lead to a complete solution to
the problem of computing the Futaki invariant.
\par
The main part of this paper is the computation of the Bando-Futaki
invariants on hypersurfaces in $\mathbb{CP}^{n}$. The value,
stated in Theorem \ref{1.1}, is in terms of the dimension $n$, the
degree of the defining polynomial of the hypersurface, and the
tangent vector field. In Theorem \ref{1.2}, we prove that Chen and
Tian's holomorphic invariants introduced in \cite[section 5]{C-T}
are the Futaki invariants. In the last section of this paper, we
study two properties of the higher order K-energy functionals. The
first property states that the higher order K-energy functionals
are independent of the choice of paths. The second property states
that they are the nonlinearizations of the Bando-Futaki
invariants. Both properties are known to experts and are proved in
\cite{BS, B-M, WB}. We reiterate the proof for the former property
in detail with an approach different from Weinkove's \cite{WB}. We
slightly generalize the condition \cite[Theorem 2]{WB} of the
latter property.
\par
Let $M$ be an $m$-dimensional compact complex manifold with
positive first Chern class $c_{1}(M)>0$. Let $\omega_{M}$ be a
K\"ahler form on $M$ such that $[\omega_{M}]\in c_{1}(M)$. Let
$\Theta$ be the curvature tensor with respect to the K\"ahler
metric $\omega_{M}$ on the holomorphic tangent bundle
$T^{1,0}(M)$. Let $c_{q}(\Theta)=P^{q}(\Theta)$ denote the $q$-th
Chern form of $\Theta$, where $P^{q}(\Theta)$ is the elementary
polynomials defined by
\[\det(\Theta+t\cdot I)=\sum_{q=0}^{m}P^{q}(\Theta)\cdot t^{m-q}.\]
Define the Chern classes $c_{q}(M)$ by
\[c_{q}(M)=\left[P^{q}(\Theta)\right]\in H_{DR}^{2q}(M),\]
where we set $c_{0}(M)=[1]\in H_{DR}^{0}(M)$. Let $Hc_{q}(\Theta)$
be the harmonic part of $c_{q}(\Theta)$ in the Hodge
decomposition. Since $M$ is K\"ahler, there exists a real
$(q-1,q-1)$ form $f_{q}$ such that
\numberwithin{equation}{section}
\begin{equation}c_{q}(\Theta)-Hc_{q}(\Theta)
=\frac{\sqrt{-1}}{2\pi}\partial\bar{\partial}f_{q},\label{1}\end{equation}
where $f_{q}$ is unique up to a $\partial\bar{\partial}$-closed
form.
\newtheorem{definition}{Definition}[section]
\begin{definition} Let
$\mathcal{F}_{q}:H^{0}(M,T^{1,0}(M))\longrightarrow \mathbb{C}$. The
$q$-th Bando-Futaki invariant is defined as
\begin{equation}\mathcal{F}_{q}(X)=\int_{M}\textrm{L}_{X}f_{q}\wedge
\omega_{M}^{m+1-q}\label{2}\end{equation} for $q=1,\cdots,m$.
\end{definition}
\par Each $\mathcal{F}_{q}$ is well-defined on the Lie algebra of holomorphic vector
fields and independent of the choice of the K\"ahler form in the
K\"ahler class. This property was proved by Bando \cite{BS} and
can also be found in Futaki's book \cite{FA2}. In particular,
$\mathcal{F}_{1}$ is known as the Futaki invariant.
\newtheorem{theorem}{Theorem}[section]
\begin{theorem}\label{1.1}
Let $M$ be a hypersurface in $\mathbb{CP}^{n}$ defined by a
homogeneous polynomial $F$ of degree $d$ with $d\leq n$. Let $X$
be a holomorphic vector field on $\mathbb{CP}^{n}$ such that
\[XF=\kappa F\] for a constant $\kappa$. Then the $q$-th
Bando-Futaki invariant is
\begin{equation*}
\mathcal{F}_{q}(X)=-(n+1-d)^{n-q}\frac{(d-1)(n+1)}{n}\sum_{j=0}^{q-1}(-d)^{j}(j+1){n\choose q-j-1}\kappa.
\end{equation*}
\end{theorem}
A hypersurface $M$ defined by the homogeneous polynomial $F=0$ of
degree $d$ has positive first Chern class if $d\leq n$. In fact, $M$
has nonzero holomorphic vector fields.
\par A summary of the proof is as follows. The first step is to
find the potential forms $f_{q}$ defined in \eqref{1} for $1\leq q
\leq n-1$. In order to do this, we compute the elementary
symmetric polynomials by using the curvature tensors of the
hypersurface in terms of local coordinates. Then we find the extra
holomorphic forms so that the potentials $f_{q}$ can be expressed
globally. The second step is to evaluate the Bando-Futaki
invariants. We take the contraction of \eqref{1} with vector field
$X$. Then we may write the resulting function as a
$\bar{\partial}$-equation of $(q-1,q-1)$-forms. That is,
\[\bar{\partial}[-q\tilde{P}^{q}(\nabla
X,\Theta,\cdots,\Theta)+q\alpha_{qq}\theta
\omega^{q-1}-i(X)\partial f_{q}]=0,\] where $\theta$ is the
Hamiltonian function of $\omega$, $\tilde{P}^{q}$ is the
polarization of the $q$-th elementary symmetric polynomial, and
$\alpha_{qq}$ is a constant. We denote that
$\omega_{M}=(n+1-d)\omega=(n+1-d)\omega_{FS}|_{M}$, where
$\omega_{FS}$ is the Fubini-Study metric on $\mathbb{CP}^{n}$. By
Hodge decomposition, we have
\[-q\tilde{P}^{q}(\nabla
X,\Theta,\cdots,\Theta)+q\alpha_{qq}\theta
\omega^{q-1}-i(X)\partial
f_{q}=\psi_{q}+\bar{\partial}\varphi_{q},\] where $\psi_{q}$ is
the harmonic part and $\bar{\partial}\varphi_{q}$ is the exact
part. In Lemma \ref{Lef}, we show that $\psi_{q}=C(q)\omega^{q-1}$
for some constant $C(q)$. In Lemma \ref{lem2.7}, we prove that
\begin{equation}
q\int_{M}\tilde{P}^{q}(\nabla X,\Theta,\cdots,\Theta)\wedge
\omega_{M}^{n-q}=0.\label{101}
\end{equation} Then we reduce
\eqref{2} to
\begin{eqnarray*}
\mathcal{F}_{q}(X)&=&\int_{M}L_{X}(f_{q})\wedge
\omega_{M}^{n-q}\\&=&
q\alpha_{qq}(n+1-d)^{n-q}\int_{M}\theta\omega^{n-1}-C(q)(n+1-q)^{n-q}\int_{M}\omega^{n-1}.\end{eqnarray*}
By \cite[Theorem 5.1]{LZ1}, we have
\begin{align*}
\int_{M}\theta\omega^{n-1}=\frac{\kappa}{n}\qquad\text{and}\qquad\int_{M}\omega^{n-1}=d.
\end{align*}
Thus the Bando-Futaki invariant is
\[\int_{M}L_{X}f_{q}\wedge\omega^{n-q}=(n+1-d)^{n-q}\left(\frac{q\alpha_{qq}\kappa}{n}-C(q)d\right).\]
It can be computed explicitly. We have overcome two obstacles in
our proof. The first obstacle is to show the equality holds in
\eqref{101}. Note that it is not true for any invariant polynomial
$\varphi$ such that
\[\int_{M}\varphi(\nabla X, \Theta, \cdots,\Theta)\wedge\omega^{n-q}=0.
\footnote{If $\varphi(\nabla X,\Theta,\cdots,\Theta)\wedge
\omega^{n-q}=\rm{div}X \cdot (c_{1}(\Theta))^{n-1}$, then
$\int_{M}\varphi(\nabla,\Theta,\cdots,\Theta)\wedge
\omega^{n-q}\neq 0$}.\] For $q=1$, it is trivial to see
\[\int_{M}\tilde{P}^{1}(\nabla X)\wedge\omega_{M}^{n-1}=\int_{M}{\rm
div}(X)\omega_{M}^{n-1}=\int_{M}L_{X}(\omega_{M}^{n-1})=0.\]
However, equation \eqref{101} is an expected result for the
polarization $\tilde{P}^{q}(\nabla X,\Theta,\cdots,\Theta)$ of the
elementary polynomial whenever $2\leq q\leq n-1$. Without this
result, it would be very difficult to find the formula. The second
obstacle is to compute $C(q)$ and write $\varphi_{q}$ as a
globally defined forms.
 \newtheorem{corollary}{Corollary}[section]
 \begin{corollary}[Lu\cite{LZ1}] Given the conditions of Theorem 1.1 and $q=1$,
 the first Bando-Futaki invariant is the same as the Futaki invariant given as
 \[\mathcal{F}_{1}(X)=-(n+1-d)^{n-1}\frac{(n+1)(d-1)}{n}\kappa.\]
\end{corollary}
Our method can be applied, in principle, to the case of complete
intersections and toric varieties (Mabuchi). However, the
computation is complete using the current notations. Therefore our
goal is to create an abstract setting from our current ideas
before applying the invariant formula to more general cases.
\par In \cite{FA2}, the Futaki invariants (and Bando-Futaki
invariants) are formularized as equivariant Chern numbers, which
can be computed using Atiyah-Bott localization formula. The method
used in this paper is completely different from Atiyah and Bott's
method and is an alternative way to compute the equivariant Chern
numbers. We believe our method can be applied to infinite
dimensions; for example, in Gromov-Witten invariants.
\newtheorem{remark}{Remark}[section]
\begin{remark}
All Bando-Futaki invariants on a hypersurface vanish if the
hypersurface is K-semistable.
\end{remark}
Here we re-state the definition \cite{LZ2}:
\begin{definition}
We say $M$ is K semistable if any holomorphic tangent vector field
$X$ on $M$ satisfies
\begin{equation}
\lim_{t\rightarrow 0}t\frac{d}{dt}M(\omega_{0},\omega_{t})\geq
0,\label{3}
\end{equation}
where $M(\omega_{0},\omega_{t})$ is the K energy with respect to $\omega$ and $\omega_{t}$ (definition is given in section 4), and $\sigma(t)^{\ast}\omega_{0}=\omega_{t}$, where the one parametor family of automorphism $\sigma(t)$ is generated by the holomorphic vector field $X$.
\end{definition}
\par
In section 4, we study the holomorphic invariants that were
introduced by Chen and Tian \cite{C-T}.
\begin{definition}
Let $M$ be an $n$-dimensional simply-connected K\"ahler manifold
with a K\"ahler form $\omega$. Since $M$ is simply-connected, there
exists a smooth function $\theta_{X}$ such that
$i(X)\omega=\frac{\sqrt{-1}}{2\pi}\bar{\partial}\theta_{X}$\footnote{In
order to keep the original definition, we use the equation
$\theta_{X}=-\theta$, which is the opposite sign from the
Hamiltonian function $\theta$ defined as above.}. Define a new
holomorphic invariant to be
\begin{eqnarray*}
\lefteqn{F_{k}(X,\omega)}\\&&=
(n-k)\int_{M}\theta_{X}\omega^{n}+(k+1)\int_{M}\Delta\theta_{X}Ric(\omega)^{k}\wedge\omega^{n-k}\\&&-(n-k)\int_{M}\theta_{X}Ric(\omega)^{k+1}\wedge\omega^{n-k-1}.
\end{eqnarray*}
\end{definition}
\begin{theorem}\label{1.2}
If $M$ is compact with a K\"ahler form $\omega\in c_{1}(M)>0$, then
the Chen-Tian's holomorphic invariants are the Futaki invariants:
\begin{equation*}
F_{k}(X,\omega)=(k+1)\int_{M}X(f_{\omega})\omega^{n},
\end{equation*}
where $f_{\omega}$ is a potential function such that
$Ric(\omega)-\omega=\frac{\sqrt{-1}}{2\pi}\partial\bar{\partial}f_{\omega}$
\end{theorem}
The K-energy is a nonlinearization of the Futaki invariant. And also
it is mentioned in \cite{C-T} that the nonlinearizations of these
holomorphic invariants are
\[E_{k,\omega}=E_{k,\omega}^{0}(\varphi)-J_{k,\omega},\]
where
\[E_{k,\omega}^{0}(\varphi)=\frac{1}{\int_{M}\omega^{n}}
\int_{M}\left(\log\frac{\omega_{\varphi}^{n}}{\omega^{n}}
-f_{\omega}\right)\left(\sum_{i=0}^{k}Ric(\omega_{\varphi})^{i}\wedge\omega^{k-i}\right)\wedge\omega_{\varphi}^{n-k},\]
and \[J_{\omega}(\varphi)=
\frac{1}{\int_{M}\omega^{n}}\sum_{i=0}^{n-1}\int_{M}\frac{i+1}{n+1}\frac{\sqrt{-1}}{2\pi}
\partial\varphi\wedge\bar{\partial}\varphi\wedge\omega^{i}\wedge_{\varphi}^{n-1-i},\]
where $\omega_{\varphi}=\omega+\frac{\sqrt{-1}}{2\pi}\partial\bar{\partial}\varphi>0$
for some smooth function $\varphi$ and $0\leq k\leq n$.
We can see that the Futaki invariants can have different nonlinearizations.
\paragraph{\textbf{Acknowledgements}}This paper will be part of my PhD. thesis.
The author thanks her advisor Z. Lu for his mathematical insights
and assistance during the preparation of this paper. She also
thanks Professor Phong for his encouragement and Professor Paul
for pointing out the effect of K semistable hypersurfaces on
Bando-Futaki invariants.
\section{Bando-Futaki invariants}
\subsection{Curvature Tensor on Hypersurfaces}
The following setting is adopted from \cite{TG1} and \cite{LZ2}.
Let $Z=[Z_{0},\cdots, Z_{n}]$ be the homogeneous coordinate of
$\mathbb{CP}^{n}$. Let $U_{0}=\{Z\in \mathbb{CP}^{n}|Z_{0}\neq
0\}$ and $z=(z_{1},\cdots,z_{n})$, where $z_{i}=Z_{i}/Z_{0}$ for
$i=1,\cdots,n$. Without loss of generality, assume that we work on
the coordinate chart $(U_{0}, z)$. Under this coordinate system,
the Fubini-Study metric is given by
\[\omega_{FS}=\frac{\sqrt{-1}}{2\pi}\sum_{i,j=1}^{n}g_{i\bar{j}}dz_{i}\wedge d\bar{z}_{j}
=\frac{\sqrt{-1}}{2\pi}\sum_{i,j=1}^{n}
\left(\frac{\delta_{ij}}{1+|z|^{2}}-\frac{\bar{z_{i}}z_{j}}{(1+|z|^{2})^{2}}\right)dz_{i}\wedge
d\bar{z}_{j},\] where $|z|^{2}=\sum_{i=1}^{n}|z_{i}|^{2}$. Let $f$
be the defining polynomial of $M$ on $U_{0}$, where
\[f(z)=F\left[1,\frac{Z_{1}}{Z_{0}},\cdots,
\frac{Z_{n}}{Z_{0}}\right]=\frac{1}{Z_{0}^{d}}F[Z_{0},\cdots,Z_{n}].\]
By the implicit function theorem, there exists an open set
$V\subseteq\{(z_{2},\cdots,z_{n})\in\mathbb{C}^{n-1}\}$ of $U_{0}$
such that $z'$ solves the equation $f(z_{1},\cdots,z_{n})=0$
whenever $\frac{\partial f}{\partial z_{1}}(z)\neq 0$. Let
$z_{1}=z'(z_{2},\cdots,z_{n})$. Therefore, $z_{1}$ is a
holomorphic function on $V$. Under the coordinate system
$(V,(z_{2},\cdots, z_{n}))$, let the restricted K\"ahler form be
\[\omega=\omega_{FS}\Big|_{M}=\sum_{i,j=2}^{n}\frac{\sqrt{-1}}{2\pi}\tilde{g}_{i\bar{j}}dz_{i}\wedge
d\bar{z}_{j},\] where
\begin{equation}
\tilde{g}_{i\bar{j}}=\frac{\delta_{ij}
+a_{i}\bar{a}_{j}}{1+|z|^{2}}-\frac{(\bar{z}_{i}+\bar{z}_{1}a_{i})(z_{j}+z_{1}\bar{a}_{j})}{(1+|z|^{2})^{2}}
\qquad\text{for $i,j=2,\cdots, n$},\label{4}\end{equation} where
$a_{i}=\frac{\partial z_{1}}{\partial z_{i}}$ for $i=2,\cdots,n$.
In order to compute the curvature tensor with respect to the
metric $\tilde{g}_{i\bar{j}}$, it is critical to find the inverse
matrix of $\tilde{g}_{i\bar{j}}$.
\newtheorem{lemma}{Lemma}[section]
\begin{lemma}\label{lem2.1}
Using the same notation as above, the inverse of
$\tilde{g}_{i\bar{j}}$ is
\begin{eqnarray*}
\tilde{g}^{i\bar{j}}&=&\frac{1}{\rho}
\left(\rho(1+|z|^{2})\delta_{ji}-a_{j}\bar{a}_{i}+\bar{z}_{j}z_{i}(1+|a|^{2})\right)\\&&
- \frac{1}{\rho}\left(a_{j}z_{i}
(\sum_{k=2}^{n}\bar{a}_{k}\bar{z}_{k}-\bar{z}_{1})+\bar{z}_{j}\bar{a}_{i}(\sum_{k=2}^{n}a_{k}z_{k}-z_{1})\right),
\end{eqnarray*}
where $
\rho=\frac{\sum_{k=0}^{n}|F_{k}|^{2}}{(1+|z|^{2})|F_{1}|^{2}}$,
$|a|^{2}=\sum_{i=2}^{n}|a_{i}|^{2}$, and $F_{k}=\frac{\partial
F}{\partial Z_{k}}$ for $k=0,1\cdots,n$.
\end{lemma}
\begin{proof}
Consider $\tilde{g}_{i\bar{j}}$ as a matrix $A_{ij}$ for $2\leq
i,j\leq n$. Since $\tilde{g}_{i\bar{j}}$ is a matrix of a linear
combination of matrices $\delta_{ij}$, $a_{i}\bar{a}_{j}$,
$\bar{z}_{i}z_{j}$, $a_{i}z_{j}$, and $\bar{z}_{i}\bar{a}_{j}$
pointwisely, its adjugate and inverse can be written as linear
combinations of $\delta_{ij}$, $a_{i}\bar{a}_{j}$,
$\bar{z}_{i}z_{j}$, $a_{i}z_{j}$, and $\bar{z}_{i}\bar{a}_{j}$
pointwisely. More clearly, let $A_{ij}=
(\gamma_{1}\delta_{ij}+\gamma_{2}a_{i}\bar{a}_{j}+\gamma_{3}\bar{z}_{i}z_{j}
+\gamma_{4}a_{i}z_{j}+\gamma_{5}\bar{z}_{i}\bar{a}_{j})$ where
\begin{multline}
(\gamma_{1},\gamma_{2},\gamma_{3},\gamma_{4},\gamma_{5})\\=
\left(\frac{1}{1+|z|^{2}},\frac{1+|z|^{2}-|z_{1}|^{2}}{(1+|z|^{2})^{2}}
,-\frac{1}{(1+|z|^{2})^{2}},-\frac{\bar{z}_{1}}{(1+|z|^{2})^{2}},
-\frac{z_{1}}{(1+|z|^{2})^{2}}\right).\label{2.2}\end{multline}

By definition, the adjugate matrix is given by
\begin{equation*}({\rm adj}A)_{ji}=(-1)^{i+j} \sum_{\sigma\in
S_{n-2}}{\rm sgn}(\sigma) A_{i_{\sigma(3)}j_{3}}\cdots
A_{i_{\sigma (n)}j_{n}},\end{equation*}
 where $S_{n-2}$ are all permutations of
 $\{i_{3},\cdots,i_{n}\}=\{2,\cdots, n\}-\{i\}$ and \newline $\{j_{3},j_{4},\cdots,j_{n}\}=\{2,\cdots,n\}-\{j\}$.
 There exists functions $\eta_{1},\cdots,\eta_{5}$  such that
\[({\rm adj}A)_{ji}=\left(\eta_{1}\delta_{ji}+\eta_{2}a_{j}\bar{a}_{i}+\eta_{3}\bar{z}_{j}z_{i}
+\eta_{4}a_{j}z_{i}+\eta_{5}\bar{z}_{j}\bar{a}_{i}\right).
\]
Apply the following formulas given in \cite{LZ2}:
\begin{align}
\begin{split}
\sum_{k=2}^{n}a_{k}z_{k}-z_{1}=\frac{F_{0}}{F_{1}},\hspace{3cm}\\
1+|a|^{2}+\big|\sum_{k=2}^{n}a_{k}z_{k}-z_{1}\big|^{2}=\frac{\sum_{k=0}^{n}|F_{k}|^{2}}{|F_{1}|^{2}},\\
\det\tilde{g}_{i\bar{j}}=\frac{1}{(1+|z|^{2})^{n}}\frac{\sum_{k=0}^{n}|F_{k}|^{2}}{|F_{1}|^{2}}.\hspace{1.5cm}
\end{split}
\label{6}\end{align} Then we can compute the coefficients
$\eta_{1},\cdots, \eta_{5}$ by solving the following linear
equation system
\begin{eqnarray}
\tilde{g}_{i\bar{j}}\tilde{g}^{k\bar{j}}&=&\tilde{g}_{i\bar{j}}\frac{1}{\det{\tilde{g}_{i\bar{j}}}}({\rm
adj}A)_{jk}\nonumber\\&=&
\frac{1}{\det{\tilde{g}_{i\bar{j}}}}(\gamma_{1}\delta_{ij}+\gamma_{2}a_{i}\bar{a}_{j}+\gamma_{3}\bar{z}_{i}z_{j}
+\gamma_{4}a_{i}z_{j}+\gamma_{5}\bar{z}_{i}\bar{a}_{j})\nonumber\\&&
\times(\eta_{1}\delta_{jk}+\eta_{2}a_{j}\bar{a}_{k}
+\eta_{3}\bar{z}_{j}z_{k}+\eta_{4}a_{j}z_{k}+\eta_{5}\bar{z}_{j}\bar{a}_{k})\nonumber\\&=&
\delta_{ik}.\label{2.3}
\end{eqnarray}
Furthermore,  we may rewrite \eqref{2.3} as
\begin{align}
\begin{split}
(\gamma_{1}+\gamma_{2}|a|^{2}+\gamma_{4}\sum_{j=2}^{n}a_{j}z_{j})\eta_{2}
+(\gamma_{2}\sum_{j=2}^{n}\bar{a}_{j}\bar{z}_{j}+\gamma_{4}\sum_{j=2}^{n}|z_{j}|^{2})\eta_{5}=&-\gamma_{2}\eta_{1}\\
(\gamma_{1}+\gamma_{3}\sum_{j=2}^{n}|z_{j}|^{2}+\gamma_{5}\sum_{j=2}^{n}\bar{a}_{j}\bar{z}_{j})\eta_{3}+
(\gamma_{5}|a|^{2}+\gamma_{3}\sum_{j=2}^{n}a_{j}z_{j})\eta_{4}=&-\gamma_{3}\eta_{1}\\
(\gamma_{2}\sum_{j=2}^{n}\bar{a}_{j}\bar{z}_{j}+\gamma_{4}\sum_{j=2}^{n}|z_{j}|^{2})\eta_{3}
+(\gamma_{1}+\gamma_{2}|a|^{2}+\gamma_{4}\sum_{j=2}^{n}a_{j}z_{j})\eta_{4}=&-\gamma_{4}\eta_{1}\\
(\gamma_{3}\sum_{j=2}^{n}a_{j}z_{j}+\gamma_{5}|a|^{2})\eta_{2}
+(\gamma_{1}+\gamma_{3}\sum_{j=2}^{n}|z_{j}|^{2}+\gamma_{5}\sum_{j=2}^{n}\bar{a}_{j}\bar{z}_{j})\eta_{5}=
&-\gamma_{5}\eta_{1},\label{2.4}
\end{split}
\end{align}
where $\eta_{1}=
\det{\tilde{g}_{i\bar{j}}}\frac{1}{\gamma_{1}}=(1+|z|^{2})$. Hence
by using \eqref{2.2} and \eqref{2.3}, the solutions for
\eqref{2.4} are
\[(\eta_{2},\eta_{3},\eta_{4},\eta_{5})=(1+|z|^{2})^{n-1}
\left(-1,1+|a|^{2},-\big(\sum_{k=2}^{n}\bar{a}_{k}\bar{z}_{k}-\bar{z}_{1}\big),
-\big(\sum_{k=2}^{n}a_{k}z_{k}-z_{1}\big)\right),\] where
$\rho=\frac{\sum_{k=0}^{n}|F_{k}|^{2}}{(1+|z|^{2})|F_{1}|^{2}}$.
Then we may obtain
\[\tilde{g}^{i\bar{j}}=\frac{1}{\det{\tilde{g}_{i\bar{j}}}}\left(\eta_{1}\delta_{ji}+\eta_{2}a_{j}\bar{a}_{i}+\eta_{3}\bar{z}_{j}z_{i}
+\eta_{4}a_{j}z_{i}+\eta_{5}\bar{z}_{j}\bar{a}_{i}\right).\]\end{proof}
The following Lemma is important for computing higher order Chern
forms of the hypersurface $M$ and for evaluating the Bando-Futaki
invariants.
\begin{lemma}\label{lem2.2}
The curvature of the hypersurface is
\[\sum_{i,j=2}^{n}R_{ki\bar{j}}^{\ell}dz_{i}\wedge d\bar{z}_{j}=
\sum_{i,j=2}^{n}\left(\delta_{k\ell}\tilde{g}_{i\bar{j}}+\delta_{i\ell}\tilde{g}_{k\bar{j}}
-\frac{1}{\rho}\frac{\partial a_{k}}{\partial z_{i}}\frac{\partial
\bar{a}_{s}}{\partial
\bar{z}_{j}}\tilde{g}^{\ell\bar{s}}\right)dz_{i}\wedge
d\bar{z}_{j}\] for $2\leq k, \ell\leq n$, where
$\rho=\frac{\sum_{k=0}^{n}|F_{k}|^{2}}{(1+|z|^{2})|F_{1}|^{2}}$.
\end{lemma}
\begin{proof}Recall the curvature formula
\begin{eqnarray}
R_{ki\bar{j}}^{\ell}&=&-\bar{\partial}_{j}(\frac{\partial
\tilde{g}_{k\bar{s}}}{\partial
z_{i}}\tilde{g}^{\ell\bar{s}})\nonumber\\&=&-\frac{\partial^{2}\tilde{g}_{k\bar{s}}}{\partial
z_{i}\partial
\bar{z}_{j}}\tilde{g}^{\ell\bar{s}}+\frac{\partial\tilde{g}_{k\bar{q}}}{\partial
z_{i}}\frac{\bar{\partial}\tilde{g}_{p\bar{s}}}{\partial
\bar{z}_{j}}\tilde{g}^{p\bar{q}}\tilde{g}^{\ell\bar{s}}
\label{7}\end{eqnarray} on K\"ahler manifolds. By using the metric
$\tilde{g}_{i\bar{j}}$ defined in \eqref{4}, we obtain
\begin{eqnarray}
\frac{\partial\tilde{g}_{k\bar{q}}}{\partial z_{i}}=
-\frac{\bar{z}_{i}+\bar{z}_{1}a_{i}}{1+|z|^{2}}\tilde{g}_{k\bar{q}}
-\frac{\bar{z}_{k}+\bar{z}_{1}a_{k}}{1+|z|^{2}}\tilde{g}_{i\bar{q}}
+\frac{\partial a_{k}}{\partial
z_{i}}\frac{\big(1+|z|^{2}-|z_{1}|^{2}\big)
\bar{a}_{q}-\bar{z}_{1}z_{q}}{(1+|z|^{2})^{2}},\label{8}\\
\frac{\partial\tilde{g}_{p\bar{s}}}{\partial \bar{z}_{j}}=
-\frac{z_{j}+z_{1}\bar{a}_{j}}{1+|z|^{2}}\tilde{g}_{p\bar{s}}
-\frac{z_{s}+z_{1}\bar{a}_{s}}{1+|z|^{2}}\tilde{g}_{p\bar{j}}
+\frac{\partial \bar{a}_{s}}{\partial \bar{z}_{j}}
\frac{\big(1+|z|^{2}-|z_{1}|^{2}\big)a_{p}-z_{1}\bar{z}_{p}}{(1+|z|^{2})^{2}}.\label{9}
\end{eqnarray}
Use \eqref{9} to compute second derivative
\begin{eqnarray}
\frac{\partial^{2}\tilde{g}_{k\bar{s}}}{\partial z_{i}\partial
\bar{z}_{j}}&=&
-\tilde{g}_{i\bar{j}}\tilde{g}_{k\bar{s}}-\tilde{g}_{i\bar{s}}\tilde{g}_{k\bar{j}}\nonumber\\&&+
\frac{(z_{j}+z_{1}\bar{a}_{j})(\bar{z}_{i}+\bar{z}_{1}a_{i})}{(1+|z|^{2})^{2}}\tilde{g}_{k\bar{s}}+
\frac{(z_{j}+z_{1}\bar{a}_{j})(\bar{z}_{k}+\bar{z}_{1}a_{k})}{(1+|z|^{2})^{2}}\tilde{g}_{i\bar{s}}\nonumber\\&&+
\frac{(z_{s}+z_{1}\bar{a}_{s})(\bar{z}_{i}+\bar{z}_{1}a_{i})}{(1+|z|^{2})^{2}}\tilde{g}_{k\bar{j}}+
\frac{(z_{s}+z_{1}\bar{a}_{s})(\bar{z}_{k}+\bar{z}_{1}a_{k})}{(1+|z|^{2})^{2}}\tilde{g}_{i\bar{j}}\nonumber\\&&-
\frac{(z_{j}+z_{1}\bar{a}_{j})}{1+|z|^{2}}\frac{\partial
a_{k}}{\partial
z_{i}}\frac{\big(1+|z|^{2}-|z_{1}|^{2}\big)\bar{a}_{s}-\bar{z}_{1}z_{s}}{(1+|z|^{2})^{2}}\nonumber\\&&-
\frac{(z_{s}+z_{1}\bar{a}_{s})}{1+|z|^{2}}\frac{\partial
a_{k}}{\partial
z_{i}}\frac{\big(1+|z|^{2}-|z_{1}|^{2}\big)\bar{a}_{j}-\bar{z}_{1}z_{j}}{(1+|z|^{2})^{2}}\nonumber\\&&+\frac{\partial
a_{k}}{\partial z_{i}}\frac{\partial \bar{a}_{s}}{\partial
\bar{z}_{j}} \frac{1+|z|^{2}-|z_{1}|^{2}}{(1+|z|^{2})^{2}}+
\frac{\partial \bar{a}_{s}}{\partial
\bar{z}_{j}}\frac{a_{k}\bar{z}_{i}-a_{i}\bar{z}_{k}}{(1+|z|^{2})^{2}}\nonumber\\&&-2\frac{\partial
\bar{a}_{s}}{\partial \bar{z}_{j}}
(\bar{z}_{i}+\bar{z}_{1}a_{i})\frac{\big(1+|z|^{2}-|z_{1}|^{2}\big)a_{k}-z_{1}\bar{z}_{k}}{(1+|z|^{2})^{3}}.\label{2.6}
\end{eqnarray}
By using \eqref{8} and \eqref{9} in \eqref{2.6}, we obtain
\begin{multline}
\frac{\partial^{2}\tilde{g}_{k\bar{s}}}{\partial z_{i}\partial
\bar{z}_{j}}=
-\tilde{g}_{i\bar{j}}\tilde{g}_{k\bar{s}}-\tilde{g}_{i\bar{s}}\tilde{g}_{k\bar{j}}
+\frac{\partial\tilde{g}_{k\bar{q}}}{\partial
z_{i}}\frac{\bar{\partial}\tilde{g}_{p\bar{s}}}{\partial
\bar{z}_{j}}\tilde{g}^{p\bar{q}}+\frac{\partial a_{k}}{\partial
z_{i}}\frac{\partial \bar{a}_{s}}{\partial \bar{z}_{j}}
\frac{1+|z|^{2}-|z_{1}|^{2}}{(1+|z|^{2})^{2}}\\-\frac{\partial
a_{k}}{\partial z_{i}}\frac{\partial \bar{a}_{s}}{\partial
\bar{z}_{j}}\left(\frac{(1+|z|^{2}-|z_{1}|^{2})a_{p}-z_{1}\bar{z}_{p}}{(1+|z|^{2})^{2}}\right)
\left(\frac{(1+|z|^{2}-|z_{1}|^{2})\bar{a}_{q}-\bar{z}_{1}z_{q}}{(1+|z|^{2})^{2}}\right)\tilde{g}^{p\bar{q}}
.\label{10}
\end{multline}
In fact, we have
\begin{multline}\left(\frac{(1+|z|^{2}-|z_{1}|^{2})a_{p}-z_{1}\bar{z}_{p}}{(1+|z|^{2})^{2}}\right)
\left(\frac{(1+|z|^{2}-|z_{1}|^{2})\bar{a}_{q}-\bar{z}_{1}z_{q}}{(1+|z|^{2})^{2}}\right)\tilde{g}^{p\bar{q}}-
\frac{1+|z|^{2}-|z_{1}|^{2}}{(1+|z|^{2})^{2}}\\=
-\frac{|F_{1}|^{2}}{(1+|z|^{2})\sum_{k=0}^{n}|F_{k}|^{2}}=-\frac{1}{\rho}.\hspace{6cm}\label{2.11}\end{multline}
By using \eqref{2.11} in \eqref{10}, we get
\[-\frac{\partial^{2}\tilde{g}_{k\bar{s}}}{\partial z_{i}\partial
\bar{z}_{j}}\tilde{g}^{\ell\bar{s}}+\frac{\partial\tilde{g}_{k\bar{q}}}{\partial
z_{i}}\frac{\bar{\partial}\tilde{g}_{p\bar{s}}}{\partial
\bar{z}_{j}}\tilde{g}^{p\bar{q}}\tilde{g}^{\ell\bar{s}}=\delta_{k\ell}\tilde{g}_{i\bar{j}}+\delta_{i\ell}\tilde{g}_{k\bar{j}}
-\frac{1}{\rho}\frac{\partial a_{k}}{\partial z_{i}}\frac{\partial
\bar{a}_{s}}{\partial \bar{z}_{j}}\tilde{g}^{\ell\bar{s}}.\]
\end{proof}
The Ricci curvature of a hypersurface has been shown in \cite{LZ2,
TG1}. It is also directly followed from Lemma \ref{lem2.2}.
\begin{remark}\label{re2.1}
Given the conditions as in Theorem \ref{1.1}, the Ricci curvature on
the hypersurface is
\[Ric(\omega)=(n+1-d)\omega-\frac{\sqrt{-1}}{2\pi}\partial\bar{\partial}\xi,\]
where \[\xi=\log\left(\frac{\sum_{k=0}^{n}|F_{k}|^{2}}{(\sum_{k=0}^{n}|Z_{k}|^{2})^{d-1}}\right).\]
\end{remark}
\begin{proof}
The trace of the curvature form is the Ricci curvature. By using
Lemma \ref{lem2.2}, compute
\begin{align*}\begin{split}
\sum_{k=2}^{n}\sum_{i,j=2}^{n}R_{ki\bar{j}}^{k}dz_{i}\wedge
d\bar{z}_{j}&=
\sum_{k=2}^{n}\sum_{i,j=2}^{n}\left(\delta_{kk}\tilde{g}_{i\bar{j}}+\delta_{ik}\tilde{g}_{k\bar{j}}
-\frac{1}{\rho}\frac{\partial a_{k}}{\partial z_{i}}\frac{\partial
\bar{a}_{s}}{\partial
\bar{z}_{j}}\tilde{g}^{k\bar{s}}\right)dz_{i}\wedge
d\bar{z}_{j}\\
&=\sum_{i,j=2}^{n}\left(n\tilde{g}_{i\bar{j}}-\frac{1}{\rho}\sum_{k,s=2}^{n}\frac{\partial
a_{k}}{\partial z_{i}} \frac{\partial\bar{a}_{s}}{\partial
\bar{z}_{j}}\tilde{g}^{k\bar{s}}\right)dz_{i}\wedge
d\bar{z}_{j},\end{split}\end{align*} where
\[\frac{1}{\rho}\frac{\sqrt{-1}}{2\pi}\sum_{k,s=2}^{n}\frac{\partial
a_{k}}{\partial z_{i}}\frac{\partial \bar{a}_{s}}{\partial
\bar{z}_{j}}\tilde{g}^{k\bar{s}}=\partial_{i}\bar{\partial}_{j}\log
\left(\frac{\sum_{k=0}^{n}|F_{k}|^{2}}{|F_{1}|^{2}}\right).\]
 We can find an extra
globally defined function such that
\begin{align}\begin{split}
\frac{\sqrt{-1}}{2\pi}\partial\bar{\partial}
\log\left(\frac{\sum_{k=0}^{n}|F_{k}|^{2}}{|F_{1}|^{2}}\right)&=\frac{\sqrt{-1}}{2\pi}\partial\bar{\partial}
\log\left(\frac{\sum_{k=0}^{n}|F_{k}|^{2}}{(\sum_{k=0}^{n}|Z_{k}|^{2})^{d-1}}\frac{(\sum_{k=0}^{n}|Z_{k}|^{2})^{d-1}}{|F_{1}|^{2}}\right)
\\&=(d-1)\omega+\frac{\sqrt{-1}}{2\pi}\partial\bar{\partial}\xi.\label{10.1}\end{split}\end{align}
\end{proof}
\subsection{Chern Forms}
\begin{lemma}\label{lem2.3} Given the conditions in Theorem \ref{1.1}, the $q$-th Chern form on a hypersurface are
\[c_{q}(\Theta)=\sum_{k=0}^{q}\alpha_{qk}\omega^{k}\wedge(\frac{\sqrt{-1}}{2\pi}
\partial\bar{\partial}\xi)^{q-k},\] where
\begin{eqnarray*}\alpha_{qq}&=&{n+1\choose
q}-d\alpha_{(q-1)(q-1)},\\
\alpha_{q(q-k)}&=&-[d\alpha_{(q-1)(q-k-1)}+\alpha_{(q-1)(q-k)}]\qquad\text{for $k=1,\cdots,q-1$}\\
\alpha_{q0}&=&(-1)^{q},
\end{eqnarray*}
for $q=1,\cdots,n-1$. Set that $\alpha_{00}=1$.
\end{lemma}
Before we prove this Lemma, we need several steps. Let
$\Theta=\sum_{k,\ell=2}^{n}\Theta_{k}^{\ell}\frac{\partial}{\partial
z_{\ell}}\otimes dz_{k}$ be the curvature tensor of $M$ such that
\[\Theta_{k}^{\ell}=
\frac{\sqrt{-1}}{2\pi}\sum_{i,j=2}^{n}R_{ki\bar{j}}^{\ell}dz_{i}\wedge
d\bar{z}_{j}\] represents a $(1,1)$-form valued matrix for $2\leq
k, \ell\leq n$. In order to save some space, we denote the
$(k,k)$-form valued matrix
$\sum_{i_{2},\cdot,i_{k}=2}^{n}\Theta_{i_{1}}^{i_{2}}\wedge\cdots\wedge\Theta_{i_{k}}^{i_{k+1}}$
by
$\sum_{i_{2},\cdot,i_{k}=2}^{n}\Theta_{i_{1}}^{i_{2}}\cdots\Theta_{i_{k}}^{i_{k+1}}$.
The $q$-th Chern polynomial defined in \cite[page 402, 417]{G-H}
is
\begin{equation}P^{q}(\Theta)=\frac{1}{q!}\sum_{\sigma\in
S_{q}}\sum_{i_{1},\cdots,i_{q}=2}^{n}{\rm
sgn}(\sigma)\Theta_{i_{1}}^{i_{\sigma(1)}}\Theta_{i_{2}}^{i_{\sigma(2)}}\cdots\Theta_{i_{q}}^{i_{\sigma(q)}}
.\label{11}
\end{equation}
In Sublemma \ref{sub2.1}, we derive a formula to compute higher
order Chern polynomials in terms of the lower order Chern
polynomials. The idea is as follows. Let
\[\phi_{j}
=\sum_{i_{1},\cdots,i_{j}=2}^{n}\Theta_{i_{1}}^{i_{2}}\Theta_{i_{2}}^{i_{3}}\cdots\Theta_{i_{j}}^{i_{1}}\]
be a $(j,j)$-form corresponding to a cycle of order $j$. For each
permutation $\sigma\in S_{q}$ not the identity is a product of
cyclic permutations which operate on disjoint indices:
$\sigma=\sigma_{1}\sigma_{2}\cdots\sigma_{k}$. If
$\sigma_{1},\cdots,\sigma_{k}$ are cycles of length
$j_{1},\cdots,j_{k}$ with respectively, then
$\Theta_{i_{1}}^{i_{\sigma(1)}}\cdots\Theta_{i_{q}}^{i_{\sigma(q)}}=\phi_{j_{1}}\cdots\phi_{j_{k}}$.
For each $1\leq \ell\leq k$, we have
$\sigma=\sigma_{\ell}\tau_{\ell}$ for the permutation
$\tau_{\ell}=\sigma_{1}\cdots\sigma_{\ell-1}\sigma_{\ell+1}\cdots\sigma_{k}\in
S_{q-j_{\ell}}$. Then ${\rm
sgn}(\tau_{\ell})\Theta_{i_{1}}^{i_{\tau_{\ell}(1)}}\cdots\Theta_{i_{q-j_{\ell}}}^{i_{\tau_{\ell}(q-j_{\ell})}}$
is one term in the expansion of $P^{q-j_{\ell}}(\Theta)$ as in
\eqref{11}. We abuse the notation that $\tau_{\ell}$ operates on
$\{1,\cdots,q-j_{\ell}\}$. In fact, $\tau_{\ell}$ operates on
$q-j_{\ell}$ many indices different from $\sigma_{\ell}$.
Therefore,
$\Theta_{i_{1}}^{i_{\sigma(1)}}\cdots\Theta_{i_{q}}^{i_{\sigma(q)}}$
can be written as one term in the expansion of ${\rm
sgn}(\sigma_{1})\phi_{j_{1}}P^{q-j_{1}}(\Theta),\cdots$, and ${\rm
sgn}(\sigma_{k})\phi_{j_{k}}P^{q-j_{k}}(\Theta)$. Then we conclude
that $P^{q}(\Theta)=
\sum_{j=1}^{q}(-1)^{j-1}a_{j}\phi_{j}P^{q-j}(\Theta)$ for some
rational number $a_{j}$ for $j=1,\cdots, q$.
\newtheorem{sublem}{Sublemma}[section]
\begin{sublem}\label{sub2.1}
The $q$-th Chern polynomial can be written as
\begin{equation}
P^{q}(\Theta)= \sum_{j=1}^{q}(-1)^{j-1}a_{j}\phi_{j}P^{q-j}(\Theta)
=\frac{1}{q}\sum_{j=1}^{q}(-1)^{j-1}\phi_{j}P^{q-j}(\Theta)\label{12}\end{equation}
for $q=1,\cdots,n-1$. Let $P^{0}(\Theta)=1$ for convention.
\end{sublem}
\begin{proof}Prove \eqref{12} by induction. For $q=1$,
$P^{1}(\Theta)=\sum_{i=1}^{n}\Theta_{i}^{i}=\phi_{1}.$ Assume that
$P^{k}(\Theta)=\frac{1}{k}\sum_{j=1}^{k}(-1)^{j-1}\phi_{j}P^{k-j}(\Theta)$
holds for $k=2,\cdots,q-1$. Actually, there are ${q\choose
j}(j-1)!$ many cycles of order $j$ in $S_{q}$. For any positive
integer numbers $j_{1},\cdots,j_{k}$ such that
$j_{1}+\cdots+j_{k}=q$, the coefficient for
$\phi_{j_{1}}\cdots\phi_{j_{k}}$ on the left hand side of
\eqref{12} is
\begin{equation*}  \left\{
\begin{aligned}
    \frac{1}{q!}{q\choose j_{1}}(j_{1}-1)!{q-j_{1}\choose
j_{2}}(j_{2}-1)!\cdots{j_{k}\choose j_{k}}(j_{k}-1)!, \quad\text{if $j_{t}\neq j_{s}$, $1\leq s,t\leq k$ ;} \\
    \frac{1}{\ell!q!}{q\choose
j_{1}}(j_{1}-1)!{q-j_{1}\choose
j_{2}}(j_{2}-1)!\cdots{j_{k}\choose j_{k}}(j_{k}-1)!, \text{if
$\exists\ell\leq k$, $j_{1}=\cdots =j_{\ell}$, ;}
\end{aligned}\right.
\end{equation*} and so on. Therefore, the coefficient for
$\phi_{q}$ is $a_{q}=\frac{1}{q!}{q\choose q}(q-1)!=\frac{1}{q}$.
For $q>j>q/2$, the coefficient for $\phi_{j}\phi_{q-j}$ on the
left hand side of \eqref{12} is $\frac{1}{q!}{q\choose
j}(j-1)!(q-j-1)!=\frac{1}{j(q-j)}$. On the other hand, since the
term $\phi_{j}\phi_{q-j}$ only appear in the expansion of
$\phi_{j}P^{q-j}(\Theta)$ and $\phi_{q-j}P^{j}(\Theta)$. By using
the induction hypothesis, compare the coefficient of
$\phi_{j}\phi_{q-j}$ on the both hand side of \eqref{12}. That is,
\begin{equation}\frac{1}{j(q-j)}=\frac{a_{j}}{q-j}+\frac{a_{q-j}}{j}.\label{12.1}\end{equation}
If $q=2k$ for some integer $k$, the coefficient of
$\phi_{k}\phi_{k}$ is
\[\frac{1}{2!q!}{q\choose k}(k-1)!(k-1)!=\frac{1}{k}a_{k}\]
since $\phi_{k}\phi_{k}$ only appear in $\phi_{k}P^{k}(\Theta)$.
Then $a_{k}=\frac{1}{2k}=\frac{1}{q}$. For $j=q-2$, since
$\phi_{q-2}\phi_{1}\phi_{1}$ only appear once in
$\phi_{q-2}P^{2}(\Theta)$ and twice in $\phi_{1}P^{q-2}(\Theta)$,
the coefficient is
\begin{equation}
\frac{1}{2(q-2)}=\frac{a_{q-2}}{2}+\frac{a_{1}}{q-1}+\frac{a_{1}}{(q-1)(q-2)}\label{12.6}
\end{equation}
Furthermore, for $q-2>j>q-j-1>1$, the term
$\phi_{j}\phi_{q-j-1}\phi_{1}$ appear twice
$\phi_{j}P^{q-j}(\Theta)$, $\phi_{q-j-1}P^{j+1}(\Theta)$ and
$\phi_{1}P^{q-1}(\Theta)$. The coefficient for
$\phi_{j}\phi_{q-j-1}\phi_{1}$ is
\begin{multline}\frac{1}{j(q-j-1)}\\
=\frac{a_{j}}{q-j}+\frac{a_{j}}{(q-j)(q-j-1)}+\frac{a_{q-j-1}}{j+1}+\frac{a_{q-j-1}}{(j+1)j}
+\frac{a_{1}}{(q-1)(q-j-1)}+\frac{a_{1}}{(q-1)j}\\=
\frac{a_{j}}{q-j-1}+\frac{a_{q-j-1}}{j}+\frac{a_{1}}{(q-j-1)j}.\hspace{6.5cm}\label{12.2}
\end{multline}
By using \eqref{12.1}, \eqref{12.6}, and \eqref{12.2}, we obtain
the relation of $a_{j}$
\begin{eqnarray}
a_{q-1}&=&t \label{12.3}\\
a_{q-j}&=&\frac{-(j-1)+j(q-1)t}{q-j}\qquad\text{for $j>q/2$}\label{12.4}\\
a_{j}&=&1-(q-1)t\qquad \text{for $j<q/2$}\label{12.5}.
\end{eqnarray}
The term $\underbrace{\phi_{1}\cdots\phi_{1}}_{q}$ only appear in
$\phi_{1}P^{q-1}(\Theta)$, the coefficient of it is
\[\frac{1}{q!q!}{q\choose 1}{q-1\choose 1}\cdots{1\choose
1}=\frac{1}{q!}=\frac{a_{1}}{(q-1)(q-2)\cdots 1}.\] We get
$a_{1}=\frac{1}{q}$. By using \eqref{12.3}, \eqref{12.4} and
\eqref{12.5}, we obtain $a_{j}=\frac{1}{q}$ for $j=2,\cdots,q$.
\end{proof}
 Therefore,
we only need to compute $\phi_{j}$ for $j=2,\cdots,q$ and use
\eqref{12} to formulate $P^{q}(\Theta)$.
\newtheorem{claim}{Claim}[section]
\begin{sublem}\label{sub2.2}
For $3\leq j\leq q$, we have
\begin{eqnarray*}\lefteqn{\sum_{i_{2},\cdots,i_{j-1}=2}^{n}\Theta_{i_{1}}^{i_{2}}\cdots\Theta_{i_{j-1}}^{i_{j}}=
\omega\sum_{i_{3},\cdots,i_{j-1}=2}^{n}\Theta_{i_{1}}^{i_{3}}\Theta_{i_{3}}^{i_{4}}\cdots\Theta_{i_{j-1}}^{i_{j}}}\\&&
\hspace{2cm}-\frac{\sqrt{-1}}{2\pi}\frac{1}{\rho}\sum_{\alpha,\beta=2}^{n}\left(\frac{\partial
a_{i_{1}}}{\partial z_{\alpha}}\frac{\partial \bar{a}_{s}}{\partial
\bar{z}_{\beta}}\tilde{g}^{i_{j}\bar{s}}dz_{\alpha}\wedge
d\bar{z}_{\beta}\right)\wedge(d\omega+\frac{\sqrt{-1}}{2\pi}\partial\bar{\partial}\xi)^{j-2},
\end{eqnarray*}  where
$\rho=\frac{\sum_{k=0}^{n}|F_{k}|^{2}}{(1+|z|^{2})|F_{1}|^{2}}$,
$\xi=\log\left(\frac{\sum_{k=0}^{n}|F_{k}|^{2}}{(\sum_{k=0}^{n}|Z_{k}|^{2})^{d-1}}\right)$,
and
$\Theta_{i}^{j}=\frac{\sqrt{-1}}{2\pi}\sum_{\alpha,\beta=2}^{n}R_{i\alpha\bar{\beta}}^{j}dz_{\alpha}\wedge
d\bar{z}_{\beta}$.
\end{sublem}
\begin{proof}
Prove by induction. For $j=3$, use Lemma \ref{lem2.2} and
\eqref{10.1} to compute
\begin{multline}
\sum_{i_{2}=2}^{n}\Theta_{i_{1}}^{i_{2}}\Theta_{i_{2}}^{i_{3}}=
\big(\frac{\sqrt{-1}}{2\pi})^{2}\sum_{i_{2},\alpha,\beta,\lambda,\eta,s,t=2}^{n}
\Big(\delta_{i_{1}i_{2}}\tilde{g}_{\alpha\bar{\beta}}+\delta_{\alpha
i_{2}}\tilde{g}_{i_{1}\bar{\beta}}-\frac{1}{\rho}\frac{\partial
a_{i_{1}}}{\partial z_{\alpha}}\frac{\partial
\bar{a}_{s}}{\partial
\bar{z}_{\beta}}\tilde{g}^{i_{2}\bar{s}}\Big)\\\hspace{1cm}
\cdot\Big(\delta_{i_{2}i_{3}}\tilde{g}_{\lambda \bar{\eta}}+
\delta_{\lambda i_{3}}\tilde{g}_{i_{2}
\bar{\eta}}-\frac{1}{\rho}\frac{\partial a_{i_{2}}}{\partial
z_{\lambda}}\frac{\partial \bar{a}_{t}}{\partial
\bar{z}_{\eta}}\tilde{g}^{i_{3}\bar{t}}\Big)dz_{\alpha}\wedge
d\bar{z}_{\beta}\wedge dz_{\lambda}\wedge d\bar{z}_{\eta}\\=
\big(\frac{\sqrt{-1}}{2\pi})^{2}\sum_{i_{2},\alpha,\beta,\lambda,\eta=2}^{n}
\bigg[\delta_{i_{1}i_{3}}\tilde{g}_{\alpha\bar{\beta}}\tilde{g}_{\lambda
\bar{\eta}}+\delta_{i_{3}\lambda}\tilde{g}_{\alpha\bar{\beta}}\tilde{g}_{i_{1}
\bar{\eta}}\hspace{1cm}\\-
\frac{1}{\rho}\sum_{t=2}^{n}\frac{\partial a_{i_{1}}}{\partial
z_{\lambda}}\frac{\partial \bar{a}_{t}}{\partial
\bar{z}_{\eta}}\tilde{g}^{i_{3}\bar{t}}\tilde{g}_{\alpha\bar{\beta}}
-\frac{1}{\rho}\sum_{s=2}^{n}\frac{\partial a_{i_{1}}}{\partial
z_{\alpha}}\frac{\partial \bar{a}_{s}}{\partial
\bar{z}_{\beta}}\tilde{g}^{i_{3}\bar{s}}\tilde{g}_{\lambda\bar{\eta}}
\hspace{0.8cm}\\\hspace{1.5cm}+\frac{1}{\rho^{2}}\sum_{s,t=2}^{n}\Big(\frac{\partial
a_{i_{1}}}{\partial z_{\alpha}}\frac{\partial
\bar{a}_{t}}{\partial
\bar{z}_{\eta}}\tilde{g}^{i_{3}\bar{t}}\Big)\Big(\frac{\partial
a_{i_{2}}}{\partial z_{\lambda}}\frac{\partial
\bar{a}_{s}}{\partial
\bar{z}_{\beta}}\tilde{g}^{i_{2}\bar{s}}\Big)\bigg]dz_{\alpha}\wedge
d\bar{z}_{\beta}\wedge dz_{\lambda}\wedge
d\bar{z}_{\eta}\\=\omega\wedge\Theta_{i_{1}}^{i_{3}}-\frac{1}{\rho}\sum_{\alpha,\beta,s=2}^{n}\frac{\partial
a_{i_{1}}}{\partial z_{\alpha}}\frac{\partial
\bar{a}_{s}}{\partial
\bar{z}_{\beta}}\tilde{g}^{i_{3}\bar{s}}dz_{\alpha}\wedge
d\bar{z}_{\beta}\wedge \omega\hspace{0.8cm}\\-
\big(\frac{\sqrt{-1}}{2\pi})^{2}\sum_{\alpha,\beta,t,\lambda,\eta=2}^{n}\frac{1}{\rho}\frac{\partial
a_{i_{1}}}{\partial z_{\alpha}}\frac{\partial
\bar{a}_{t}}{\partial \bar{z}_{\eta}}\tilde{g}^{i_{3}\bar{t}}\Big(
\partial_{\lambda}\bar{\partial}_{\beta}\log\frac{\sum_{k=0}^{n}|F_{k}|^{2}}{|F_{1}|^{2}}\Big)dz_{\alpha}\wedge
d\bar{z}_{\eta}\wedge dz_{\lambda}\wedge d\bar{z}_{\beta}\\=
\omega\wedge
\Theta_{i_{1}}^{i_{3}}-\frac{\sqrt{-1}}{2\pi}\frac{1}{\rho}\sum_{\alpha,\eta,t=2}^{n}\Big(\frac{\partial
a_{i_{1}}}{\partial z_{\alpha}}\frac{\partial
\bar{a}_{t}}{\partial
\bar{z}_{\eta}}\tilde{g}^{i_{3}\bar{t}}dz_{\alpha}\wedge
d\bar{z}_{\eta}\Big)\wedge(d\omega+\frac{\sqrt{-1}}{2\pi}\partial\bar{\partial}\xi).\label{2.13}\end{multline}
Suppose that the statement is true for $j-1$. By induction
hypothesis, we get
\begin{eqnarray*}
\lefteqn{\sum_{i_{2},\cdots,i_{j-1}=2}^{n}\Theta_{i_{1}}^{i_{2}}\cdots\Theta_{i_{j-1}}^{i_{j}}
=\Big(\omega\wedge\sum_{i_{3},\cdots,i_{j-1}=2}^{n}\Theta_{i_{1}}^{i_{3}}\Theta_{i_{3}}^{i_{4}}\cdots
\Theta_{i_{j-2}}^{i_{j-1}}\Big)\Theta_{i_{j-1}}^{i_{j}}}\\&&
-\left(\frac{\sqrt{-1}}{2\pi}\frac{1}{\rho}\sum_{s,\alpha,\beta=2}^{n}\Big(\frac{\partial
a_{i_{1}}}{\partial z_{\alpha}}\frac{\partial \bar{a}_{s}}{\partial
\bar{z}_{\beta}}\tilde{g}^{i_{j-1}\bar{s}}dz_{\alpha}\wedge
d\bar{z}_{\beta}\Big)\wedge(d\omega+\frac{\sqrt{-1}}{2\pi}\partial\bar{\partial}\xi)^{j-3}\right)\wedge\Theta_{i_{j-1}}^{i_{j}}\\&&=
\omega\wedge\sum_{i_{3},\cdots,i_{j-1}=2}^{n}\Theta_{i_{1}}^{i_{3}}\Theta_{i_{3}}^{i_{4}}\cdots\Theta_{i_{j-1}}^{i_{j}}
\\&&\qquad-\frac{\sqrt{-1}}{2\pi}\frac{1}{\rho}\sum_{\alpha,\beta=2}^{n}\Big(\frac{\partial
a_{i_{1}}}{\partial z_{\alpha}}\frac{\partial \bar{a}_{t}}{\partial
\bar{z}_{\beta}}\tilde{g}^{i_{j}\bar{t}}dz_{\alpha}\wedge
d\bar{z}_{\beta}\Big)\wedge(d\omega+\frac{\sqrt{-1}}{2\pi}\partial\bar{\partial}\xi)^{j-2},
\end{eqnarray*}
where
\begin{multline}
\sum_{i_{j-1},\alpha,\beta,s=2}^{n}\Big(\frac{\partial
a_{i_{1}}}{\partial z_{\alpha}}\frac{\partial
\bar{a}_{s}}{\partial
\bar{z}_{\beta}}\tilde{g}^{i_{j-1}\bar{s}}dz_{\alpha}\wedge
d\bar{z}_{\beta}\Big)\wedge\Theta_{i_{j-1}}^{i_{j}}\\=
\sum_{\alpha,\beta,s=2}^{n}\Big(\frac{\partial a_{i_{1}}}{\partial
z_{\alpha}}\frac{\partial \bar{a}_{s}}{\partial
\bar{z}_{\beta}}\tilde{g}^{i_{j}\bar{s}}dz_{\alpha}\wedge
d\bar{z}_{\beta}\Big)\wedge(d\omega+\frac{\sqrt{-1}}{2\pi}\partial\bar{\partial}\xi)\label{2.14}
\end{multline}
is part of \eqref{2.13}. We omit the detail of the proof.
\end{proof}
\begin{theorem}\label{thm2.1}
With the curvature given in Lemma \ref{lem2.2}, the trace of the
wedge product of $j$ many curvature tensors on the hypersurface $M$
is
\[
\Theta_{i_{1}}^{i_{2}}\cdots\Theta_{i_{j-1}}^{i_{j}}\Theta_{i_{j}}^{i_{1}}
=(n+1)\omega^{j}-(d\omega+\frac{\sqrt{-1}}{2\pi}\partial\bar{\partial}\xi)^{j}\] for $2\leq j\leq q$.
\end{theorem}
\begin{proof}
Prove by induction. For $j=2$, compute it directly from
\eqref{2.13}. That is
\begin{eqnarray*}
\sum_{i_{1},i_{2}=2}^{n}\Theta_{i_{1}}^{i_{2}}\Theta_{i_{2}}^{i_{1}}
&=&\omega\wedge\big((n+1-d)\omega-\frac{\sqrt{-1}}{2\pi}\partial\bar{\partial}\xi\big)
\\&&-\big((d-1)\omega+\frac{\sqrt{-1}}{2\pi}\partial\bar{\partial}\xi\big)\big(d\omega+\frac{\sqrt{-1}}{2\pi}\partial\bar{\partial}\xi\big)\\&=&
(n+1)\omega^{2}-(d\omega+\frac{\sqrt{-1}}{2\pi}\partial\bar{\partial}\xi)^{2}.\end{eqnarray*}
Suppose the statement is true for $j-1$. By Sublemma \ref{sub2.2},
the induction hypothesis and \eqref{2.14}, we get
\begin{eqnarray*}
\lefteqn{\left(\Theta_{i_{1}}^{i_{2}}\Theta_{i_{2}}^{i_{3}}\cdots\Theta_{i_{j-1}}^{i_{j}}\right)\Theta_{i_{j}}^{i_{1}}}\\&&=
\omega\wedge\Theta_{i_{1}}^{i_{3}}\Theta_{i_{3}}^{i_{4}}\cdots\Theta_{i_{j-1}}^{i_{j}}\Theta_{i_{j}}^{i_{1}}\\&&\qquad-
\frac{\sqrt{-1}}{2\pi}\sum_{\alpha,\beta,s=2}^{n}\frac{\partial
a_{i_{1}}}{\partial z_{\alpha}}\frac{\partial
\bar{a}_{s}}{\partial
\bar{z}_{\beta}}\tilde{g}^{i_{j}\bar{s}}dz_{\alpha}\wedge
d\bar{z}_{\beta}\wedge(d\omega+\frac{\sqrt{-1}}{2\pi}\partial\bar{\partial}\xi)^{j-2}\wedge\Theta_{i_{j}}^{i_{1}}\\&&=
\omega\wedge[(n+1)\omega^{j-1}-(d\omega+\frac{\sqrt{-1}}{2\pi}\partial\bar{\partial}\xi)^{j-1}]
\\&&\qquad-[(d-1)\omega+\frac{\sqrt{-1}}{2\pi}\partial\bar{\partial}\xi]
\wedge(d\omega+\frac{\sqrt{-1}}{2\pi}\partial\bar{\partial}\xi)^{j-1}\\&&=
(n+1)\omega^{j}-(d\omega+\frac{\sqrt{-1}}{2\pi}\partial\bar{\partial}\xi)^{j}.
\end{eqnarray*}
\end{proof}
\begin{proof}[Proof of Lemma \ref{lem2.3}]
Prove by induction. For $q=1$,
\begin{equation*}P^{1}(\Theta)=(n+1-d)\omega-\frac{\sqrt{-1}}{2\pi}\partial\bar{\partial}\xi
=\alpha_{10}\frac{\sqrt{-1}}{2\pi}\partial\bar{\partial}\xi+\alpha_{11}\omega,\end{equation*}
where $\alpha_{10}=(-1)$, $\alpha_{11}={n+1\choose
1}-d\alpha_{00}$ and $\alpha_{00}=1$. Suppose that
\[P^{i}(\Theta)=\sum_{k=0}^{i}\alpha_{ik}\omega^{k}\wedge\big(\frac{\sqrt{-1}}{2\pi}\partial\bar{\partial}\xi\big)^{i-k}\]
 holds for
$2\leq i\leq q-1$. By Sublemma \ref{sub2.1}, Theorem \ref{thm2.1}
and induction hypothesis, we have
\begin{multline}P^{q}(\Theta)=\sum_{j=1}^{q}(-1)^{j-1}\phi_{j}P^{q-j}(\Theta)\\
=\sum_{j=1}^{q}(-1)^{j-1}[(n+1)\omega^{j}-(d\omega+\frac{\sqrt{-1}}{2\pi}\partial\bar{\partial}\xi)^{j}]\sum_{k=0}^{q-j}
\alpha_{(q-j)k}\omega^{k}\wedge\big(\frac{\sqrt{-1}}{2\pi}\partial\bar{\partial}\xi\big)^{q-j-k},\label{P}
\end{multline}
where
\begin{eqnarray}
\alpha_{(q-j)k}&=&\left\{%
\begin{array}{ll}
    -(d\alpha_{(q-j-1)(k-1)}+\alpha_{(q-j-1)k}), & \hbox{if $k<q-j$;} \\
    {n+1\choose q-j}-d\alpha_{(q-j-1)(q-j-1)}, & \hbox{$k=q-j$.} \\
\end{array}%
\right. \nonumber\\&=&(-1)^{q-j}\sum_{\ell=0}^{k}{q-j-\ell\choose
k-\ell}(-1)^{\ell}d^{k-\ell}{n+1\choose \ell}\label{2.15}
 \end{eqnarray}
 for $j=1,\cdots,q-1$. We rewrite \eqref{P} as
 $P^{q}(\Theta)=\sum_{k=0}^{q}\alpha_{qk}\omega^{k}\big(\frac{\sqrt{-1}}{2\pi}\partial\bar{\partial}\xi\big)^{q-k}$,
 where
\begin{eqnarray}\alpha_{qk}&=&
\frac{1}{q}\sum_{j=1}^{k}(-1)^{j-1}(n+1-d^{j})\alpha_{(q-j)(k-j)}
\nonumber\\&&+\frac{1}{q}\sum_{i=1}^{q-k}\sum_{j=i}^{k+i}(-1)^{j}d^{j-i}
{j\choose j-i}\alpha_{(q-j)(k-j+i)}.\label{2.151}\end{eqnarray} By
using \eqref{2.15}, we compute \eqref{2.151}
\begin{eqnarray}
\lefteqn{\alpha_{qk}=
\frac{n+1}{q}(-1)^{q-1}\sum_{j=1}^{k}\sum_{\ell=0}^{k-j}{q-j-\ell\choose
k-j-\ell}(-1)^{\ell}d^{k-j-\ell}{n+1\choose \ell}}\label{2.16}\\&&
+\frac{(-1)^{q}}{q}\sum_{i=0}^{q-k}\sum_{j=i}^{k+i}\sum_{\ell=0}^{k-j+i}(-1)^{\ell}d^{k-\ell}
{j\choose j-i}{q-j-\ell\choose k-j+i-\ell}{n+1\choose
\ell}\label{2.18}\\&&
-\frac{(-1)^{q}}{q}\sum_{\ell=0}^{k}{q-\ell\choose k}{n+1\choose
\ell}d^{k-\ell}(-1)^{\ell}\label{2.17}.
\end{eqnarray}
In fact, we have
\begin{equation}
\sum_{\ell=0}^{\lambda-1}(-1)^{\ell}{n+1\choose
\ell}=(-1)^{\lambda-1}{n\choose \lambda
-1}.\label{2.161}\end{equation}
 We change the index in \eqref{2.16} and use \eqref{2.161} to
get
\begin{eqnarray}
\lefteqn{\frac{n+1}{q}(-1)^{q-1}\sum_{j=1}^{k}\sum_{\ell=0}^{k-j}{q-j-\ell\choose
k-j-\ell}(-1)^{\ell}d^{k-j-\ell}{n+1\choose \ell}}\nonumber\\&&
=\frac{n+1}{q}(-1)^{q-1}\sum_{\lambda=1}^{k}{q-\lambda\choose
k-\lambda}d^{k-\lambda}\sum_{\ell=0}^{\lambda-1}(-1)^{\ell}{n+1\choose
\ell}\nonumber\\&&
=\frac{n+1}{q}(-1)^{q-1}\sum_{\lambda=1}^{k}{q-\lambda\choose
k-\lambda}d^{k-\lambda}(-1)^{\lambda-1}{n\choose
\lambda-1}\nonumber\\&&=\frac{1}{q}(-1)^{q}\sum_{\lambda=1}^{k}{q-\lambda\choose
k-\lambda}d^{k-\lambda}(-1)^{\lambda}\lambda{n+1\choose
\lambda}.\label{2.20}
\end{eqnarray}
Change the summing order in \eqref{2.18} and get
\begin{multline}
\frac{(-1)^{q}}{q}\sum_{i=0}^{q-k}\sum_{j=i}^{k+i}\sum_{\ell=0}^{k-j+i}(-1)^{\ell}d^{k-\ell}
{j\choose j-i}{q-j-\ell\choose k-j+i-\ell}{n+1\choose \ell}
\\=\frac{(-1)^{q}}{q}\sum_{j=0}^{k}(-1)^{\ell}d^{k-\ell}
\sum_{i=0}^{q-k}\sum_{\ell=0}^{j}{k+i-j\choose
k-j}{q-k-i+j-\ell\choose j-\ell}{n+1\choose
\ell}\\=\frac{(-1)^{q}}{q}\sum_{\ell=0}^{k}(-1)^{\ell}d^{k-\ell}
\sum_{j=\ell}^{k}\sum_{i=0}^{q-k}{k+i-j\choose
k-j}{q-k-i+j-\ell\choose j-\ell}{n+1\choose \ell}.\label{2.21}
\end{multline}
In fact, we have \begin{equation}\sum_{i=0}^{d}{s+i\choose
s}={s+d+1\choose d+1}.\label{2.211}\end{equation} Re-assemble the
following term:
\begin{multline}
\sum_{i=0}^{q-k}{k+i-j\choose k-j}{q-k-i+j-\ell\choose
j-\ell}=\sum_{i=0}^{q-k}{k+i-j\choose
k-j}\\+\sum_{s=1}^{q-k}\sum_{i=0}^{q-k-s}{k+i-j\choose
k-j}\left[{s+j-\ell\choose j-\ell}-{s-1+j-\ell\choose
j-\ell}\right].\label{2.212}
\end{multline}
Then we may apply \eqref{2.211} to each assorted item in
\eqref{2.212} to get
\begin{multline}
 \sum_{i=0}^{q-k}{k+i-j\choose k-j}{q-k-i+j-\ell\choose
j-\ell}\\
=\sum_{s=0}^{q-k}{k+s-(j-1)\choose
k-(j-1)}{q-k-s+(j-1)-\ell\choose
(j-1)-\ell}.\label{2.23}\end{multline} If we repeat the procedure
in \eqref{2.23} for $r$ times for $r=1,\cdots,j-\ell$, we get
\begin{multline*}\sum_{i=0}^{q-k}{k+i-j\choose
k-j}{q-k-i+j-\ell\choose
j-\ell}\\=\sum_{i=0}^{q-k}{k+i-(j-r)\choose
k-(j-r)}{q-k-i+(j-r)-\ell\choose (j-r)-\ell}.\end{multline*} For
each $j=r+\ell$, we have
\begin{equation*}\sum_{i=0}^{q-k}{k+i-j\choose k-j}{q-k-i+j-\ell\choose
j-\ell}=\sum_{i=0}^{q-k}{k+i-\ell\choose k-\ell}{q-k-i\choose
0}={q-\ell+1\choose k-\ell+1}.\end{equation*} Then we obtain
\begin{equation}\sum_{j=\ell}^{k}\sum_{i=0}^{q-k}{k+i-j\choose
k-j}{q-k-i+j-\ell\choose j-\ell}=(k-\ell+1){q-\ell+1\choose
k-\ell+1}.\label{2.24}\end{equation} Use \eqref{2.24} in
\eqref{2.21} to get
\begin{multline}
\frac{(-1)^{q}}{q}\sum_{i=0}^{q-k}\sum_{j=i}^{k+i}\sum_{\ell=0}^{k-j+i}(-1)^{\ell}d^{k-\ell}
{j\choose j-i}{q-j-\ell\choose k-j+i-\ell}{n+1\choose \ell}\\=
\frac{(-1)^{q}}{q}\sum_{\ell=0}^{k}(-1)^{\ell}d^{k-\ell}(k-\ell+1){q-\ell
+1\choose k-\ell+1}{n+1\choose \ell}\\=
\frac{(-1)^{q}}{q}\sum_{\ell=0}^{k}(-1)^{\ell}d^{k-\ell}(q-\ell+1){q-\ell
\choose k-\ell}{n+1\choose \ell}.\hspace{1.35cm}\label{2.25}
\end{multline}
By adding \eqref{2.20}, \eqref{2.25} and \eqref{2.17}, we can get
the coefficient
\begin{eqnarray}
\alpha_{qk}&=&\frac{(-1)^{q}}{q}\sum_{\ell=1}^{k}{q-\ell\choose
k-\ell}d^{k-\ell}(-1)^{\ell}\ell{n+1\choose \ell}\nonumber\\&&+
\frac{(-1)^{q}}{q}d^{k}q{q\choose k}+
\frac{(-1)^{q}}{q}\sum_{\ell=1}^{k}(-1)^{\ell}d^{k-\ell}(q-\ell+\ell){q-\ell
\choose k-\ell}{n+1\choose
\ell}\nonumber\\&&-\frac{(-1)^{q}}{q}\sum_{\ell=0}^{k}(-1)^{\ell}d^{k-\ell}(q-\ell){q-\ell
\choose k-\ell}{n+1\choose \ell}\nonumber\\&=&
(-1)^{q}\sum_{\ell=0}^{k}(-1)^{\ell}d^{k-\ell}{q-\ell \choose
k-\ell}{n+1\choose \ell}.\label{2.251}
\end{eqnarray}
Comparing \eqref{2.251} and \eqref{2.21}, we get
\begin{equation*}\alpha_{qk}=\left\{%
\begin{aligned}
    {n+1\choose q}-d\alpha_{(q-1)(q-1)},  \qquad\text{if $k=q$;} \\
    -(d\alpha_{(q-1)(k-1)}+\alpha_{(q-1)k}), \qquad \text{if $0\leq k<q$.} \\
\end{aligned}%
\right.
\end{equation*}
\end{proof}
\subsection{Computation of the Bando-Futaki invariants}
The Ricci curvature represents the first Chern class. In Remark
\ref{re2.1}, we have
\[c_{1}(\Theta)=(n+1-d)\omega-\frac{\sqrt{-1}}{2\pi}\partial\bar{\partial}\xi.\]
It is clear that the harmonic part of $c_{1}(\Theta)$ is
$(n+1-d)\omega$. In Lemma \ref{lem2.3}, we compute the $q$-th
Chern form $c_{q}(\Theta)$. In order to obtain the Hodge
decomposition of $c_{q}(\Theta)$ in the de Rham cohomology
$H^{2q}_{DR}(M)$, we need the following Lemma.
\begin{lemma}\label{Lef}
Let $M$ be a hypersurface in $\mathbb{CP}^{n}$. Then we have
\begin{equation*}
\dim H^{2q}_{DR}(M)=1
\end{equation*}
for $q=1,\cdots,n-1$.
\end{lemma}
\begin{proof}
 Consider the Lefschetz Hyperplane
Theorem: Let $\Omega^{q}_{M}$ be the sheaf of germs of holomorphic
$p$-forms on $M$. The map
\begin{equation*}
H^{q}(M,\Omega^{p}_{M})\longrightarrow H^{q}(\mathbb{CP}^{n},
\Omega^{p}_{\mathbb{CP}^{n}})\cong\left\{
                                \begin{array}{ll}
                                  0, & \hbox{if $p\neq q$;} \\
                                  \mathbb{C}, & \hbox{ if $p=q$}
                                \end{array}
                              \right.\end{equation*}  is an isomorphism for $p+q\leq n-2$ and injective for $p+q=n-1$.
                              Therefore we can compute the Hodge number
\begin{eqnarray*}
h^{p,q}(M)=\dim H^{q}(M,\Omega_{M}^{p})=\left\{
             \begin{array}{ll}
               0, & \hbox{ if $p\neq q$ and $p+q\leq n-2$.;}\\
               1, & \hbox{ if $p=q\leq \frac{n-2}{2}$.} \\
               \end{array}
           \right.
\end{eqnarray*}
By Kodaira-Serre Theorem, we have
$h^{p,q}(M)=h^{n-1-p,n-1-q}(M)=0$ if $p\neq q$ and $p+q\geq n$;
$h^{p,p}(M)=1$ for $p\geq \frac{n}{2}$. Then the Hodge
decomposition gives the Betti number
\begin{equation*} b_{r}(M)=\dim H_{DR}^{r}(M)=\sum_{p+q=r}h^{p,q}(M)=\left\{
                          \begin{array}{ll}
                            0, & \hbox{if $r$ is odd;} \\
                            1, & \hbox{if $r\neq n-1$ and $r=2p=2q$.}
                          \end{array}
                        \right.
\end{equation*} Let $\mathcal{H}^{r}=\{\varphi\in
\wedge^{r} T^{\ast}(M)|\triangle \varphi=0\}$ be the vector space
of harmonic $r$-forms on $M$. For $n=2p+1$, consider the map
\[L:\mathcal{H}^{n-1}(M)\longrightarrow
\mathcal{H}^{n+1}(M),\] where $L(\phi)=\omega\wedge\phi$ for
$\phi\in \mathcal{H}^{n-1}(M)$. The map $L$ is well-defined since
we have $[L,\triangle]=0$ on compact K\"ahler manifolds, where
$\triangle=d\delta+\delta d$. Therefore,
$L(\phi)=\omega\wedge\phi$ is also harmonic. By Hodge Theorem and
Lefschetz Hyperplane Theorem, we have $\dim
\mathcal{H}^{n-1}(M)=\dim H^{n-1}_{DR}(M)\geq 1$. Suppose that
there exists a harmonic $(n-1)$-form $\phi\neq
c\omega^{\frac{n-1}{2}}$ for all $c\in\mathbb{C}$. Since
$\omega\wedge \phi\neq c\omega^{\frac{n+1}{2}}$, this implies
$\dim \mathcal{H}^{n+1}(M)\geq 2$. It contradicts to the fact
$\dim \mathcal{H}^{n+1}(M)=1$.
\end{proof}
Let  For $1\leq q\leq n-1$, $\omega^{q}$ is harmonic on compact
K\"ahler manifolds. The Hodge Theorem says that for any $(q,q)$
form $\phi$ and $[\phi]\in H_{DR}^{2q}(M)$, there exists a unique
harmonic form representing $[\phi]$. Since $\dim
\mathcal{H}^{2q}(M)=\dim H_{DR}^{2q}(M)=1$, there exists some
constant $c$ such that $[\phi]=[c\omega^{q}]\in H^{2q}_{DR}(M)$.
That is,
\begin{corollary}\label{cor2.1}
For $1\leq q\leq n-1$, the harmonic part of $c_{q}(\Theta)$ is
proportional to $\omega^{q}$. \end{corollary} By Corollary
\ref{cor2.1} and the Hodge decomposition theory, we have that
\begin{equation}c_{q}(\Theta)=\alpha_{qq}\omega^{q}
+\frac{\sqrt{-1}}{2\pi}\partial\bar{\partial}f_{q},\label{2.27}\end{equation}
where
$f_{q}=\sum_{k=0}^{q-1}\alpha_{qk}\xi\omega^{q}\wedge\big(\frac{\sqrt{-1}}{2\pi}\partial\bar{\partial}\xi\big)^{q-k-1}$.
 However, we will not compute the $q$-th Bando-Futaki invariant directly:
\[\mathcal{F}_{q}(X)=\int_{M}\mathcal{L}_{X}f_{q}\wedge\omega_{M}^{n-q-1},\]
where $\omega_{M}=(n+1-d)\omega\in c_{1}(M)$ and
$\omega=\omega_{FS}|_{M}$. \vspace{0.3cm} First, take the
contraction map on \eqref{2.27} with $X$:
\begin{equation}
i(X)c_{q}(\Theta)-i(X)\alpha_{qq}\omega_{q}=\frac{\sqrt{-1}}{2\pi}\bar{\partial}i(X)\partial
f_{q}.\label{13}\end{equation} Let us compute each term in
\eqref{13} separately. Take the contraction of $\omega^{q}$ with
$X$:
\[i(X)\omega^{q}=qi(X)(\omega)\omega^{q-1}=q(-\frac{\sqrt{-1}}{2\pi}\bar{\partial}\theta)\omega^{q-1},\]
where $i(X)\omega=-\frac{\sqrt{-1}}{2\pi}\bar{\partial}\theta$.
More precisely, we can express a holomorphic vector field
\begin{equation}\tilde{X}=\sum_{i=0}^{n}\lambda_{i}Z_{i}\frac{\partial}{\partial
Z_{i}}\label{2.301}\end{equation} over $\mathbb{CP}^{n}$ with
$\sum_{k=0}^{n}\lambda_{k}=0$. The
 restriction of $\tilde{X}$ on $M\cap V$ is given by
\[\tilde{X}|_{V}=X=\sum_{i=2}^{n}(\lambda_{i}-\lambda_{0})z_{i}\left(a_{i}\frac{\partial}{\partial
z_{1}}+\frac{\partial}{\partial z_{i}}\right).\] Hence the Hamilton
function can be expressed explicitly
\begin{equation}\theta=-\tilde{X}\log(\sum_{k=0}^{n}|Z_{k}|^{2})
=-\frac{\sum_{k=0}^{n}\lambda_{k}|Z_{k}|^{2}}{\sum_{k=0}^{n}|Z_{k}|^{k}}
=-\frac{\sum_{k=1}^{n}\lambda_{k}|z_{k}|^{2}}{1+|z|^{2}}-\lambda_{0}.\label{15}\end{equation}
Let $\tilde{P}\underbrace{(\Theta,\cdots,\Theta)}_{q}$ be the
polarization of its elementary invariant polynomial
$P^{q}(\Theta)$. Take the contraction of the curvature with $X$:
\[i(X)\Theta_{k}^{\ell}=i(X)\frac{\sqrt{-1}}{2\pi}\sum_{i,j=2}^{n}R^{\ell}_{ki\bar{j}}dz_{i}\wedge
d\bar{z}_{j}=\frac{\sqrt{-1}}{2\pi}
\sum_{i,j=2}^{n}X^{i}R^{\ell}_{ki\bar{j}}d\bar{z}_{j}
=-\frac{\sqrt{-1}}{2\pi}\bar{\partial}X^{\ell}_{k},\] where
\begin{equation}X^{\ell}_{k}=\frac{\partial X^{\ell}}{\partial
z_{k}}+\sum_{i}X^{i}\Gamma^{\ell}_{ik}=-\sum_{j=2}^{n}\tilde{g}^{\ell
\bar{j}}\partial_{k}\bar{\partial}_{j}\theta\label{2.311}\end{equation}
for $2\leq k,\ell\leq n$. Let \[\nabla X=\sum_{k,\ell}
X^{\ell}_{k} dz_{k}\otimes \frac{\partial}{\partial
z_{\ell}}=\sum_{k,\ell}(\frac{\partial X^{\ell}}{\partial
z_{k}}+\sum_{i}X^{i}\Gamma^{\ell}_{ki})dz_{k}\otimes
\frac{\partial}{\partial z_{\ell}}.\] Then we can take the
contraction map on the $q$-th Chern form
\[i(X)c_{q}(\Theta)=q
\tilde{P}^{q}(i(X)\Theta,\Theta,\cdots,\Theta)=
-q\frac{\sqrt{-1}}{2\pi}\bar{\partial} \tilde{P}^{q}(\nabla
X,\Theta,\cdots,\Theta).\] We may re-write \eqref{13} as
\[\bar{\partial}[-q\tilde{P}^{q}(\nabla
X,\Theta,\cdots,\Theta)+q\alpha_{qq}\theta
\omega^{q-1}-i(X)\partial f_{q,\Theta}]=0.\] By Hodge
Decomposition Theorem, we get
\begin{equation}-q\tilde{P}^{q}(\nabla
X,\Theta,\cdots,\Theta)+q\alpha_{qq}\theta
\omega^{q-1}-i(X)\partial f_{q,\Theta}=\psi_{q}+\bar{\partial}
\varphi_{q},\label{2.30}\end{equation} where $\psi_{q}$ is the
harmonic part of the left-hand side of \eqref{2.30} and
$\varphi_{q}$ is a $2(q-1)-1$ form. Since the right hand side is
of $(q-1,q-1)$ form, $\varphi_{q}$ is of $(q-1,q-2)$ form. By
Lemma \ref{Lef}, there exists a constant $C(q)$ such that
$\psi_{q}=C(q)\omega^{q-1}$. Instead of computing the Bando-Futaki
invariants directly, we will compute the following:
\begin{multline} \int_{M}L_{X}f_{q}\wedge
\omega_{M}^{n-q}=\int_{M}(di(X)f_{q}+i(X)\partial f_{q})\wedge
\omega_{M}^{n-q}\\\qquad=\int_{M}\left(-q\tilde{P}^{q}(\nabla
X,\Theta,\cdots,\Theta)+q\alpha_{qq}\theta\omega^{q-1}-C(q)\omega^{q-1}-\bar{\partial}\varphi_{q}\right)\wedge
\omega_{M}^{n-q}.\label{2.31}\end{multline} Finding $C(q)$ and
showing $q\int_{M}\tilde{P}^{q}(\nabla
X,\Theta,\cdots,\Theta)\wedge \omega_{M}^{n-q}=0$ are the next two
steps to compute \eqref{2.31}. In order to evaluate $C(q)$, it is
necessary to express $\tilde{P}^{q}(\nabla
X,\Theta,\cdots,\Theta)$ explicitly.
\begin{lemma}\label{lem2.4}
The covariant derivative of the polarization of the elementary
polynomial $P^{q}$ is given by
\begin{multline}
q\tilde{P}^{q}(\nabla X,\Theta,\cdots,\Theta)=-{\rm
div}(X)\gamma_{q1}+\theta
\gamma_{q2}\\-\sum_{j=2}^{q}(-1)^{j-1}\eta_{q-j}
\wedge\Big(\big(\omega^{j-2}+
n\zeta_{j-2}\big)\wedge\big(\partial\bar{\partial}\theta)-\zeta_{j-2}\wedge
\big(\partial\bar{\partial}\Delta\theta\big)\Big)
,\label{2.32}\end{multline} where
\[\begin{array}{lcl}
\gamma_{q1}&=&\sum_{k=0}^{q-1}\alpha_{qk}(q-k)\omega^{k}\wedge(\frac{\sqrt{-1}}{2\pi}\partial\bar{\partial}\xi)^{q-1-k},\\
\gamma_{q2}&=&\sum_{k=0}^{q-1}\left((q-k)(n+1-d)\alpha_{qk}+(k+1)\alpha_{q(k+1)}\right)
\omega^{k}\wedge(\frac{\sqrt{-1}}{2\pi}\partial\bar{\partial}\xi)^{q-1-k},\\
\eta_{q-j}&=&\sum_{k=0}^{q-j}\alpha_{(q-j)k}\frac{\sqrt{-1}}{2\pi}
\omega^{k}\wedge(\frac{\sqrt{-1}}{2\pi}\partial\bar{\partial}\xi)^{q-j-k},\\
\zeta_{j-2}&=&\sum_{k=0}^{j-2}(d\omega+\frac{\sqrt{-1}}{2\pi}\partial\bar{\partial}\xi)^{j-2-k}\wedge\omega^{k}, \\
\Delta\theta&=&-\sum_{\alpha,\beta=2}^{n}\tilde{g}^{\alpha\bar{\beta}}\partial_{\alpha}\bar{\partial}_{\beta}\theta.
\end{array}\]
\end{lemma}\begin{proof} According to Sublemma \ref{sub2.1}, we have
\begin{multline}q\tilde{P}^{q}(\nabla
X,\Theta,\cdots,\Theta)=\frac{q}{q!}\sum_{\sigma\in
S_{q}}\sum_{i_{1},\cdots,i_{q}=2}^{n}{\rm
sgn}(\sigma)X_{i_{1}}^{i_{\sigma(1)}}\Theta_{i_{2}}^{i_{\sigma(2)}}\cdots\Theta_{i_{2}}^{i_{\sigma(2)}}
\\=\sum_{i_{1}=2}^{n}X_{i_{1}}^{i_{1}}P^{q-1}(\Theta)-\sum_{i_{1},i_{2}}^{n}X_{i_{1}}^{i_{2}}\Theta_{i_{2}}^{i_{1}}P^{q-2}(\Theta)
+\sum_{j=3}^{q}(-1)^{j-1}E_{j}P^{q-j}(\Theta),\label{2.33}\end{multline}
where
\[E_{j}=\sum_{i_{1},\cdots,i_{j}=2}^{n}X_{i_{1}}^{i_{2}}\Theta_{i_{2}}^{i_{3}}\cdots\Theta_{i_{j}}^{i_{1}}\]
for $j=3,\cdots,q$. Let $E_{1}=\sum_{i=2}^{n}X_{i}^{i}$ and
$E_{2}=\sum_{i,j=2}^{n}X_{i}^{j}\Theta_{i}^{j}$. To formularize
$E_{j}$, we need the following:
\begin{sublem}\label{sub2.3}
For $j=3,\cdots,q$, we have a regression relation for $E_{j}$:
\begin{eqnarray}
E_{j}&=&\omega\wedge E_{j-1}+\Phi
\wedge(d\omega+\frac{\sqrt{-1}}{2\pi}\partial\bar{\partial}\xi)^{j-2}\label{2.34}\\&=&
\omega^{j-2}\wedge
E_{2}+\Phi\wedge\sum_{k=0}^{j-3}(d\omega+\frac{\sqrt{-1}}{2\pi}\partial\bar{\partial}\xi)^{j-2-k}\omega^{k},\label{2.35}
\end{eqnarray}
where
\begin{equation}\Phi=-\frac{\sqrt{-1}}{2\pi}\frac{1}{\rho}\sum_{k,\ell,p,q=2}^{n}X_{k}^{\ell}\frac{\partial
a_{\ell}}{\partial z_{p}}\frac{\partial \bar{a}_{s}}{\partial
\bar{z}_{q}}\tilde{g}^{k\bar{s}}dz_{p}\wedge
d\bar{z}_{q}.\label{2.36}\end{equation}
\end{sublem}
\begin{proof} \eqref{2.34} can be proved by
induction. For $j=3$, use Sublemma \ref{sub2.2} to compute
\begin{eqnarray*}
\lefteqn{E_{3}=\sum_{i,j,k=2}^{3}X_{i}^{j}\Theta_{j}^{k}\Theta_{k}^{i}=
 \omega\wedge \sum_{i,j=2}^{n}(X_{i}^{j}\Theta_{j}^{i})}\\&&\qquad-
\frac{\sqrt{-1}}{2\pi}\frac{1}{\rho}\sum_{\alpha,\beta=2}^{n}X_{i}^{j}\left(\frac{\partial
a_{j}}{\partial z_{\alpha}}\frac{\partial \bar{a}_{s}}{\partial
\bar{z}_{\beta}}\tilde{g}^{i\bar{s}}dz_{\alpha}\wedge
d\bar{z}_{\beta}\right)\wedge(d\omega+\frac{\sqrt{-1}}{2\pi}\partial\bar{\partial}\xi)\\&&=
\omega\wedge E_{2}+\Phi\wedge
(d\omega+\frac{\sqrt{-1}}{2\pi}\partial\bar{\partial}\xi).\end{eqnarray*}
 Assume that \eqref{2.34} holds for $j-1$. By Sublemma \ref{sub2.2}, the induction
 hypothesis, and \eqref{2.14}, compute that
\begin{multline*}
E_{j}=\sum_{i_{1},\cdots,i_{j}=2}^{n}X_{i_{1}}^{i_{2}}\left(\Theta_{i_{2}}^{i_{3}}\cdots\Theta_{i_{j-1}}^{i_{j}}\right)\Theta_{i_{j}}^{i_{1}}\\=
\omega\sum_{i_{1},i_{2},i_{4}\cdots,i_{j}=2}^{n}X_{i_{1}}^{i_{2}}
\Theta_{i_{2}}^{i_{4}}\Theta_{i_{4}}^{i_{5}}\cdots\Theta_{i_{j-1}}^{i_{j}}\Theta_{i_{j}}^{i_{1}}\hspace{5cm}\\
-\frac{\sqrt{-1}}{2\pi}\frac{1}{\rho}\sum_{i_{1},i_{2},i_{j},\alpha,\beta=2}^{n}X_{i_{1}}^{i_{2}}\left(\frac{\partial
a_{i_{2}}}{\partial z_{\alpha}}\frac{\partial
\bar{a}_{s}}{\partial
\bar{z}_{\beta}}\tilde{g}^{i_{j}\bar{s}}dz_{\alpha}\wedge
d\bar{z}_{\beta}\right)\wedge(d\omega+\frac{\sqrt{-1}}{2\pi}\partial\bar{\partial}\xi)^{j-3}\wedge
\Theta_{i_{j}}^{i_{1}}\\
=\omega\wedge
E_{j-1}-\frac{\sqrt{-1}}{2\pi}\frac{1}{\rho}\sum_{i_{1},i_{2},i_{j},\alpha,\beta=2}^{n}
X_{i_{1}}^{i_{2}}\left(\frac{\partial a_{i_{2}}}{\partial
z_{\alpha}}\frac{\partial \bar{a}_{s}}{\partial
\bar{z}_{\beta}}\tilde{g}^{i_{1}\bar{s}}dz_{\alpha}\wedge
d\bar{z}_{\beta}\right)\wedge(d\omega+\frac{\sqrt{-1}}{2\pi}\partial\bar{\partial}\xi)^{j-2}\\
=\omega\wedge
E_{j-1}+\Phi\wedge(d\omega+\frac{\sqrt{-1}}{2\pi}\partial\bar{\partial}\xi)^{j-2}.\hspace{6.6cm}
\end{multline*}
Equation \eqref{2.35} follows directly from \eqref{2.34}.
\end{proof}Before we start compute $E_{j}$, we need to find
$\Phi$.
\begin{sublem}\label{sub2.4}
We compute $\Phi$ explicitly as follows:
\begin{eqnarray*}
\Phi&=&{\rm
div}(X)\left((d-1)w+\frac{\sqrt{-1}}{2\pi}\partial\bar{\partial}\xi\right)
-n\frac{\sqrt{-1}}{2\pi}\partial\bar{\partial}\theta+\frac{\sqrt{-1}}{2\pi}\partial\bar{\partial}\Delta\theta
\\&&-(n+1)\theta\left((d-1)w+\frac{\sqrt{-1}}{2\pi}\partial\bar{\partial}\xi\right).
\end{eqnarray*}
\end{sublem}
\begin{proof}
Using \eqref{2.311} in \eqref{2.36}, $\Phi$ can be computed
alternately by
\[\Phi=\frac{\sqrt{-1}}{2\pi}\frac{1}{\rho}\sum_{k,\ell,j,p,q=2}^{n}
\tilde{g}^{\ell\bar{j}}\partial_{k}\bar{\partial}_{j}\theta\frac{\partial
a_{\ell}}{\partial z_{p}}\frac{\partial \bar{a}_{s}}{\partial
\bar{z}_{q}}\tilde{g}^{k\bar{s}}dz_{p}\wedge d\bar{z}_{q}.\] By
using the definition of $\theta$ in \eqref{15}, we get
\[\bar{\partial}_{j}\theta=
-\frac{(\lambda_{j}-\lambda_{0})z_{j}}{1+|z|^{2}}-
\frac{(\lambda_{1}-\lambda_{0})z_{1}\bar{a}_{j}}{1+|z|^{2}}
-(\theta+\lambda_{0})\frac{z_{j}+z_{1}\bar{a}_{j}}{1+|z|^{2}}.\]
Then we compute
\begin{eqnarray}
\partial_{k}\bar{\partial}_{j}\theta &=&
-\frac{\delta_{kj}(\lambda_{j}-\lambda_{0})}{1+|z|^{2}}
+\frac{(\lambda_{j}-\lambda_{0})z_{j}(\bar{z}_{k}+\bar{z}_{1}a_{k})}{(1+|z|^{2})^{2}}\nonumber\\&&
-\frac{(\lambda_{1}-\lambda_{0})a_{k}\bar{a}_{j}}{1+|z|^{2}}
+\frac{(\lambda_{1}-\lambda_{0})z_{1}\bar{a}_{j}(\bar{z}_{k}+\bar{z}_{1}a_{k})}{(1+|z|^{2})^{2}}\nonumber\\&&
+\frac{(\lambda_{j}-\lambda_{0})\bar{z}_{k}}{1+|z|^{2}}\frac{z_{j}+z_{1}\bar{a}_{j}}{1+|z|^{2}}
+\frac{(\lambda_{1}-\lambda_{0})\bar{z}_{1}z_{k}}{1+|z|^{2}}\frac{z_{j}+z_{1}\bar{a}_{j}}{1+|z|^{2}}\nonumber\\&&
+2(\theta+\lambda_{0})\frac{(\bar{z}_{k}+\bar{z}_{1}a_{k})(z_{j}+z_{1}\bar{a}_{j})}{(1+|z|^{2})^{2}}
-(\theta+\lambda_{0})\frac{\delta_{kj}+a_{k}\bar{a}_{j}}{1+|z|^{2}}.\label{2.39}
\end{eqnarray}
Then we obtain that
\begin{eqnarray*}
\tilde{g}^{\ell\bar{j}}\partial_{k}\bar{\partial}_{j}\theta &=&
-\delta_{\ell
k}(\lambda_{\ell}-\lambda_{0})+\frac{(\lambda_{\ell}-\lambda_{0})z_{\ell}(\bar{z}_{k}+\bar{z}_{1}a_{k})}{1+|z|^{2}}
-\delta_{\ell
k}(\theta+\lambda_{0})\\&&+\bar{a}_{\ell}\frac{|F_{1}|^{2}}{\sum_{i=0}^{n}|
F_{i}|^{2}}a_{k}(\lambda_{k}-\lambda_{0})
-(\lambda_{1}-\lambda_{0})a_{k}\bar{a}_{\ell}\frac{|F_{1}|^{2}}{\sum_{i=0}^{n}|
F_{i}|^{2}}\\&& -z_{\ell}\frac{|F_{1}|^{2}}{\sum_{i=0}^{n}|
F_{i}|^{2}}\frac{\bar{F}_{0}}{\bar{F}_{1}}a_{k}\big((\lambda_{1}-\lambda_{0})-(\lambda_{k}-\lambda_{0})\big).
\end{eqnarray*}
Hence the Laplace of $\theta$ is
\begin{eqnarray}
\triangle\theta&=&-\tilde{g}^{k\bar{j}}\partial_{k}\bar{\partial}_{j}\theta=
-\left(\lambda_{0}+\frac{|F_{1}|^{2}}{\sum_{i=0}^{2}|
F_{i}|^{2}}\sum_{k=1}^{n}(\lambda_{k}-\lambda_{0})|a_{k}|^{2}-n\theta\right)\nonumber\\&=&
-\left(\frac{\sum_{k=0}^{n}\lambda_{k}|F_{k}|^{2}}{\sum_{i=0}^{n}|
F_{i}|^{2}}-n\theta\right).\label{2.41}
\end{eqnarray}
Denote that $\rm{div}(X)=\sum_{i=2}^{n}X_{i}^{i}=\triangle\theta$.
Observe that
\begin{eqnarray*}
\partial_{p}\bar{\partial}_{q}\triangle\theta&=&n\partial_{p}\bar{\partial}_{q}\theta-
\frac{|F_{1}|^{2}}{\sum_{i=0}^{n}|
F_{i}|^{2}}\sum_{k=2}^{n}(\lambda_{k}-\lambda_{0})\frac{\partial
a_{k}}{\partial
z_{p}}\frac{\partial\bar{a}_{k}}{\partial\bar{z}_{q}}\\&&
+\frac{|F_{1}|^{2}}{\sum_{i=0}^{n}|
F_{i}|^{2}}\sum_{k=2}^{n}(\lambda_{k}-\lambda_{0})\frac{\partial
a_{k}}{\partial
z_{p}}\bar{a}_{k}\left(\sum_{s=2}^{n}\frac{\partial\bar{a}_{s}}{\partial\bar{z}_{q}}a_{s}+
\frac{F_{0}}{F_{1}}\sum_{s=2}^{n}\frac{\partial\bar{a}_{s}}{\partial\bar{z}_{q}}\bar{z}_{s}\right)\\&&
+\frac{|F_{1}|^{2}}{\sum_{i=0}^{n}|
F_{i}|^{2}}\sum_{k=2}^{n}(\lambda_{k}-\lambda_{0})\frac{\partial
\bar{a}_{k}}{\partial
\bar{z}_{q}}a_{k}\left(\sum_{s=2}^{n}\frac{\partial
a_{s}}{\partial z_{p}}\bar{a}_{s}+
\frac{\bar{F}_{0}}{\bar{F}_{1}}\sum_{s=2}^{n}\frac{\partial
a_{s}}{\partial z_{p}}z_{s}\right)\\&&
-\sum_{k=1}^{n}(\lambda_{k}-\lambda_{0})|a_{k}|^{2}(\frac{|F_{1}|^{2}}{\sum_{i=0}^{n}|
F_{i}|^{2}})^{3}\\&&\quad\cdot\left(\sum_{s=2}^{n}\frac{\partial
a_{s}}{\partial z_{p}}\bar{a}_{s}+
\frac{\bar{F}_{0}}{\bar{F}_{1}}\sum_{s=2}^{n}\frac{\partial
a_{s}}{\partial z_{p}}z_{s}\right)
\left(\sum_{s=2}^{n}\frac{\partial\bar{a}_{s}}{\partial\bar{z}_{q}}a_{s}+
\frac{F_{0}}{F_{1}}\sum_{s=2}^{n}\frac{\partial\bar{a}_{s}}{\partial\bar{z}_{q}}\bar{z}_{s}\right)\\&&
+\sum_{k=1}^{n}(\lambda_{k}-\lambda_{0})|a_{k}|^{2}\frac{|F_{1}|^{2}}{\sum_{i=0}^{n}|
F_{i}|^{2}}\partial_{p}\bar{\partial}_{q}\log|\frac{\nabla
F}{F_{1}}|^{2}.
\end{eqnarray*}
By using the definition of $\tilde{g}^{k\bar{s}}$ given in Lemma
\ref{lem2.1}, we get
\begin{multline}
\rho\frac{\partial a_{\ell}}{\partial z_{p}}\frac{\partial
\bar{a}_{s}}{\partial\bar{z}_{q}}\tilde{g}^{k\bar{s}}=
\frac{\sum_{i=0}^{n}| F_{i}|^{2}}{|F_{1}|^{2}}\frac{\partial
a_{\ell}}{\partial
z_{p}}\frac{\partial\bar{a}_{k}}{\partial\bar{z}_{q}}-
\frac{\partial a_{\ell}}{\partial
z_{p}}\bar{a_{k}}\sum_{s=2}^{n}\frac{\partial\bar{a}_{s}}{\partial\bar{z}_{q}}a_{s}
+(1+|a|^{2})\frac{\partial a_{\ell}}{\partial
z_{p}}z_{k}\sum_{s=2}^{n}\frac{\partial\bar{a}_{s}}{\partial\bar{z}_{q}}\bar{z}_{s}
\\-\frac{\bar{F_{0}}}{\bar{F_{1}}}\frac{\partial a_{\ell}}{\partial
z_{p}}z_{k}\sum_{s=2}^{n}\frac{\partial\bar{a}_{s}}{\partial\bar{z}_{q}}a_{s}
-\frac{F_{0}}{F_{1}}\frac{\partial a_{\ell}}{\partial
z_{p}}\bar{a_{k}}\sum_{s=2}^{n}\frac{\partial\bar{a}_{s}}{\partial\bar{z}_{q}}\bar{z}_{s}.\hspace{2cm}\label{2.40}
\end{multline}
 Then multiple \eqref{2.39} and \eqref{2.40} together to get
\begin{eqnarray*}
\frac{1}{\rho}\tilde{g}^{\ell\bar{j}}\partial_{k}\bar{\partial}_{j}\theta
\frac{\partial a_{\ell}}{\partial z_{p}}\frac{\partial
\bar{a}_{s}}{\partial\bar{z}_{q}}\tilde{g}^{k\bar{s}}&=&
\partial_{p}\bar{\partial}_{q}\triangle\theta-n\partial_{p}\bar{\partial}_{q}\theta-
(\theta+\lambda_{0})\partial_{p}\bar{\partial}_{q}\log\frac{\sum_{i=0}^{n}|
F_{i}|^{2}}{|F_{1}|^{2}}\\&&-\sum_{i=1}^{n}(\lambda_{i}-\lambda_{0})|a_{i}|^{2}\frac{|F_{1}|^{2}}{\sum_{i=0}^{n}|
F_{i}|^{2}}|^{2}\partial_{p}\bar{\partial}_{q}\log\frac{\sum_{i=0}^{n}|
F_{i}|^{2}}{|F_{1}|^{2}}.
\end{eqnarray*}
Therefore, we obtain that \begin{eqnarray*}
\Phi&=&\frac{\sqrt{-1}}{2\pi}\partial\bar{\partial}\triangle\theta-
\frac{\sqrt{-1}}{2\pi}n\partial\bar{\partial}\theta
-(\theta+\lambda_{0})\partial\bar{\partial}\log\frac{\sum_{i=0}^{2}|
F_{i}|^{2}}{|F_{1}|^{2}}\\&&-\sum_{i=1}^{n}(\lambda_{i}-\lambda_{0})|a_{i}|^{2}\frac{|F_{1}|^{2}}{\sum_{i=0}^{2}|
F_{i}|^{2}}\partial\bar{\partial}\frac{\sum_{i=0}^{2}|
F_{i}|^{2}}{|F_{1}|^{2}}\\&=&
\frac{\sqrt{-1}}{2\pi}\partial\bar{\partial}\triangle\theta-
\frac{\sqrt{-1}}{2\pi}n\partial\bar{\partial}\theta
-(n+1)\theta\Big((d-1)\omega+\frac{\sqrt{-1}}{2\pi}\partial\bar{\partial}\xi\Big)
\\&&+\Big(n\theta-\sum_{i=1}^{n}(\lambda_{i}-\lambda_{0})|a_{i}|^{2}\frac{|F_{1}|^{2}}{\sum_{i=0}^{2}|
F_{i}|^{2}}-\lambda_{0}\Big)\Big((d-1)\omega+\frac{\sqrt{-1}}{2\pi}\partial\bar{\partial}\xi\Big).
\end{eqnarray*}
Then by using \eqref{2.41} and $\Delta\theta={\rm div}(X)$, we get
$\Phi$.
\end{proof}
\begin{sublem}\label{sub2.5}
For $j=3,\cdots,q$, formularize $E_{j}$ explicitly as following:
\begin{eqnarray}
E_{j}&=& {\rm div}(X)\psi_{j-1}
-(n+1)\theta(\psi_{j-1}-\omega^{j-1}) \nonumber\\&&
-\big(\frac{\sqrt{-1}}{2\pi}\partial\bar{\partial}\theta\big)\wedge(n\zeta_{j-2}+\omega^{j-2})
+\big(\frac{\sqrt{-1}}{2\pi}\partial\bar{\partial}\Delta\theta\big)\wedge\zeta_{j-2},\label{21}
\end{eqnarray}
where
$\zeta_{i}=\sum_{k=0}^{i}(d\omega+\frac{\sqrt{-1}}{2\pi}\partial\bar{\partial}\xi)^{i-k}\omega^{k}$
and
$\psi_{i}=(d\omega+\frac{\sqrt{-1}}{2\pi}\partial\bar{\partial}\xi)^{i}$.
\end{sublem}
\begin{proof}
$E_{2}$ can be obtained directly from $\Phi$:
\[E_{2}=\sum_{i,j=2}^{n}X_{i}^{j}\Theta_{j}^{i}=
(\rm{div}(X)\omega-\frac{\sqrt{-1}}{2\pi}\partial\bar{\partial}\theta
+\Phi).\]
 Then by using \eqref{2.35} and Sublemma \ref{sub2.4}, we get
\begin{eqnarray*}
E_{j}&=&E_{2}\wedge\omega^{j-2}
+\Phi\wedge\sum_{k=0}^{j-3}(d\omega+\frac{\sqrt{-1}}{2\pi}\partial\bar{\partial}\xi)^{j-2-k}\omega^{k}\\
&=&\Phi\wedge\sum_{k=0}^{j-2}(d\omega+\frac{\sqrt{-1}}{2\pi}\partial\bar{\partial}\xi)^{j-2-k}\omega^{k}
+\rm{div}(X)\omega^{j-1}-\frac{\sqrt{-1}}{2\pi}\partial\bar{\partial}\theta\wedge\omega^{j-2}\\&=&
\rm{div}(X)(d\omega+\frac{\sqrt{-1}}{2\pi}\partial\bar{\partial}\xi)^{j-1}
-(n+1)\theta\left((d\omega+\frac{\sqrt{-1}}{2\pi}\partial\bar{\partial}\xi)^{j-1}-\omega^{j-1}\right)\\&&
+\frac{\sqrt{-1}}{2\pi}\partial\bar{\partial}\theta
\wedge\left(n\sum_{k=0}^{j-2}(d\omega+\frac{\sqrt{-1}}{2\pi}\partial\bar{\partial}\xi)^{j-2-k}\omega^{k}-\omega^{j-2}\right)\\&&
+\frac{\sqrt{-1}}{2\pi}\partial\bar{\partial}\triangle\theta
\wedge\sum_{k=0}^{j-2}(d\omega+\frac{\sqrt{-1}}{2\pi}\partial\bar{\partial}\xi)^{j-2-k}\omega^{k}.
\end{eqnarray*}
\end{proof}
Back to the proof of Lemma 2.4. Denote that
$E_{1}=\sum_{i=2}^{n}X_{i}^{i}=\rm{div}X=\triangle\theta$. In
fact, we can write down $E_{2}$ explicitly:
\begin{eqnarray*}
E_{2}&=&{\rm
div}X\left(d\omega+\frac{\sqrt{-1}}{2\pi}\partial\bar{\partial}\xi\right)-(n+1)\frac{\sqrt{-1}}{2\pi}\partial\bar{\partial}\theta
\\&&+\frac{\sqrt{-1}}{2\pi}\partial\bar{\partial}\triangle\theta-(n+1)\theta
\left((d-1)\omega+\frac{\sqrt{-1}}{2\pi}\partial\bar{\partial}\xi\right).\end{eqnarray*}
Then we can use Sublemma \ref{sub2.5} to compute
$\sum_{j=1}^{q}(-1)^{j-1}E_{j}P^{q-j}(\Theta)$. The coefficient
for $\rm{div}(X)$ in \eqref{2.33} is given by
\begin{multline}
\sum_{j=1}^{q}(-1)^{j-1}\psi_{j-1}\sum_{k=0}^{q-j}\alpha_{(q-j)k}\omega^{k}
\big(\frac{\sqrt{-1}}{2\pi}\partial\bar{\partial}\xi\big)^{q-j-k}\\
=\sum_{j=1}^{q}(-1)^{j-1}(d\omega+\frac{\sqrt{-1}}{2\pi}\partial\bar{\partial}\xi)^{j-1}\sum_{k=0}^{q-j}\alpha_{(q-j)k}\omega^{k}
\big(\frac{\sqrt{-1}}{2\pi}\partial\bar{\partial}\xi\big)^{q-j-k}\hspace{1cm}\\=
\sum_{t=1}^{q}\sum_{\ell=0}^{t-1}\sum_{j=t-\ell}^{q-\ell}\alpha_{(q-j)(q-j-\ell)}(-1)^{j-1}d^{j+\ell-t}{j-1\choose
j+\ell-t}\omega^{q-t}\big(\frac{\sqrt{-1}}{2\pi}\partial\bar{\partial}\xi\big)^{t-1}.\label{2.42}\end{multline}
We derive several formulas before further computation. For each
fixed index $j\geq t$, use \eqref{2.15} to get
\begin{multline}
\sum_{\ell=0}^{t-1}\alpha_{(q-j)(q-j-\ell)}(-1)^{j-1}d^{j+\ell-t}{j-1\choose
j+\ell-t}\\=\sum_{\ell=1}^{t-1}(-1)\Big(d\alpha_{(q-j-1)(q-j-\ell-1)}+\alpha_{(q-j-1)(q-j-\ell)}\Big)
(-1)^{j-1}d^{j+\ell-t}{j-1\choose
j+\ell-t}\\+(-1)^{j-1}d^{j-t}\Big({n+1\choose
q-j}-d\alpha_{(q-j-1)(q-j-1)}\Big){j-1\choose
j-t}.\hspace{2cm}\label{2.43}
\end{multline}
 Rewrite \eqref{2.43} to get
\begin{multline}
\sum_{\ell=1}^{t-1}(-1)\Big(d\alpha_{(q-j-1)(q-j-\ell-1)}+\alpha_{(q-j-1)(q-j-\ell)}\Big)
(-1)^{j-1}d^{j+\ell-t}{j-1\choose j+\ell-t}\\=
\sum_{\ell=1}^{t-2}\alpha_{(q-j-1)(q-j-\ell-1)}(-1)^{j}d^{j+1+\ell-t}\Big[{j\choose
j+1+\ell-t}-{j-1\choose
j+\ell+1-t}\Big]\\+\alpha_{(q-j-1)(q-j-t)}(-d)^{j}+\alpha_{(q-j-1)(q-j-1)}
(-1)^{j}d^{j+1-t}{j-1\choose
j+1-t}\\+\sum_{\ell=1}^{t-2}\alpha_{(q-j-1)(q-j-\ell-1)}
(-1)^{j}d^{j+\ell+1-t}{j-1\choose
j+\ell+1-t}\hspace{1.7cm}\\+(-1)^{j-1}d^{j-t}\Big({n+1\choose
q-j}-d\alpha_{(q-j-1)(q-j-1)}\Big){j-1\choose
j-t}.\hspace{1.7cm}\label{2.44}\end{multline} Cancelling out the
same terms in \eqref{2.44}, we obtain
\begin{multline}
\sum_{\ell=0}^{t-1}\alpha_{(q-j)(q-j-\ell)}(-1)^{j-1}d^{j+\ell-t}{j-1\choose
j+\ell-t}\\=
\sum_{\ell=0}^{t-1}\alpha_{(q-j-1)(q-j-\ell-1)}(-1)^{j}d^{j+1+\ell-t}{j\choose
j+1+\ell-t}\\+(-1)^{j-1}d^{j-t}{n+1\choose q-j}{j-1\choose
j-t}\qquad\text{for $j\geq t$}.\hspace{2cm}\label{2.46}
\end{multline}
Similarly to \eqref{2.46}, we have
\begin{multline}
\sum_{\ell=t-j}^{t-1}\alpha_{(q-j)(q-j-\ell)}(-1)^{j-1}d^{j+\ell-t}{j-1\choose
j+\ell-t}\\=
\sum_{\ell=t-j-1}^{t-1}\alpha_{(q-j-1)(q-j-\ell-1)}(-1)^{j}d^{j+1+\ell-t}{j\choose
j+1+\ell-t}\quad\text{for $j<t$}.\label{2.461}
\end{multline}
For each $1\leq s\leq q-2t+1$, we divide the following summation
into two parts:
 \begin{multline}
\lefteqn{\sum_{\ell=0}^{t-1}\sum_{j=t-\ell}^{q-\ell}\alpha_{(q-j)(q-j-\ell)}(-1)^{j-1}d^{j+\ell-t}{j-1\choose
j+\ell-t}}\\=\sum_{\ell=0}^{t-1}\sum_{j=t+s}^{q-\ell}\alpha_{(q-j)(q-j-\ell)}(-1)^{j-1}d^{j+\ell-t}{j-1\choose
j+\ell-t}\\+\sum_{j=1}^{t+s-1}\sum_{\ell=\max\{0,t-j\}}^{t-1}\alpha_{(q-t)(q-t-\ell)}(-1)^{j-1}d^{j+\ell-t}{j-1\choose
j+\ell-t}.\label{2.462}
\end{multline}
Then use \eqref{2.461} and \eqref{2.46} to get
\begin{multline}
\sum_{j=1}^{t}\sum_{\ell=0}^{t-1}\alpha_{(q-t)(q-t-\ell)}(-1)^{j-1}d^{j+\ell-t}{j-1\choose
j+\ell-t}\\=t\sum_{\ell=0}^{t-1}\alpha_{(q-t)(q-t-\ell)}(-1)^{t-1}d^{\ell}{t-1\choose
\ell}\hspace{3.6cm}\\=t\sum_{\ell=0}^{t-1}\alpha_{(q-t-1)(q-t-1-\ell)}(-1)^{t}d^{\ell+1}{t\choose
\ell+1}+t(-1)^{t-1}{n+1\choose q-t}.\label{2.47}
\end{multline}
By \eqref{2.47}, we have
\begin{multline}
\sum_{j=1}^{t+s-1}\sum_{\ell=\max\{0,t-j\}}^{t-1}\alpha_{(q-t)(q-t-\ell)}(-1)^{j-1}d^{j+\ell-t}{j-1\choose
j+\ell-t}\\=
t\sum_{\ell=0}^{t-1}\alpha_{(q-t-1)(q-t-1-\ell)}(-1)^{t}d^{\ell+1}{t\choose
\ell+1}+t(-1)^{t-1}{n+1\choose
q-t}\\+\sum_{j=t+1}^{t+s-1}\sum_{\ell=0}^{t-1}\alpha_{(q-t)(q-t-\ell)}(-1)^{j-1}d^{j+\ell-t}{j-1\choose
j+\ell-t}\\=
(t+1)\sum_{\ell=0}^{t-1}\alpha_{(q-(t+1))(q-(t+1)-\ell)}(-1)^{t}d^{\ell+1}{t\choose
\ell+1}+t(-1)^{t-1}{n+1\choose
q-t}\\+\sum_{j=t+2}^{t+s-1}\sum_{\ell=0}^{t-1}\alpha_{(q-t)(q-t-\ell)}(-1)^{j-1}d^{j+\ell-t}{j-1\choose
j+\ell-t}.\label{2.4712}
\end{multline}
Substituting \eqref{2.46} in \eqref{2.4712}, we obtain
\begin{multline}
\sum_{j=1}^{t+s-1}\sum_{\ell=\max\{0,t-j\}}^{t-1}\alpha_{(q-t)(q-t-\ell)}(-1)^{j-1}d^{j+\ell-t}{j-1\choose
j+\ell-t}\\=
(t+2)\sum_{\ell=0}^{t-1}\alpha_{(q-(t+2))(q-(t+2)-\ell)}(-1)^{t+1}d^{\ell+2}{t+1\choose
\ell+2}\\+\sum_{j=t}^{t+1}(-1)^{j-1}d^{j-t}{n+1\choose
q-j}{j-1\choose
j-t}\\+\sum_{j=t+3}^{t+s-1}\sum_{\ell=0}^{t-1}\alpha_{(q-t)(q-t-\ell)}(-1)^{j-1}d^{j+\ell-t}{j-1\choose
j+\ell-t}.\label{2.4711}
\end{multline}
 Repeat the procedure
from \eqref{2.4712} to \eqref{2.4711} and put it back in
\eqref{2.462} to get
\begin{multline}
\sum_{\ell=0}^{t-1}\sum_{j=t-\ell}^{q-\ell}\alpha_{(q-j)(q-j-\ell)}(-1)^{j-1}d^{j+\ell-t}{j-1\choose
j+\ell-t}\\=\sum_{\ell=0}^{t-1}\sum_{j=t+s}^{q-\ell}\alpha_{(q-j)(q-j-\ell)}(-1)^{j-1}d^{j+\ell-t}{j-1\choose
j+\ell-t}\hspace{2cm}\\+(t+s-1)\sum_{\ell=0}^{t-1}\alpha_{(q-(t+s-1))(q-(t+s-1)-\ell)}(-1)^{t-s-2}d^{s+\ell-1}{t-s-2\choose
s+\ell-1}
\\+\sum_{j=t}^{t+s-2}j(-1)^{j-1}d^{j-t}{n+1\choose
q-j}{j-1\choose j-t}\qquad\text {for $1\leq s\leq
q-2t+1$}.\label{2.471}
\end{multline}
For $s=q-2t+1$, we change indices in the first term on the right
hand side of \eqref{2.471} to get
\begin{multline}
\sum_{\ell=0}^{t-1}\sum_{j=q-t+1}^{q-\ell}\alpha_{(q-j)(q-j-\ell)}(-1)^{j-1}d^{j+\ell-t}{j-1\choose
j+\ell-t}\\=\sum_{j=0}^{t-1}\sum_{\ell=0}^{t-1-j}\alpha_{(t-1-j)(t-1-j-\ell)}(-1)^{q-t+j}d^{q-2t+j+\ell+1}{q-t+j\choose
q-2t+j+\ell+1}.\label{2.473}
\end{multline}
With the same procedure as we do in \eqref{2.46}, we have
\begin{multline}
\sum_{\ell=0}^{t-1-\lambda}\alpha_{(t-1-\lambda)(t-1-\lambda-\ell)}(-1)^{q-t+\lambda}d^{q-2t+\lambda+\ell+1}{q-t+\lambda\choose
q-2t+\lambda+\ell+1}\\=\sum_{\ell=0}^{t-2-\lambda}\alpha_{(t-2-\lambda)(t-2-\lambda-\ell)}(-1)^{q-t+\lambda+1}d^{q-2t+\lambda+\ell+2}{q-t+\lambda+1\choose
q-2t+\lambda+\ell+2}\\+(-1)^{q-t+j}d^{q-2t+\lambda+\ell+1}{q-t+\lambda\choose
q-2t+\lambda+\ell+1}{n+1\choose t-1-\lambda}.\label{2.472}
\end{multline}
for $j=0,1,\cdots,t-1$. Put \eqref{2.473} back to \eqref{2.471}
and use \eqref{2.472} to repeat the procedure as \eqref{2.471}
\begin{multline}
\sum_{\ell=0}^{t-1}\sum_{j=t-\ell}^{q-\ell}\alpha_{(q-j)(q-j-\ell)}(-1)^{j-1}d^{j+\ell-t}{j-1\choose
j+\ell-t}\\=\sum_{j=0}^{t-1}\sum_{\ell=0}^{t-1-j}\alpha_{(t-1-j)(t-1-j-\ell)}(-1)^{q-t+j}d^{q-2t+j+\ell+1}{q-t+j\choose
q-2t+j+\ell+1}\\+(q-t)\sum_{\ell=0}^{t-1}\alpha_{(t-1)(t-1-\ell)}(-1)^{q-t}d^{q-2t+\ell+1}{q-t\choose
q-2t+\ell+1}\\ +\sum_{j=t}^{q-t}j(-1)^{j-1}d^{j-t}{n+1\choose
q-j}{j-1\choose j-t}\hspace{2cm}\\=
\sum_{j=\lambda}^{t-1}\sum_{\ell=0}^{t-1-j}\alpha_{(t-1-j)(t-1-j-\ell)}(-1)^{q-t+j}d^{q-2t+j+\ell+1}{q-t+j\choose
q-2t+j+\ell+1}\\+(q-t+\lambda)\sum_{\ell=0}^{t-1-\lambda}\alpha_{(t-1-\lambda)(t-1-\lambda-\ell)}
(-1)^{q-t+\lambda}d^{q-2t+\lambda+\ell+1}{q-t+\lambda\choose
q-2t+\lambda+\ell+1}\\
+\sum_{j=t}^{q-t+\lambda}j(-1)^{j-1}d^{j-t}{n+1\choose
q-j}{j-1\choose j-t}\hspace{3cm}\label{2.474}
\end{multline}
 for $\lambda=1,\cdots,t-1$. By substituting $\lambda=t-1$ in
 \eqref{2.474}, we get
\begin{multline}
\sum_{\ell=0}^{t-1}\sum_{j=t-\ell}^{q-\ell}\alpha_{(q-j)(q-j-\ell)}(-1)^{j-1}d^{j+\ell-t}{j-1\choose
j+\ell-t}\\=
 \sum_{j=t}^{q}j(-1)^{j-1}d^{j-t}{n+1\choose
q-j}{j-1\choose j-t}=-t\alpha_{q(q-t)}.\hspace{3.5cm}\label{2.48}
\end{multline}
By substituting \eqref{2.48} in \eqref{2.42}, we get the
coefficient of $\rm{div}(X)$
\begin{multline}
\sum_{j=1}^{q}(-1)^{j-1}\psi_{j-1}\sum_{k=0}^{q-j}\alpha_{(q-j)k}\omega^{k}
\big(\frac{\sqrt{-1}}{2\pi}\partial\bar{\partial}\xi\big)^{q-j-k}\\=
\sum_{t=1}^{q}\sum_{\ell=0}^{t-1}\sum_{j=t-\ell}^{q-\ell}\alpha_{(q-j)(q-j-\ell)}(-1)^{j-1}d^{j+\ell-t}{j-1\choose
j+\ell-t}\omega^{q-t}\big(\frac{\sqrt{-1}}{2\pi}\partial\bar{\partial}\xi\big)^{t-1}\hspace{1cm}\\
=-\sum_{t=1}^{q}t\alpha_{q(q-t)}\omega^{q-t}\big(\frac{\sqrt{-1}}{2\pi}\partial\bar{\partial}\xi\big)^{t-1}
=-\sum_{t=0}^{q-1}(q-t)\alpha_{qt}\omega^{t}\big(\frac{\sqrt{-1}}{2\pi}\partial\bar{\partial}\xi\big)^{q-t-1}.\label{2.49}
\end{multline}
The coefficient for $\theta$ in \eqref{2.33} is
\begin{equation}
-(n+1)\sum_{j=2}^{q}(-1)^{j-1}(\psi_{j-1}-\omega^{j-1})\sum_{k=0}^{q-j}\alpha_{(q-j)k}\omega^{k}
\big(\frac{\sqrt{-1}}{2\pi}\partial\bar{\partial}\xi\big)^{q-j-k}.\label{2.50}
\end{equation}
Use \eqref{2.15} in the second part of
 \eqref{2.50} and we get\begin{multline}
(n+1)\sum_{j=1}^{q}(-1)^{j}\sum_{k=0}^{q-j}\alpha_{(q-j)k}\omega^{k+j-1}
\big(\frac{\sqrt{-1}}{2\pi}\partial\bar{\partial}\xi\big)^{q-j-k}
\\=-(n+1)\sum_{t=0}^{q-1}\sum_{k=0}^{t}\alpha_{(q-t-1+k)k}(-1)^{t-k}\omega^{t}
\big(\frac{\sqrt{-1}}{2\pi}\partial\bar{\partial}\xi\big)^{q-t-1}\hspace{2.5cm}\\=
(-1)^{q}(n+1)\sum_{t=0}^{q-1}\sum_{k=0}^{t}{q-1-k\choose
t-k}{n\choose k}(-1)^{k}d^{t-k}\omega^{t}
\big(\frac{\sqrt{-1}}{2\pi}\partial\bar{\partial}\xi\big)^{q-t-1}
.\label{2.51}\end{multline} Use \eqref{2.15} to compute
\begin{multline}
d(q-t)\alpha_{qt}-(t+1)\alpha_{q(t+1)}\\=(-1)^{q}d(q-t)\sum_{k=0}^{t}{q-k\choose
t-k}{n+1\choose
k}d^{t-k}(-1)^{k}\hspace{4cm}\\-(-1)^{q}(t+1)\sum_{k=0}^{t+1}{q-k\choose
t+1-k}{n+1\choose k}d^{t+1-k}(-1)^{k}\\=(-1)^{q}
\sum_{k=0}^{t-1}(-1)^{k+1}\left((q-t){q-1-k\choose
t-1-k}-(t+1){q-1-k\choose t-k}\right){n+1\choose k+1}\\
+(-1)^{q+t}(t+1){q-1-t\choose 0}{n+1\choose t+1}\hspace{1.5cm}\\
=(-1)^{q}\sum_{k=0}^{t}(-1)^{k}(n+1){q-1-k\choose t-k}{n\choose
k}.\hspace{5cm}\label{2.52}
\end{multline}
Comparing the coefficient in \eqref{2.52} and \eqref{2.51}, we get
\begin{multline}
(n+1)\sum_{j=1}^{q}(-1)^{j}\sum_{k=0}^{q-j}\alpha_{(q-j)k}\omega^{k+j-1}
\big(\frac{\sqrt{-1}}{2\pi}\partial\bar{\partial}\xi\big)^{q-j-k}
\\=-\sum_{t=0}^{q-1}\left(d(q-t)\alpha_{qt}-(t+1)\alpha_{q(t+1)}\right)\omega^{t}
\big(\frac{\sqrt{-1}}{2\pi}\partial\bar{\partial}\xi\big)^{q-t-1}.\label{2.53}
\end{multline} Substituting \eqref{2.49} and \eqref{2.53} in
\eqref{2.51}, the coefficient for $\theta$ in \eqref{2.33} is
\begin{equation}
\sum_{t=0}^{q-1}\left((n+1-d)(q-t)\alpha_{qt}+(t+1)\alpha_{q(t+1)}\right)\omega^{t}
\big(\frac{\sqrt{-1}}{2\pi}\partial\bar{\partial}\xi\big)^{q-t-1}.\label{2.54}\end{equation}
 By using Sublemma \ref{sub2.4} and \eqref{2.35}, the coefficient for
$-\partial\bar{\partial}\theta$ in \eqref{2.33} is
\begin{multline}
\frac{\sqrt{-1}}{2\pi}\sum_{j=2}^{q}(-1)^{j-1}\left(n\sum_{\ell=0}^{j-2}\big(d\omega+
\frac{\sqrt{-1}}{2\pi}\partial\bar{\partial}\xi\big)^{j-2-\ell}\omega^{\ell}+\omega^{j-2}\right)\\
\cdot\sum_{k=0}^{q-j}\alpha_{(q-j)k}\omega^{k}
\big(\frac{\sqrt{-1}}{2\pi}\partial\bar{\partial}\xi\big)^{q-j-k}
=\sum_{j=2}^{q}(-1)^{j-1}\big(n\zeta_{j-2}+\omega^{j-2}\big)\wedge\eta_{q-j},
\label{2.55}
\end{multline}
where
\begin{eqnarray*}
\begin{array}{lcl}
\eta_{q-j}&=&\sum_{k=0}^{q-j}\alpha_{(q-j)k}\omega^{k}
\big(\frac{\sqrt{-1}}{2\pi}\partial\bar{\partial}\xi\big)^{q-j-k}\\
\zeta_{j-2}&=&\frac{\sqrt{-1}}{2\pi}\sum_{k=0}^{j-2}
\big(d\omega+\frac{\sqrt{-1}}{2\pi}\partial\bar{\partial}\xi\big)^{j-2-k}\omega^{k}.
\end{array}
\end{eqnarray*}
By using Sublemma \ref{sub2.4} and \eqref{2.35}, the coefficient
for $\partial\bar{\partial}\triangle\theta$ in \eqref{2.33} is
\begin{multline}
\frac{\sqrt{-1}}{2\pi}\sum_{j=2}^{q}(-1)^{j-1}\sum_{\ell=0}^{j-2}\big(d\omega+
\frac{\sqrt{-1}}{2\pi}\partial\bar{\partial}\xi\big)^{j-2-\ell}\omega^{\ell}\sum_{k=0}^{q-j}\alpha_{(q-j)k}\omega^{k}
\big(\frac{\sqrt{-1}}{2\pi}\partial\bar{\partial}\xi\big)^{q-j-k}\\
=\sum_{j=2}^{q}(-1)^{j-1}\zeta_{j-2}\wedge\eta_{q-j}.\hspace{7cm}\label{2.56}
\end{multline}
By adding
$-\rm{div}(X)\eqref{2.49}+\theta\eqref{2.54}+\partial\bar{\partial}\theta
\eqref{2.55}+\partial\bar{\partial}\triangle\theta\eqref{2.56}$,
we obtain \newline $q\tilde{P}^{q}(\nabla
X,\Theta,\cdots,\Theta)$.
\end{proof}
\begin{lemma}\label{lem2.5}The Hodge decomposition of equation \eqref{2.30} can be computed as follows:
\begin{equation}-q\tilde{P}^{q}(\nabla X,\Theta,\cdots,\Theta)+q\alpha_{qq}\theta\omega^{q-1}-i(X)\partial
f_{q}=-\kappa\alpha_{q(q-1)}\omega^{q-1}+\bar{\partial}\varphi_{q},\label{2.57}\end{equation}
where
\begin{eqnarray*}\varphi_{q}
&=&\kappa\frac{\sqrt{-1}}{2\pi}\sum_{k=0}^{q-2}(q-k)
\omega^{k}\wedge(\partial\xi)\wedge\big(\frac{\sqrt{-1}}{2\pi}\partial\bar{\partial}\xi\big)^{q-k-2}\\&&
+\frac{\sqrt{-1}}{2\pi}\sum_{k=1}^{q-1}k\alpha_{qk}\theta(\partial\xi)\wedge
\omega^{k-1}\wedge\big(\frac{\sqrt{-1}}{2\pi}\partial\bar{\partial}\xi\big)^{q-k-1}\\&&-
\frac{\sqrt{-1}}{2\pi}\sum_{k=0}^{q-1}\alpha_{qk}(q-k-1)X(\xi)\partial\xi
\omega^{k}\wedge\big(\frac{\sqrt{-1}}{2\pi}\partial\bar{\partial}\xi\big)^{q-k-2}\\&&
-\sum_{j=1}^{q}(-1)^{j-1}\left[(\partial\theta)
\wedge\big(\omega^{j-1} +n\zeta_{j-2}\big)
-(\partial\triangle\theta)\wedge
\zeta_{j-2}\right]\wedge\eta_{q-j}
\end{eqnarray*} is a globally defined form.
\end{lemma}
\begin{proof} In \eqref{2.27}, the potential form is obtained
\[f_{q}=\sum_{k=0}^{q-1}\alpha_{qk}\xi
\omega^{k}\wedge\big(\frac{\sqrt{-1}}{2\pi}\partial\bar{\partial}\xi\big)^{q-k-1}\]
for $q=1,\cdots,n-1$. By direction computation, we get
\begin{eqnarray}
\lefteqn{i(X)\partial f_{q}=\sum_{k=0}^{q-1}
\alpha_{qk}X(\xi)
\omega^{k}\wedge\big(\frac{\sqrt{-1}}{2\pi}\partial\bar{\partial}\xi\big)^{q-k-1}}
\nonumber\\&&+\sum_{k=1}^{q-1}k\alpha_{qk}\frac{\sqrt{-1}}{2\pi}(\partial\xi)
\wedge(\bar{\partial}\theta)
\wedge\omega^{k-1}(\frac{\sqrt{-1}}{2\pi}\partial\bar{\partial}\xi)^{q-k-1}\nonumber\\&&-
\sum_{k=0}^{q-2}(q-k-1)\alpha_{qk}\frac{\sqrt{-1}}{2\pi}(\partial\xi)\wedge
\Big(\bar{\partial}X(\xi)\Big)\wedge\omega^{k}\wedge\big(\frac{\sqrt{-1}}{2\pi}\partial\bar{\partial}\xi\big)^{q-k-2}.
\label{2.58}
\end{eqnarray}
Use Lemma \ref{lem2.4} and \eqref{2.58}, we obtain
\begin{eqnarray}
\lefteqn{-q\tilde{P}^{q}(\nabla
X,\Theta,\cdots,\Theta)+q\alpha_{qq}\theta\omega^{q-1}-i(X)\partial
f_{q}}\nonumber\\&&=
\rm{div}(X)\sum_{k=0}^{q-1}(q-k)\alpha_{qk}\omega^{k}
\wedge\big(\frac{\sqrt{-1}}{2\pi}\partial\bar{\partial}\xi\big)^{q-1-k}\nonumber\\&&
-\theta\sum_{k=0}^{q-1}\left[(q-k)(n+1-d)\alpha_{qk}+(k+1)\alpha_{q(k+1)}\right]
\omega^{k}\wedge\big(\frac{\sqrt{-1}}{2\pi}\partial\bar{\partial}\xi\big)^{q-1-k}\nonumber\\&&
+\sum_{j=1}^{q}(-1)^{j-1}\left[(\partial\bar{\partial}\theta)
\wedge\big(\omega^{j-1} +n\zeta_{j-2}\big)
-(\partial\bar{\partial}\triangle\theta)\wedge
\zeta_{j-2}\right]\wedge\eta_{q-j}\nonumber\\&&
+q\alpha_{qq}\theta\omega^{q-1}
-X(\xi)\sum_{k=0}^{q-1}[(q-k)-(q-k-1)]\alpha_{qk}\omega^{k}\wedge\big(\frac{\sqrt{-1}}{2\pi}\partial\bar{\partial}\xi\big)^{q-1-k}\nonumber\\&&
-\frac{\sqrt{-1}}{2\pi}\sum_{k=1}^{q-1}k\alpha_{qk}(\partial\xi)
\wedge(\bar{\partial}\theta)
\wedge\omega^{k-1}(\frac{\sqrt{-1}}{2\pi}\partial\bar{\partial}\xi)^{q-k-1}\nonumber\\&&+
\frac{\sqrt{-1}}{2\pi}\sum_{k=0}^{q-2}(q-k-1)\alpha_{qk}(\partial\xi)\wedge
\Big(\bar{\partial}X(\xi)\Big)\wedge\omega^{k}\wedge\big(\frac{\sqrt{-1}}{2\pi}\partial\bar{\partial}\xi\big)^{q-k-2}.
\label{2.59}
\end{eqnarray}
In fact, we have the formula \cite[Theorem 4.1]{LZ1}
\begin{equation}\rm{div}(X)-X(\xi)-(n-d+1)\theta=-\kappa.\label{2.60}\end{equation}
By using \eqref{2.60} in \eqref{2.59}, we have
\begin{eqnarray}
\lefteqn{-q\tilde{P}^{q}(\nabla
X,\Theta,\cdots,\Theta)+q\alpha_{qq}\theta\omega^{q-1}-i(X)\partial
f_{q}}\nonumber\\&&=
-\kappa\sum_{k=0}^{q-1}(q-k)\alpha_{qk}\omega^{k}
\wedge\big(\frac{\sqrt{-1}}{2\pi}\partial\bar{\partial}\xi\big)^{q-1-k}\nonumber\\&&
-\theta\sum_{k=1}^{q-1}k\alpha_{qk}\omega^{k-1}\wedge
\big(\frac{\sqrt{-1}}{2\pi}\partial\bar{\partial}\xi\big)^{q-1-k}\nonumber\\&&
-\frac{\sqrt{-1}}{2\pi}\sum_{k=1}^{q-1}k\alpha_{qk}(\partial\xi)
\wedge(\bar{\partial}\theta)
\wedge\omega^{k-1}(\frac{\sqrt{-1}}{2\pi}\partial\bar{\partial}\xi)^{q-k-1}\nonumber\\&&
+\sum_{j=1}^{q}(-1)^{j-1}\left[(\partial\bar{\partial}\theta)
\wedge\big(\omega^{j-1} +n\zeta_{j-2}\big)
-(\partial\bar{\partial}\triangle\theta)\wedge
\zeta_{j-2}\right]\wedge\eta_{q-j}\nonumber\\&&
+X(\xi)\sum_{k=0}^{q-2}(q-k-1)\alpha_{qk}\omega^{k}
\wedge\big(\frac{\sqrt{-1}}{2\pi}\partial\bar{\partial}\xi\big)^{q-1-k}\nonumber\\&&+
\frac{\sqrt{-1}}{2\pi}\sum_{k=0}^{q-2}(q-k-1)\alpha_{qk}(\partial\xi)\wedge
\Big(\bar{\partial}X(\xi)\Big)\wedge\omega^{k}\wedge\big(\frac{\sqrt{-1}}{2\pi}\partial\bar{\partial}\xi\big)^{q-k-2}.\label{2.61}
\end{eqnarray}
From Lemma \ref{Lef}, the harmonic part on the left hand side of
\eqref{2.57} is proportional to $\omega^{q-1}$. Therefore, the
Hodge decomposition of \eqref{2.61} is
$-\kappa\alpha_{q(q-1)}\omega^{q-1}+\bar{\partial}\varphi_{q}$.
\end{proof}
\begin{lemma}\label{lem2.7}We define $q\tilde{P}^{q}(\nabla
X,\Theta,\cdots,\Theta)$ in \eqref{2.33}. For $q=1,\cdots$, $n-1$,
we have
\begin{equation}
\int_{M}q\tilde{P}^{q}(\nabla
X,\Theta,\cdots,\Theta)\wedge\omega_{M}^{n-q}=0,\label{2.62}
\end{equation}
where $\omega_{M}=(n+1-d)\omega$.
\end{lemma}
\begin{proof}
 Consider
\begin{eqnarray}
0&=&\int_{M}L_{X}(P^{q-1}(\Theta)\wedge\omega^{n-q})\nonumber\\&
=&\int_{M}\big(L_{X}P^{q-1}(\Theta)\big)\wedge\omega^{n-q}
+\int_{M}P^{q-1}(\Theta)\wedge L_{X}(\omega^{n-q})\nonumber\\&=&
-(q-1)\frac{\sqrt{-1}}{2\pi}\int_{M}\partial\bar{\partial}\tilde{P}^{q-1}(\nabla
X,\Theta,\cdots,\Theta)\wedge\omega^{n-q}\nonumber\\&&
-(n-q)\frac{\sqrt{-1}}{2\pi}\int_{M}P^{q-1}(\Theta)\wedge\partial\bar{\partial}\theta\wedge\omega^{n-q-1}.\label{2.63}
\end{eqnarray}
By Stoke's Theorem, we have
$\int_{M}\partial\bar{\partial}\tilde{P}^{q-1}(\nabla
X,\Theta,\cdots,\Theta)\wedge\omega^{n-q}=0$. By using \eqref{11}
and \eqref{2.311}, compute that
\begin{multline}
-(n-q)\frac{\sqrt{-1}}{2\pi}\partial\bar{\partial}\theta\wedge
P^{q-1}(\Theta)\wedge\omega^{n-q-1}\\=-\frac{(n-q)(n-1-q)!}{(q-1)!(n-1)!}
\sum_{\sigma\in S_{q-1},\tau\in
S_{q}}\sum_{i_{1},j_{1},k_{1}\cdots,k_{q-1},i_{q},j_{q}=2}^{n}\rm{sgn}(\sigma)\rm{sgn}(\tau)\\
\cdot\partial_{i_{1}}\bar{\partial}_{j_{1}}\theta
R_{k_{1}i_{2}\bar{j}_{2}}^{k_{\sigma(1)}}\cdots
R_{k_{q-1}i_{q}\bar{j}_{q}}^{k_{\sigma(q-1)}}\tilde{g}^{i_{\tau(1)}\bar{j}_{1}}\cdots\tilde{g}^{i_{\tau(q)}\bar{j}_{q}}
\omega^{n-1}\\=\frac{(n-q)!}{(q-1)!(n-1)!} \sum_{\sigma\in
S_{q-1},\tau\in
S_{q}}\sum_{i_{1},k_{1},s_{1},\cdots,s_{q-1},i_{q},k_{q-1}=2}^{n}\rm{sgn}(\sigma)\rm{sgn}(\tau)\\
\qquad\cdot X_{i_{1}}^{i_{\tau(1)}}
R_{i_{2}k_{1}\bar{s}_{1}}^{i_{\tau(2)}}\cdots
R_{i_{q}k_{q-1}\bar{s}_{q-1}}^{i_{\tau(q)}}\tilde{g}^{k_{\sigma(1)}\bar{s}_{1}}\cdots\tilde{g}^{k_{\sigma(q-1)}\bar{s}_{q-1}}
\omega^{n-1}.\label{2.64}
\end{multline}
By definition, we have \begin{multline} \tilde{P}^{q}(\nabla
X,\Theta,\cdots,\Theta)=\frac{1}{q!}\big(\frac{\sqrt{-1}}{2\pi}\big)^{q-1}\sum_{\tau\in
S_{q}}\sum_{i_{1},k_{1},s_{1}\cdots,k_{q-1},s_{q-1},i_{q}=2}^{n}{\rm sgn}(\tau)\\
\cdot
X_{i_{1}}^{i_{\tau(1)}}R_{i_{2}k_{1}\bar{s}_{1}}^{i_{\tau(2)}}\cdots
R_{i_{q}k_{q-1}\bar{s}_{q-1}}^{i_{\tau(q)}}dz_{k_{1}}\wedge
d\bar{z}_{s_{1}}\wedge\cdots dz_{k_{q-1}}\wedge
d\bar{z}_{s_{q-1}}.\label{2.641}\end{multline} By comparing
\eqref{2.64} and \eqref{2.641}, we obtain that
\begin{equation*}-(n-q)\frac{\sqrt{-1}}{2\pi}\partial\bar{\partial}\theta\wedge
P^{q-1}(\Theta)\wedge\omega^{n-q-1}=q \tilde{P}^{q}(\nabla
X,\Theta,\cdots,\Theta)\wedge
\omega^{n-q}.\end{equation*}Therefore, we have
\begin{equation*}
\int_{M}q\tilde{P}^{q}(\nabla
X,\Theta,\cdots,\Theta)\wedge\omega_{M}^{n-q}=
\int_{M}L_{X}(P^{q-1}(\Theta)\wedge\omega_{M}^{n-q})=0.
\end{equation*}
\end{proof}
For an extension discussion of Lemma \ref{lem2.7}, please see the
Appendix.
\begin{proof}[Proof of Theorem \ref{1.1}.]
From \eqref{2.31}, we can compute the $q$-th Bando-Futaki
invariant by
\begin{multline}
\int_{M}L_{X}f_{q}\wedge\omega_{M}^{n-q}=(n+1-d)^{n-q}\left(-C(q)\int_{M}\omega^{n-1}+q\alpha_{qq}\int_{M}\theta
\omega^{n-1}\right)\\-q\int_{M}\tilde{P}^{q}(\nabla
X,\Theta,\cdots,\Theta)\wedge\omega_{M}^{n-q}
-\int_{M}\bar{\partial}\varphi_{q}\wedge \omega_{M}^{n-q}.
\end{multline}
From Lemma \ref{lem2.5}, we have $C(q)=-\kappa\alpha_{q(q-1)}$.
 By
\cite[Theorem 5.1]{LZ1}, we have
$\int_{M}\theta\omega^{n-1}=\frac{\kappa}{n}$. From Lemma
\ref{lem2.7}, we have $q\int_{M}\tilde{P}^{q}(\nabla
X,\Theta,\cdots,\Theta)\wedge\omega_{M}^{n-q}=0$. In addition, the
divergence theorem implies that
$\int_{M}\bar{\partial}\varphi_{q}\wedge \omega_{M}^{n-q}=0$. With
these information, we get
\[\int_{M}L_{X}f_{q}\wedge\omega_{M}^{n-q}=\kappa(n+1-d)^{n-q}(\alpha_{qq}\frac{q}{n}+
d\alpha_{q(q-1)}).\] Using \eqref{2.15}, the $q$-th Bando-Futaki
invariant can be computed directly by
\begin{eqnarray*}
\mathcal{F}_{q}(X)&=&-(n+1-d)^{n-q}\frac{(d-1)}{n}\sum_{j=0}^{q-1}(-d)^{j}(j+1)(q-j){n+1\choose q-j}\kappa\\
&=&-(n+1-d)^{n-q}\frac{(d-1)(n+1)}{n}\sum_{j=0}^{q-1}(-d)^{j}(j+1){n\choose q-1-j}\kappa.
\end{eqnarray*}
\end{proof}
\newtheorem{example}{Example}[section]
\begin{corollary}
With the notation as in Theorem 1.1, if $d=1$, all the
Bando-Futaki invariants are zero.
\end{corollary}
\begin{corollary}
With the notation as in Theorem 1.1, if $M$ is the cubic surface
in $\mathbb{CP}^{3}$, then the Bando-Futaki invariants are
\[\mathcal{F}_{1}(X)=-\frac{8}{3}\kappa,
\text{and}\qquad\mathcal{F}_{2}(X)=16\kappa.\]
\end{corollary}
The first example is due to \cite{D-T} and \cite{LZ1}.
\begin{example}
Let
\[M=\{Z\in\mathbb{CP}^{3}|F(Z)=Z_{0}Z_{1}^{2}+Z_{2}^{2}Z_{3}-Z_{2}Z_{3}^{2}=0\}\] be a hypersurface in $\mathbb{CP}^{3}$,
where $Z=[Z_{0},Z_{1},Z_{2},Z_{3}]$ is homogeneous coordinates of
$\mathbb{CP}^{3}$. $M$ is a K\"ahler orbifold with $c_{1}(M)>0$.
Let \[X=-7Z_{0}\frac{\partial}{\partial
Z_{0}}+5Z_{1}\frac{\partial}{\partial
Z_{1}}+Z_{2}\frac{\partial}{\partial
Z_{2}}+Z_{3}\frac{\partial}{\partial Z_{3}}\] be a holomorphic
vector field in $\mathbb{CP}^{3}$, which satisfies $XF=3F$. The
Bando-Futaki invariants are
\[\mathcal{F}_{1}(X)=-8, \text{and}\qquad\mathcal{F}_{2}(X)=48.\]
\end{example}
\begin{example}
Let
\[M=\{Z\in\mathbb{CP}^{4}|F(Z)=Z_{0}^{2}+Z_{1}^{2}+Z_{2}^{2}+Z_{3}^{2}+Z_{4}^{2}=0\}\]
be a hypersurface in $\mathbb{CP}^{4}$. Let
\[X=Z_{0}\frac{\partial}{\partial
Z_{0}}+Z_{1}\frac{\partial}{\partial
Z_{1}}+Z_{2}\frac{\partial}{\partial
Z_{2}}+Z_{3}\frac{\partial}{\partial
Z_{3}}+Z_{4}\frac{\partial}{\partial Z_{4}}\] be a holomorphic
vector field on $\mathbb{CP}^{4}$, which satisfies \[X(F)=2F.\]
The Bando-Futaki invariants are
\[\mathcal{F}_{1}(X)=-270,\quad\mathcal{F}_{2}(X)=135,\quad\text{and}F_{3}(X)=-\frac{285}{2}.\]
\end{example}
\section{Chen and Tian's holomorphic invariants}
The holomorphic invariants were introduced by Chen and Tian \cite{C-T}. We prove that they are the Futaki invariants.
\begin{definition}
Let $M$ be an $n$-dimensional simply-connected K\"ahler manifold
with a K\"ahler form $\omega$. There exists a smooth function
$\theta_{X}$ such that
$i(X)\omega=\frac{\sqrt{-1}}{2\pi}\bar{\partial}\theta_{X}$\footnote{In
order to maintain the definition as the original paper, we have
opposite sign of the notation for $\theta_{X}$ and $\Delta$ that
we used in previous section.}. Define
\begin{eqnarray}\mathcal{F}_{k}(X,\omega)&=&
(n-k)\int_{M}\theta_{X}\omega^{n}
+(k+1)\int_{M}\Delta\theta_{X}Ric(\omega)^{k}\wedge \omega^{n-k}
\nonumber\\&&-(n-k)\int_{M}\theta_{X}Ric(\omega)^{k+1}\wedge
\omega^{n-k-1}.\label{3.1}
\end{eqnarray}
\end{definition}
\begin{proof}[Proof of \text{Theorem 1.2}]
These new holomorphic invariants are independent of the choices of
the K\"ahler metrics in the K\"ahler class $[\omega]$, which were
shown in \cite{C-T}. There exists a constant $\alpha$, such that
$\alpha \omega\in c_{1}(M)$. Therefore, there exists a smooth real
valued function $f$ over $M$, such that $Ric(\omega)-\alpha
\omega=\frac{\sqrt{-1}}{2\pi}\partial\bar{\partial}f$. Take inner
derivative on both sides, we have
\begin{equation}divX+\alpha\theta_{X}+X(f)=\beta,\label{30}\end{equation}
where $\beta$ is a constant if $M$ is compact. We need the
following two formulas
\begin{eqnarray}
0&=&\int_{M}\left(i(X)[\partial f(\partial\bar{\partial}f)^{j-1}\omega^{n-j+1}]\right)\nonumber\\
&=&j\int_{M}X(f)(\partial\bar{\partial}f)^{j-1}\omega^{n-j+1}+(n-j+1)\int_{M}\bar{\partial}\theta_{X}\partial
f(\partial\bar{\partial}f)^{j-1}\omega^{n-j}\nonumber\\&=&
j\int_{M}X(f)(\partial\bar{\partial}f)^{j-1}\omega^{n-j+1}+(n-j+1)\int_{M}\theta_{X}(\partial\bar{\partial}f)^{j}\omega^{n-j}
\label{31}
\end{eqnarray}for $1\leq j\leq k+1$, and
\begin{equation}divX=\Delta\theta_{X}+c,\label{32}\end{equation} where $c$ is a constant and
$\Delta\theta_{X}=g^{i\bar{j}}\partial_{i}\bar{\partial}_{j}\theta_{X}$
if $\omega=\frac{\sqrt{-1}}{2\pi}\sum_{i,j=1}^{n}g_{i\bar{j}}dz_{i}\wedge d\bar{z}_{j}$.
The new holomorphic invariants are
\begin{eqnarray*}
\mathcal{F}_{k}(X,\omega)&=&(n-k)\int_{M}\theta_{X} \omega^{n}+
(k+1)\int_{M}\Delta\theta_{X} Ric(\omega)^{k}\wedge
\omega^{n-k}\\&&-(n-k)\int_{M}\theta_{X} Ric(\omega)^{k+1}\wedge
\omega^{n-k-1}\\&=& (n-k)(1-\alpha^{k+1})\int_{M}\theta_{X}
\omega^{n}\\&&+(k+1)\int_{M}\Delta\theta_{X}\sum_{i=1}^{k}{k\choose
i}(\frac{\sqrt{-1}}{2\pi}\partial\bar{\partial}f)^{i}\alpha^{k-i}\omega^{n-i}\\&&
-(n-k)\int_{M}\theta_{X}\sum_{i=1}^{k+1}{k+1\choose i}(\alpha
\omega)^{k+1-i}(\frac{\sqrt{-1}}{2\pi}\partial\bar{\partial}f)^{i}\omega^{n-k-1}.
\end{eqnarray*}
For each $1\leq i\leq k+1$, break the number into two terms
$(n-i+1)=(n-k)+(k-i-1)$. By using \eqref{31}, we get
\begin{multline}-(n-k)\int_{M}\theta_{X}(\frac{\sqrt{-1}}{2\pi}\partial\bar{\partial}f)^{i}\omega^{n-i}\\=
(k-i-1)\int_{M}\theta_{X}(\frac{\sqrt{-1}}{2\pi}\partial\bar{\partial}f)^{i}\omega^{n-i}+i
\int_{M}X(f)(\partial\bar{\partial}f)^{i-1}\omega^{n-i+1}.\label{3.5}
\end{multline}
 By using \eqref{32} and \eqref{3.5}, we have
\begin{eqnarray}
\mathcal{F}_{k}(X,\omega)&=&(n-k)(1-\alpha^{k+1})\int_{M}\theta_{X}
\omega^{n}\nonumber\\&&+(k+1)\int_{M}divX\sum_{i=1}^{k}{k\choose
i}(\frac{\sqrt{-1}}{2\pi}\partial\bar{\partial}f)^{i}\alpha^{k-i}\omega^{n-i}\nonumber\\&&+
\int_{M}\sum_{i=1}^{k+1}i{k+1\choose
i}X(f)(\frac{\sqrt{-1}}{2\pi}\partial\bar{\partial}f)^{i-1}\alpha^{k+1-i}\omega^{n-i+1}\nonumber\\&&
+\int_{M}\sum_{i=1}^{k+1}(k-i+1){k+1\choose
i}\theta_{X}(\frac{\sqrt{-1}}{2\pi}\partial\bar{\partial}f)^{i}\alpha^{k+1-i}\omega^{n-i}.\label{3.6}
\end{eqnarray}
Since $\sum_{i=1}^{k+1}i{k+1\choose
i}=\sum_{i=0}^{k}(i+1){k+1\choose
i+1}=(k+1)+\sum_{i=1}^{k}(k+1){k\choose i}$, we devide the term of
$X(f)$ in \eqref{3.6} into two parts and combine one part with
$\rm{div}(X)$. That is,
\begin{eqnarray}
\mathcal{F}_{k}(X,\omega)&=&(n-k)(1-\alpha^{k+1})\int_{M}\theta_{X}
\omega^{n}\nonumber\\&&+(k+1)\int_{M}\left(divX+X(f)\right)\sum_{i=1}^{k}{k\choose
i}(\frac{\sqrt{-1}}{2\pi}\partial\bar{\partial}f)^{i}\alpha^{k-i}\omega^{n-i}\nonumber\\&&+
(k+1)\alpha^{k}\int_{M}X(f)\omega^{n}\nonumber\\&&+
\int_{M}\sum_{i=1}^{k+1}(k+1){k\choose
i}\theta_{X}(\frac{\sqrt{-1}}{2\pi}\partial\bar{\partial}f)^{i}\alpha^{k+1-i}\omega^{n-i}.\label{3.7}
\end{eqnarray}
By using \eqref{30} in \eqref{3.7}, we get
\begin{eqnarray*}
\mathcal{F}_{k}(X,\omega)&=&(n-k)(1-\alpha^{k+1})\int_{M}\theta_{X}
\omega^{n}+
(k+1)\alpha^{k}\int_{M}X(f)\omega^{n}\\&&-(k+1)\int_{M}\theta_{X}\sum_{i=1}^{k}{k\choose
i}(\frac{\sqrt{-1}}{2\pi}\partial\bar{\partial}f)^{i}\alpha^{k-i+1}\omega^{n-i}\\&&+
\int_{M}\sum_{i=1}^{k}(k+1){k\choose
i}\theta_{X}(\frac{\sqrt{-1}}{2\pi}\partial\bar{\partial}f)^{i}\alpha^{k+1-i}\omega^{n-i}\\&
=&(n-k)(1-\alpha^{k+1})\int_{M}\theta_{X}
w^{n}+(k+1)\alpha^{k}\int_{M}X(f)\omega^{n}.
\end{eqnarray*}
If the K\"ahler form is normalized, then we can choose $\alpha=1$.
These holomorphic invariants are simply the Futaki invariants on
compact K\"ahler manifolds with $c_{1}(M)>0$. Indeed Kobayashi
\cite{KS1} proved that a compact K\"ahler manifold with positive
Ricci curvature is simply connected.
\end{proof}
The generalized energy functionals introduced in the same paper
are the nonlinearizations of these holomorphic invariants. The
Futaki invariant can have different nonlinearizations.
\section{Higher order K-energy Functionals}
In 1986, Mabuchi first introduced K-energy as the nonlinearization
of the Futaki invariant \cite{MT}. The critical point of the
K-energy functional is the K\"ahler-Einstein form. K-energy are
studied to understand the stability of K\"ahler manifolds by Tian
\cite{TG1,TG2, TG3}, Phong, and Sturm \cite{P-S1, P-S2}.
Furthermore, Lu \cite{LZ2} provided the K-energy in an explicit
formula for the hypersurface in the projective spaces. Phong and
Sturm \cite{P-S2} formulized it on complete intersections using
the Deligne pairing technique. Moreover, Bando and Mabuchi
constructed higher-order K-energy functionals \cite{B-M}, which
are considered as nonlinearizations of the Bando-Futaki invariants
\cite{B-M}. (~cf. Theorem 2 of Weinkove's \cite{WB}) However, we
can remove Weinkove's assumption, which states that the
$qth$-Chern form $c_{q}(\omega)$ is in the same cohomology class
as $\mu_{q}[\omega^{q}]\in\mathrm{H}^{2q}(M,\mathbb{Z})$ where
$\omega$ is the K\"ahler form and $\mu_{q}$. Most importantly, he
\cite{WB} derived the higher order K-energy as a generalization of
Tian's formula of K-energy \cite{TG1}. Bando and Mabuchi's proof
\cite{B-M} is discussed in detail in the following proof
concerning the independence of the choice of paths of higher order
K-energy functionals in the K\"ahler class by using Mabuchi's
method \cite{MT}.
\begin{definition} Let $M$ be a connected compact $n$-dimensional K\"ahler manifold with positive first Chern class.
Let $\Omega$ be the K\"ahler class which represents the first Chern form.
For any $\omega_{0}, \omega_{1}\in \Omega$,
let $\omega_{t}$, $0\leq t\leq 1$, be a curve joining $\omega_{0}$ and $\omega_{1}$.
Since $M$ is K\"ahler, there exists a smooth real valued function $\varphi_{t}$
such that $\omega_{t}=\omega_{0}+\frac{\sqrt{-1}}{2\pi}\partial\bar{\partial}\varphi_{t}$
with $\int_{M}\frac{d}{dt}\varphi_{t}\omega_{t}^{n}=0$.
Define higher order K-energy functionals as
\begin{equation}M_{q}(\omega_{0},\omega_{1})
=\frac{1}{V}\int_{0}^{1}\int_{M}\frac{d \varphi_{t}}{dt}(c_{q}(\omega_{t})-Hc_{q}(\omega_{t}))\wedge
\omega_{t}^{n-q}dt,\label{33}\end{equation} where
$V=\int_{M}\omega_{FS}|_{M}^{n}$ and $Hc_{q}(\omega_{t})$ is the
harmonic part of $c_{q}(\omega)$.
\end{definition} The independence of path choosing in the K\"ahler class for K-energy functionals was proved Mabuchi
\cite{MT}, and for the higher order K-energy functionals was proved by Bando and Mabuchi \cite{B-M} when
$1\leq q\leq n$. Recently, Weinkove gave an alternative derivation of the proof by using Bott-Chern forms.
\par
Let us re-prove the argument of Bando-Mabuchi in detail.
First, we need
\begin{sublem}[Bando \cite{BS}]
$Hc_{q}(\omega)\wedge \omega^{n-q}$ is harmonic if $\omega\in\Omega$.
\end{sublem}
\begin{proof}
We may use Lefschetz decomposition theorem,
\[Hc_{q}(\omega)=\sum_{k=0}^{q}\omega^{k}\wedge\varphi_{k},\]
where $\varphi_{k}\in \mathcal{H}^{2q-2k}(M,\mathbb{C})$ is the
primitive $2(q-k)$-form of $Hc_{q}(\omega)$. Therefore,
$Hc_{q}(\omega)\wedge\omega^{n-q}=\sum_{k=0}^{q}\omega^{n-q+k}\wedge\varphi_{k}$
is harmonic. We can also see this in a different method. Let
$L\eta=\omega\wedge\eta$. We know $[\Delta, L]=0$. Since $\Delta
Hc_{q}(\omega)=0$, we have
\begin{eqnarray*}
\Delta(Hc_{q}(\omega)\wedge\omega^{n-q})
&=&\Delta(\omega^{n-q}\wedge Hc_{q}(\omega))\\
&=&\omega^{n-q}\wedge\Delta(Hc_{q}(\omega))\\
&=&0.
\end{eqnarray*}
Since $\mathrm{dim}(\mathcal{H}^{2n}(M,\mathbb{C}))=1$ and
$Hc_{q}(\omega)\wedge\omega^{n-q}\in \mathcal{H}^{2n}(M)$,
$Hc_{q}(\omega)\wedge\omega^{n-q}=\lambda_{q}\omega^{n}$. Since $M$
is compact, $\lambda_{q}$ must be a constant.
\end{proof}
Hence,
\begin{eqnarray*}\int_{M}c_{q}(\omega)\wedge\omega^{n-q}&=&
\int_{M}(Hc_{q}(\omega)+\frac{\sqrt{-1}}{2\pi}\partial\bar{\partial}f_{q})\wedge\omega^{n-q}
\\&=&\int_{M}Hc_{q}(\omega)\wedge\omega^{n-q}=\lambda_{q}\int_{M}\omega^{n}.\end{eqnarray*}
We can conclude that Definition 4.1 is the same as in \cite{B-M}
\begin{equation}\frac{1}{V}\int_{0}^{1}\int_{M}\frac{d\varphi_{t}}{dt}(c_{q}(\omega_{t})-Hc_{q}(\omega_{t}))\wedge
\omega_{t}^{n-q}dt=\frac{1}{V}\int_{0}^{1}\int_{M}\frac{d\varphi_{t}}{dt}(c_{q}(\omega_{t})\wedge\omega_{t}^{n-q}
-\lambda_{q}\omega_{t}^{n})dt.\label{34}\end{equation}
\begin{sublem}
$\int_{0}^{1}\int_{M}\frac{d\varphi_{t}}{dt}\lambda_{q}\omega_{t}^{n}dt$
is independent of path choosing in the K\"ahler class.
\end{sublem}
It is trivial by the following method.
\begin{sublem}
$\int_{0}^{1}\int_{M}\frac{d\varphi_{t}}{dt}c_{q}(\omega_{t})\wedge\omega^{n-q}$ is independent of path choosing in the K\"ahler class.
\end{sublem}
\begin{proof}
Let $\omega_{0}=(n-d+1)\omega_{FS}|_{M}$,
$\omega_{s,t}=\omega_{0,0}+\frac{\sqrt{-1}}{2\pi}\psi_{s,t}$ and
$\psi_{s,t}=s\varphi_{t}$, where $\varphi_{t}(z)\in
C^{\infty}([0,1]\times M)$. Let
\begin{equation}\Psi_{s,t}^{q}=\left(\int_{M}\frac{\partial
\psi_{s,t}}{\partial s}c_{q}(\omega_{s,t})\wedge
\omega_{s,t}^{n-q}\right)ds+\left(\int_{M}\frac{\partial
\psi_{s,t}}{\partial t}c_{q}(\omega_{s,t})\wedge
\omega_{s,t}^{n-q}\right)dt.\label{40}\end{equation} Use Stoke's
theorem
\begin{eqnarray}
\int_{0}^{1}\int_{0}^{1}d\Psi_{s,t}^{q}&=&
-\int_{0}^{1}\left(\int_{M}\frac{\partial \psi_{s,t}}{\partial
s}c_{q}(\omega_{s,t})\wedge
\omega_{s,t}^{n-q}\right)ds\Big|_{t=0}^{t=1}\nonumber\\&&+
\int_{0}^{1}\left(\int_{M}\frac{\partial \psi_{s,t}}{\partial
t}c_{q}(\omega_{s,t})\wedge
\omega_{s,t}^{n-q}\right)\Big|_{s=0}^{s=1}dt\nonumber\\&=&
-\int_{0}^{1}\left(\int_{M}\varphi_{t}c_{q}(\omega_{s,t})\wedge
\omega_{s,t}^{n-q}\right)ds\Big|_{t=0}^{t=1}\nonumber\\&&-
\int_{0}^{1}\left(\int_{M}\dot{\varphi}_{t}c_{q}(\omega_{s,t})\wedge
\omega_{s,t}^{n-q}\right)dt.\label{35}
\end{eqnarray}
\begin{sublem}\label{sub4.4}
For the one form $\Psi_{s,t}^{q}$ in \eqref{40}, we have $d\Psi_{s,t}^{q}=0$.
\end{sublem}
Suppose that Sublemma \ref{sub4.4} is true. Let
$\varphi_{0}=\varphi_{1}$ in \eqref{36} to get
\[\int_{0}^{1}\left(\int_{M}\dot{\varphi}_{t}c_{q}(\omega_{s,t})\wedge
\omega_{s,t}^{n-q}\right)dt=-\int_{0}^{1}\left(\int_{M}\varphi_{t}c_{q}(\omega_{s,t})\wedge
\omega_{s,t}^{n-q}\right)ds\Big|_{t=0}^{t=1}=0\] This shows that it
is independent of path choosing of $\omega_{t}$.
\end{proof}
\begin{proof}[Proof of the Sublemma \ref{sub4.4}]
By further computation, we have
\begin{eqnarray}d\Psi^{q}&=&-\int_{M}\frac{\partial}{\partial t
}\left(\frac{\partial \psi_{s,t}}{\partial
s}c_{q}(\omega_{s,t})\wedge \omega_{s,t}^{n-q}\right)ds\wedge dt
\nonumber\\&&+\int_{M}\frac{\partial}{\partial
s}\left(\frac{\partial \psi_{s,t}}{\partial
t}c_{q}(\omega_{s,t})\wedge \omega_{s,t}^{n-q}\right)ds\wedge
dt\nonumber\\&=& -\int_{M}\frac{\partial \psi_{s,t}}{\partial
s}\tilde{P}^{q}(-\bar{\partial}(\nabla \frac{\partial
h_{s,t}}{\partial t} h_{s,t}^{-1}),R_{s,t},\cdots,R_{s,t})\wedge
\omega_{s,t}^{n-q}ds\wedge
dt\nonumber\\&&-(n-q)\int_{M}\frac{\partial \psi_{s,t}}{\partial s}
c_{q}(\omega_{s,t})\wedge
\partial\bar{\partial}\frac{\partial \psi_{s,t}}{\partial
t}\omega_{s,t}^{n-1-q}ds\wedge dt\nonumber\\&&
+\int_{M}\frac{\partial \psi_{s,t}}{\partial
t}\tilde{P}^{q}(-\bar{\partial}(\nabla \frac{\partial
h_{s,t}}{\partial s} h_{s,t}^{-1}),R_{s,t},\cdots,R_{s,t})\wedge
\omega_{s,t}^{n-q}ds\wedge
dt\nonumber\\&&-(n-q)\int_{M}\frac{\partial \psi_{s,t}}{\partial
t}c_{q}(\omega_{s,t}) \wedge
\partial\bar{\partial}\frac{\partial \psi_{s,t}}{\partial
s}\omega_{s,t}^{n-1-q}ds\wedge dt\nonumber\\&=&
-\int_{M}\partial\bar{\partial}\frac{\partial \psi_{s,t}}{\partial
s}\tilde{P}^{q}(\frac{\partial h_{s,t}}{\partial
t}h_{s,t}^{-1},R_{s,t},\cdots,R_{s,t})\wedge
\omega_{s,t}^{n-q}ds\wedge
dt\nonumber\\&&+(n-q)\int_{M}\partial\frac{\partial
\psi_{s,t}}{\partial s}c_{q}(\omega_{s,t})\wedge
\bar{\partial}\frac{\partial \psi_{s,t}}{\partial
t}\omega_{s,t}^{n-1-q}ds\wedge
dt\nonumber\\&&+\int_{M}\partial\bar{\partial}\frac{\partial
\psi_{s,t}}{\partial t}\tilde{P}^{q}(\frac{\partial
h_{s,t}}{\partial s}h_{s,t}^{-1},R_{s,t},\cdots,R_{s,t})\wedge
\omega_{s,t}^{n-q}ds\wedge
dt\nonumber\\&&+(n-q)\int_{M}\bar{\partial}\frac{\partial
\psi_{s,t}}{\partial t}c_{q}(\omega_{s,t})\wedge
\partial\frac{\partial \psi_{s,t}}{\partial
s}\omega_{s,t}^{n-1-q}ds\wedge dt,\label{36}
\end{eqnarray}
where
$\omega_{s,t}=\frac{\sqrt{-1}}{2\pi}\sum_{\alpha,\beta}(h_{s,t})_{\alpha\bar{\beta}}dz_{\alpha}\wedge
d\bar{z}_{\beta}$ and
$R_{s,t}=\frac{\sqrt{-1}}{2\pi}\bar{\partial}[(\partial
h_{s,t})h_{s,t}^{-1})]$ is the curvature form with respect to metric
$\omega_{s,t}$. We need to show that
\begin{multline*}\int_{M}\partial\bar{\partial}\frac{\partial
\psi_{s,t}}{\partial s}\tilde{P}^{q}(\frac{\partial
h_{s,t}}{\partial t}h_{s,t}^{-1},R_{s,t},\cdots,R_{s,t})\wedge
\omega_{s,t}^{n-q}\\=\int_{M}\partial\bar{\partial}\frac{\partial
\psi_{s,t}}{\partial t}\tilde{P}^{q}(\frac{\partial
h_{s,t}}{\partial s}h_{s,t}^{-1},R_{s,t},\cdots,R_{s,t})\wedge
\omega_{s,t}^{n-q}\end{multline*} to conclude $d\Psi^{q}=0$. Compute
\begin{eqnarray}
\lefteqn{\partial\bar{\partial}\frac{\partial \psi_{s,t}}{\partial
s}\tilde{P}^{q}(\frac{\partial h_{s,t}}{\partial
t}h_{s,t}^{-1},R_{s,t},\cdots,R_{s,t})\wedge
\omega_{s,t}^{n-q}}\nonumber\\&&= \tilde{P}^{q}(\frac{\partial
h_{s,t}}{\partial t}h_{s,t}^{-1}\partial\bar{\partial}\frac{\partial
\psi_{s,t}}{\partial s},R_{s,t},\cdots,R_{s,t})\wedge
\omega_{s,t}^{n-q}\nonumber\\&&=\lefteqn{\frac{1}{q!}
\sum_{\sigma,\tau\in S_{q}} sgn(\sigma)sgn(\tau)\frac{\partial
(h_{s,t})_{i_{1}\bar{j}}}{\partial
t}h_{s,t}^{i_{\sigma(1)}\bar{j}}\partial_{\alpha_{1}}\bar{\partial}_{\beta_{1}}\frac{\partial
\psi_{s,t}}{\partial
s}}\nonumber\\&&\qquad\times\left((R_{s,t})_{i_{2}\alpha_{2}\bar{\beta}_{2}}^{i_{\sigma(2)}}\cdots
(R_{s,t})_{i_{q}\alpha_{q}\bar{\beta}_{q}}^{i_{\sigma(q)}}\right)
(h_{s,t})^{\alpha_{\tau(1)}\bar{\beta}_{1}}\cdots
(h_{s,t})^{\alpha_{\tau(q)}\bar{\beta}_{q}}\omega_{s,t}^{n}.\label{4.5}
\end{eqnarray}
Since $\frac{\partial (h_{s,t})_{i_{1}\bar{j}}}{\partial t}
=\partial_{i_{1}}\bar{\partial}_{j}\frac{\partial
\psi_{s,t}}{\partial t}$, \eqref{4.5} is equal to
\begin{eqnarray*}
\lefteqn{ \frac{1}{q!}\sum_{\sigma,\tau\in
S_{q}}sgn(\sigma)sgn(\tau)\partial_{i_{1}}\bar{\partial}_{j}\frac{\partial
\psi_{s,t}}{\partial
t}(h_{s,t})^{i_{\sigma(1)}\bar{j}}\partial_{\alpha_{1}}\bar{\partial}_{\beta_{1}}\frac{\partial
\psi_{s,t}}{\partial
s}}\\&&\times\left((R_{s,t})_{i_{2}\bar{\eta}_{2}\alpha_{2}\bar{\beta}_{2}}(h_{s,t})^{i_{\sigma(2)}\bar{\eta}_{2}}\cdots
(R_{s,t})_{i_{q}\bar{\eta}_{q}\alpha_{q}\bar{\beta}_{q}}(h_{s,t})^{i_{\sigma(q)}\bar{\eta}_{q}}
\right)\\&&\times(h_{s,t})^{\alpha_{\tau(1)}\bar{\beta}_{1}}
\cdots(h_{s,t})^{\alpha_{\tau(k)}\bar{\beta}_{q}}\omega_{s,t}^{n}\\&&=
\frac{1}{q!}\sum_{\sigma,\tau\in
S_{q}}sgn(\sigma)sgn(\tau)\partial_{i_{1}}\bar{\partial}_{j}\frac{\partial
\psi_{s,t}}{\partial
t}(h_{s,t})^{i_{\sigma(1)}\bar{j}}\partial_{\alpha_{1}}\bar{\partial}_{\beta_{1}}\frac{\partial
\psi_{s,t}}{\partial
s}\\&&\times\left((R_{s,t})_{\alpha_{2}\bar{\beta}_{2}i_{2}\bar{\eta}_{2}}(h_{s,t})^{i_{\sigma(2)}\bar{\eta}_{2}}\cdots
(R_{s,t})_{\alpha_{q}\bar{\beta}_{q}i_{q}\bar{\eta}_{q}}(h_{s,t})^{i_{\sigma(q)}\bar{\eta}_{q}}
\right)\\&&\times(h_{s,t})^{\alpha_{\tau(1)}\bar{\beta}_{1}}
\cdots(h_{s,t})^{\alpha_{\tau(k)}\bar{\beta}_{q}}\omega_{s,t}^{n}\\&&=
\frac{1}{q!}\sum_{\sigma,\tau\in
S_{q}}sgn(\sigma)sgn(\tau)\partial_{i_{1}}\bar{\partial}_{j}\frac{\partial
\psi_{s,t}}{\partial
t}\partial_{\alpha_{1}}\bar{\partial}_{\beta_{1}}\frac{\partial
\psi_{s,t}}{\partial
s}(h_{s,t})^{\alpha_{\tau(1)}\bar{\beta_{1}}}\\&&\times\left((R_{s,t})_{\alpha_{2}i_{2}\bar{\eta}_{2}}^{\alpha_{\tau(2)}}\cdots
(R_{s,t})_{\alpha_{q}i_{q}\bar{\eta}_{q}}^{\alpha_{\tau(q)}}\right)(h_{s,t})^{i_{\sigma(1)}\bar{j}}(h_{s,t})^{i_{\sigma(2)}\bar{\eta}_{2}}\cdots(h_{s,t})^{i_{\sigma(q)}\bar{\eta}_{q}}w_{s,t}^{n}\\&&=\partial\bar{\partial}\frac{\partial
\psi_{s,t}}{\partial t}\tilde{P}^{q}(\frac{\partial
h_{s,t}}{\partial s}h_{s,t}^{-1},R_{s,t},\cdots,R_{s,t})\wedge
\omega_{s,t}^{n-q}.
\end{eqnarray*}
\end{proof}
We restate and clarify as follows:
\begin{lemma} $($\cite{B-M, WB}$)$
Higher order K-energy functionals are the nonlinearizations of Bando-Futaki invariants.
\begin{equation}
\frac{2}{V}\mathrm{Re}(\mathcal{F}_{q}(X))=(n+1-q)\frac{d}{dt}M_{q}(\omega_{0},\omega_{t})\label{37}
\end{equation}
\end{lemma}
Let $M$ be an $n$-dimensional compact connected K\"ahler manifold
in $\mathbb{CP}^{N}$ with positive first Chern class. There exists
a constant $\alpha>0$ such that $\alpha \omega_{FS}|_{M}\in
c_{1}(M)$, where $\omega_{FS}$ is the Fubini-Study metric in
$\mathbb{CP}^{N}$. Let $\sigma_{t}$ be a one-parameter family of
automorphism of $\mathbb{CP}^{N}$ and $X$ be the holomorphic
vector field induced by $\sigma_{t}$. We may write
\[\sigma_{t}[Z_{0},\cdots, Z_{N}]=[e^{\lambda_{0}t}Z_{0},\cdots,e^{\lambda_{N}t}Z_{N}]\]
for integers $\lambda_{0},\cdots,\lambda_{N}$ with
$\sum_{i=0}^{N}\lambda_{i}=0$. Then $\omega_{t}=\alpha
\sigma_{t}^{\ast}\omega_{FS}|_{M}$ restricts a family of metrics
on $M$, such that $w_{0}=\alpha \omega_{FS}|_{M}$. Recall
$\omega_{FS}=\frac{\sqrt{-1}}{2\pi}\partial\bar{\partial}\log(\sum_{i=0}^{N}|Z_{i}|^{2})$.
Hence,
\\ $\sigma_{t}^{\ast}\omega_{FS}=\frac{\sqrt{-1}}{2\pi}\partial\bar{\partial}\log(\sum_{i=0}^{N}|e^{\lambda_{i}t}Z_{i}|^{2})$.
Let
\[\varphi_{t}=\alpha\log\left(\frac{\sum_{i=0}^{N}|e^{\lambda_{i}t}Z_{i}|^{2}}{\sum_{i=0}^{N}|Z_{i}|^{2}}\right).\]
It follows
\[\omega_{t}-\omega_{0}=\frac{\sqrt{-1}}{2\pi}\partial\bar{\partial}\varphi_{t}.\]
Then
\[\frac{d\varphi_{t}}{dt}=\frac{2\alpha\mathrm{Re}\sum_{i=0}^{N}\lambda_{i}e^{\lambda_{i}t}Z_{i}\overline{e^{\lambda_{i}t}Z_{i}}}{\sum_{i=0}^{N}e^{\lambda_{i}t}Z_{i}}=-2\mathrm{Re}(\alpha\theta\circ\sigma_{t}),\]
where
$i(X)\omega_{FS}=-\frac{\sqrt{-1}}{2\pi}\bar{\partial}\theta$, and
$\theta=-\frac{\sum_{i=0}^{N}|\lambda_{i}Z_{i}|^{2}}{\sum_{i=0}^{N}|Z_{i}|^{2}}$.
From \cite{BS} and Lemma 4.1 in \cite{WB}, the Bando-Futaki
invariants can be written as
\[\mathcal{F}_{q}(X)=-(n+1-q)\int_{M}\alpha\theta(c_{q}(\omega)-Hc_{q}(\omega))\wedge\omega^{n-q},\]
where $\omega=\alpha\omega_{FS}|_{M}$.
\begin{eqnarray*}
\lefteqn{(n+1-q)\frac{d}{dt}M_{q}(\omega,\omega_{t})}\\&&=
(n+1-q)\frac{1}{V}\int_{M}\frac{d\varphi_{t}}{dt}(c_{q}(\omega_{t})-Hc_{q}(\omega_{t}))\wedge\omega_{t}^{n-q}\\&&=
-(n+1-q)\frac{1}{V}\int_{M}2\mathrm{Re}(\alpha\theta\circ\sigma_{t})
(c_{q}(\omega_{t})-Hc_{q}(\omega_{t}))\wedge\omega_{t}^{n-q}\\&&=
-(n+1-q)\frac{2}{V}\mathrm{Re}\left(\int_{M}\alpha\theta(c_{q}(\omega)-Hc_{q}(\omega))\wedge\omega^{n-q}\right)\\&&=
\frac{2}{V}\mathrm{Re}(\mathcal{F}_{q}(X)),
\end{eqnarray*}
since Bando-Futaki invariants are independent of the choices of
metrics in the K\"ahler class.
  \section{APPENDIX}  

  \begin{lemma}\label{lem5.1}For $q=2,\cdots,n-1$,
let \[Y=\sum_{i_{1}=2}^{n}Y^{i_{1}}\frac{\partial}{\partial
z_{i_{1}}}\in
T^{1,0}(M)\otimes\wedge^{q-1}(T^{1,0}(M)^{\ast}\otimes
\overline{T^{1,0}(M)}^{\ast})\] be a holomorphic vector field with
$(q-1,q-1)$-valued form, where
\[Y^{i_{1}}=X^{i_{1}}P^{q-1}(\Theta)+\sum_{j=2}^{q}\sum_{i_{2},\cdots,i_{j}=2}^{n}
(-1)^{j-1}X^{i_{2}}\Theta_{i_{2}}^{i_{3}}\cdots\Theta_{i_{j}}^{i_{1}}P^{q-j}(\Theta).\]
Then we have $q\tilde{P}^{q}(\nabla
X,\Theta,\cdots,\Theta)=\rm{div}(Y)$. For $q=1$, we have $Y=X$ and
$\tilde{P}^{1}(\nabla X)={\rm div}(X)$.
\end{lemma}
\begin{proof}
We have \begin{eqnarray} {\rm
div}(Y)&=&\sum_{i_{1}=2}^{n}\nabla_{i_{1}}Y^{i_{1}}\nonumber\\&=&
\sum_{i_{1}=2}^{n}X_{i_{1}}^{i_{1}}P^{q-1}+\sum_{i_{1}=2}^{n}X^{i_{1}}\nabla_{i_{1}}P^{q-1}\nonumber\\&&
+\sum_{j=2}^{q}\sum_{i_{1},\cdots,i_{j}=2}^{n}
(-1)^{j-1}X_{i_{1}}^{i_{2}}\Theta_{i_{2}}^{i_{3}}\cdots\Theta_{i_{j}}^{i_{1}}P^{q-j}(\Theta)\nonumber
\\&&+\sum_{j=2}^{q}\sum_{i_{1},\cdots,i_{j}=2}^{n}
(-1)^{j-1}X^{i_{2}}\nabla_{i_{1}}\big(\Theta_{i_{2}}^{i_{3}}\cdots\Theta_{i_{j}}^{i_{1}}\big)P^{q-j}(\Theta)\nonumber
\\&&+\sum_{j=2}^{q}\sum_{i_{1},\cdots,i_{j}=2}^{n}
(-1)^{j-1}X^{i_{2}}\Theta_{i_{2}}^{i_{3}}\cdots\Theta_{i_{j}}^{i_{1}}\nabla_{i_{1}}(P^{q-j}(\Theta)).\label{5.1}
\end{eqnarray}
By using the definition of $q\tilde{P}^{q}(\nabla
X,\Theta,\cdots,\Theta)$ in \eqref{2.33}, we get
\begin{eqnarray}
{\rm div}(Y)&=& q\tilde{P}^{q}(\nabla
X,\Theta,\cdots,\Theta)+\sum_{i_{1}=2}^{n}X^{i_{1}}\nabla_{i_{1}}P^{q-1}\nonumber\\&&
-\sum_{j=2}^{q}\sum_{i_{1},\cdots,i_{j}=2}^{n}(-1)^{j}X^{i_{2}}\nabla_{i_{1}}(\Theta_{i_{2}}^{i_{3}}\cdots\Theta_{i_{j}}^{i_{1}})P^{q-j}(\Theta)
\nonumber\\&&-\sum_{j=2}^{q}\sum_{i_{1},\cdots,i_{j}=2}^{n}(-1)^{j}
X^{i_{2}}\Theta_{i_{2}}^{i_{3}}\cdots\Theta_{i_{j}}^{i_{1}}\nabla_{i_{1}}P^{q-j}(\Theta).\label{5.2}
\end{eqnarray}
It is equivalent to show the following Sublemma.
\begin{sublem}\label{sub5.1}
For $q=2,\cdots,n-1$, we have
\begin{eqnarray}
\sum_{i_{1}=2}^{n}X^{i_{1}}\nabla_{i_{1}}P^{q-1}&=&
\sum_{j=2}^{q}\sum_{i_{1},\cdots,
i_{j}=2}^{n}(-1)^{j}X^{i_{2}}\nabla_{i_{1}}(\Theta_{i_{2}}^{i_{3}}\cdots
\Theta_{i_{j}}^{i_{1}})P^{q-j}(\Theta)\nonumber\\&&+
\sum_{j=2}^{q}\sum_{i_{1},\cdots,i_{j}=2}^{n}(-1)^{j}
X^{i_{2}}\Theta_{i_{2}}^{i_{3}}\cdots\Theta_{i_{j}}^{i_{1}}\nabla_{i_{1}}P^{q-j}(\Theta).\label{5.3}
\end{eqnarray} For $q=1$, \eqref{5.3} is zero.
\end{sublem}
We need two formulas before proving Sublemma \ref{sub5.1}.
\begin{sublem}\label{sub5.2}
For $2\leq i\leq n$ and $q=1,\cdots,n-1$, we compute the covariant
derivative on the curvature form and the Chern forms
\begin{eqnarray}
\nabla_{i}\Theta_{j}^{k}&=&\nabla_{j}\Theta_{i}^{k},\label{5.4}\\
\nabla_{i}P^{q}(\Theta)&=&\sum_{\ell=1}^{q}\sum_{i_{1},\cdots,i_{\ell}=2}^{n}(-1)^{\ell-1}
\frac{1}{\ell}\nabla_{i}
(\Theta_{i_{1}}^{i_{2}}\Theta_{i_{2}}^{i_{3}}\cdots\Theta_{i_{\ell}}^{i_{1}})P^{q-\ell}(\Theta),\label{5.5}
\end{eqnarray}
where
$\Theta_{j}^{k}=\frac{\sqrt{-1}}{2\pi}\sum_{s,t=2}^{n}R_{js\bar{t}}^{k}dz_{s}\wedge
d\bar{z}_{t}=-\frac{\sqrt{-1}}{2\pi}\sum_{s,t=2}^{n}\bar{\partial}_{t}\Gamma_{sj}^{k}dz_{s}\wedge
d\bar{z}_{t}$. Since $M$ is K\"ahler,
$\Gamma_{sj}^{k}=\Gamma_{js}^{k}=\frac{\partial\tilde{g}_{s\bar{\beta}}}{\partial
z_{j}}\tilde{g}^{k\bar{\beta}}$.
\end{sublem}
\begin{proof}
First, we compute
\begin{eqnarray}\nabla_{i}R_{js\bar{t}}^{k}&=&
-\partial_{i}\bar{\partial}_{t}\Gamma_{js}^{k}+\Gamma_{i\lambda}^{k}R_{js\bar{t}}^{\lambda}-
\Gamma_{ij}^{\lambda}R_{\lambda
s\bar{t}}^{k}-\Gamma_{is}^{\lambda}R_{j\lambda
\bar{t}}^{k}\nonumber\\&=&
-\bar{\partial}_{t}\left(\frac{\partial^{2}\tilde{g}_{s\bar{\beta}}}{\partial
z_{i}\partial
z_{j}}\tilde{g}^{k\bar{\beta}}\right)+\bar{\partial}_{t}\left(\frac{\partial\tilde{g}_{s\bar{\beta}}}{\partial
z_{j}}\frac{\partial\tilde{g}_{\alpha\bar{\eta}}}{\partial
z_{i}}\tilde{g}^{k\bar{\eta}}\tilde{g}^{\alpha\bar{\beta}}\right)\nonumber\\&&+
\Gamma_{i\lambda}^{k}R_{js\bar{t}}^{\lambda}-
\Gamma_{ij}^{\lambda}R_{\lambda
s\bar{t}}^{k}-\Gamma_{is}^{\lambda}R_{j\lambda
\bar{t}}^{k}\nonumber\\&=&
-\bar{\partial}_{t}\left(\frac{\partial^{2}\tilde{g}_{s\bar{\beta}}}{\partial
z_{i}\partial
z_{j}}\tilde{g}^{k\bar{\beta}}\right)-\Gamma_{i\alpha}^{k}R_{js\bar{t}}^{\alpha}
-\Gamma_{js}^{\alpha}R_{i\alpha \bar{t}}^{k}\nonumber\\&&+
\Gamma_{i\lambda}^{k}R_{js\bar{t}}^{\lambda}-
\Gamma_{ij}^{\lambda}R_{\lambda
s\bar{t}}^{k}-\Gamma_{is}^{\lambda}R_{j\lambda
\bar{t}}^{k}.\label{5.6}
\end{eqnarray}
Similarly, we get
\begin{eqnarray}
\nabla_{j}R_{is\bar{t}}^{k}&=&
-\bar{\partial}_{t}\left(\frac{\partial^{2}\tilde{g}_{s\bar{\beta}}}{\partial
z_{j}\partial
z_{i}}\tilde{g}^{k\bar{\beta}}\right)-\Gamma_{j\alpha}^{k}R_{is\bar{t}}^{\alpha}
-\Gamma_{is}^{\alpha}R_{j\alpha \bar{t}}^{k}\nonumber\\&&+
\Gamma_{j\lambda}^{k}R_{is\bar{t}}^{\lambda}-
\Gamma_{ji}^{\lambda}R_{\lambda
s\bar{t}}^{k}-\Gamma_{js}^{\lambda}R_{i\lambda
\bar{t}}^{k}.\label{5.7}
\end{eqnarray}
Since $M$ is K\"ahler,
$\Gamma_{ij}^{\lambda}=\Gamma_{ji}^{\lambda}$. Therefore,
\eqref{5.6} and \eqref{5.7} are equal. Then we get \eqref{5.4}.
Next, we prove \eqref{5.5} by induction. For $q=1$,
$\nabla_{i}P^{1}(\Theta)=\nabla_{i}\sum_{i_{1}=2}^{n}\Theta_{i_{1}}^{i_{1}}$.
Assume that \eqref{5.5} is true for $2\leq k\leq q-1$. By using
\eqref{12}, compute the covariant derivative on $q$-chern form:
\begin{eqnarray}
\nabla_{i}P^{q}(\Theta)&=&\frac{1}{q}\sum_{j=1}^{q}(-1)^{j-1}
\nabla_{i}(\Theta_{i_{1}}^{i_{2}}\Theta_{i_{2}}^{i_{3}}\cdots\Theta_{i_{j}}^{i_{1}})P^{q-j}(\Theta)
\nonumber\\&&+\frac{1}{q}\sum_{j=1}^{q}(-1)^{j-1}
\Theta_{i_{1}}^{i_{2}}\Theta_{i_{2}}^{i_{3}}\cdots\Theta_{i_{j}}^{i_{1}}\nabla_{i}P^{q-j}(\Theta).\label{5.8}
\end{eqnarray}
By using induction hypothesis on the second term of \eqref{5.8},
we get \begin{eqnarray}
\lefteqn{\frac{1}{q}\sum_{j=1}^{q}(-1)^{j-1}
\Theta_{i_{1}}^{i_{2}}\Theta_{i_{2}}^{i_{3}}\cdots\Theta_{i_{j}}^{i_{1}}\nabla_{i}P^{q-j}(\Theta)}\nonumber\\&&=
\frac{1}{q}\sum_{j=1}^{q}(-1)^{j-1}\Theta_{i_{1}}^{i_{2}}\cdots\Theta_{i_{j}}^{i_{1}}
\sum_{\ell=1}^{q-j}(-1)^{\ell-1}\frac{1}{\ell}
\nabla_{i}(\Theta_{i_{j+1}}^{i_{j+2}}\cdots\Theta_{i_{j+\ell}}^{i_{j+1}})P^{q-j-\ell}(\Theta)\nonumber\\&&=
\frac{1}{q}\sum_{\ell=1}^{q}(-1)^{\ell-1}\frac{1}{\ell}\nabla_{i}
(\Theta_{i_{j+1}}^{i_{j+2}}\cdots\Theta_{i_{j+\ell}}^{i_{j+1}})\sum_{j=1}^{q-\ell}(-1)^{j-1}\Theta_{i_{1}}^{i_{2}}\cdots\Theta_{i_{j}}^{i_{1}}P^{q-j-\ell}(\Theta)\nonumber\\&&=
\frac{1}{q}\sum_{\ell=1}^{q}(-1)^{\ell-1}\frac{q-\ell}{\ell}\nabla_{i}
(\Theta_{i_{1}}^{i_{2}}\cdots\Theta_{i_{\ell}}^{i_{1}})P^{q-\ell}(\Theta)\label{5.9}
\end{eqnarray}
By substituting \eqref{5.9} back to \eqref{5.8}, we obtain
\eqref{5.5}.
\end{proof}

\begin{proof}[Proof of Sublemm \ref{sub5.1}.]
Prove by induction. For $q=2$, by using \eqref{5.3}, we get
\[X^{i}\nabla_{i}\Theta_{j}^{j}=X^{i}\nabla_{j}\Theta_{i}^{j}=X^{j}\nabla_{i}\Theta_{j}^{i}.\]
Assume that \eqref{5.3} holds for $k=2,\cdots,q-1$. By using
\eqref{5.5}, we compute
\begin{multline}
\sum_{i_{1}=2}^{n}X^{i_{1}}\nabla_{i_{1}}P^{q-1}(\Theta)
-\sum_{j=2}^{q}\sum_{i_{1},\cdots,i_{j}=2}^{n}(-1)^{j}
X^{i_{2}}\nabla_{i_{1}}(\Theta_{i_{2}}^{i_{3}}\cdots\Theta_{i_{j}}^{i_{1}})P^{q-j}(\Theta)\\-
\sum_{j=2}^{q}\sum_{i_{1},\cdots,i_{j}=2}^{n}(-1)^{j}
X^{i_{2}}\Theta_{i_{2}}^{i_{3}}\cdots\Theta_{i_{j}}^{i_{1}}\nabla_{i_{1}}P^{q-j}(\Theta)\\=
\sum_{j=2}^{q}\sum_{i_{1},\cdots,i_{j}=2}^{n}(-1)^{j}\frac{1}{j-1}X^{i_{1}}
\nabla_{i_{1}}(\Theta_{i_{2}}^{i_{3}}\cdots\Theta_{i_{j}}^{i_{2}})P^{q-j}(\Theta)\hspace{2.9cm}\\
-\sum_{j=2}^{q}\sum_{i_{1},\cdots,i_{j}=2}^{n}(-1)^{j}X^{i_{2}}\nabla_{i_{1}}
(\Theta_{i_{2}}^{i_{3}}\cdots\Theta_{i_{j}}^{i_{1}})P^{q-j}(\Theta)\\-
\sum_{j=2}^{q}\sum_{i_{1},\cdots,i_{j}=2}^{n}(-1)^{j}X^{i_{2}}\Theta_{i_{2}}^{i_{3}}\cdots\Theta_{i_{j}}^{i_{1}}
\sum_{\ell=1}^{q-j}(-1)^{\ell-1}\frac{1}{\ell}\nabla_{i_{1}}(\Theta_{i_{j+1}}^{i_{j+2}}
\cdots\Theta_{i_{j+\ell}}^{i_{j+1}})P^{q-j}(\Theta).\label{5.10}
\end{multline}
In fact, we have \begin{equation}
\nabla_{i_{1}}\big(\sum_{i_{2},\cdots,i_{j}=2}^{n}\Theta_{i_{2}}^{i_{3}}\cdots\Theta_{i_{j}}^{i_{2}}\big)
=(j-1)\sum_{i_{2},\cdots,i_{j}=2}^{n}(\nabla_{i_{1}}\Theta_{i_{2}}^{i_{3}})\cdots\Theta_{i_{j}}^{i_{2}},\label{5.11}
\end{equation}
and \begin{multline}\nabla_{i_{1}}
(\Theta_{i_{2}}^{i_{3}}\cdots\Theta_{i_{j}}^{i_{1}})=
\big(\nabla_{i_{1}}\Theta_{i_{2}}^{i_{3}}\big)\Theta_{i_{3}}^{i_{4}}\cdots\Theta_{i_{j}}^{i_{1}}+
\sum_{k=3}^{j-1}\Theta_{i_{2}}^{i_{3}}\cdots\Theta_{i_{k-1}}^{i_{k}}
\big(\nabla_{i_{1}}\Theta_{i_{k}}^{i_{k+1}}\big)\cdots\Theta_{i_{j}}^{i_{1}}\\+
\Theta_{i_{2}}^{i_{3}}\cdots\Theta_{i_{j-1}}^{i_{j}}\big(\nabla_{i_{1}}\Theta_{i_{j}}^{i_{1}}\big).\hspace{4cm}\label{5.111}
\end{multline}
Use \eqref{5.11} in the first term on the right hand side of
\eqref{5.10} and use \eqref{5.111} in the second term of
\eqref{5.10}. We get
\begin{multline}
\sum_{i_{1},i_{2}=2}^{n}X^{i_{1}}\nabla_{i_{1}}\Theta_{i_{2}}^{i_{2}}P^{q-2}(\Theta)+
\sum_{j=3}^{q}\sum_{i_{1},\cdots,i_{j}=2}^{n}(-1)^{j}
X^{i_{1}}\big(\nabla_{i_{1}}\Theta_{i_{2}}^{i_{3}}\big)\cdots\Theta_{i_{j}}^{i_{2}}P^{q-j}\\
-\sum_{j=2}^{q}\sum_{i_{1},\cdots,i_{j}=2}^{n}(-1)^{j}X^{i_{2}}
\big(\nabla_{i_{1}}\Theta_{i_{2}}^{i_{3}}\big)\cdots\Theta_{i_{j}}^{i_{1}}P^{q-j}(\Theta)\hspace{4cm}\\-
\sum_{j=2}^{q}\sum_{i_{1},\cdots,i_{j}=2}^{n}(-1)^{j}\sum_{k=3}^{j-1}
X^{i_{2}}\Theta_{i_{2}}^{i_{3}}\cdots\Theta_{i_{k-1}}^{i_{k}}
\big(\nabla_{i_{1}}\Theta_{i_{k}}^{i_{k+1}}\big)\cdots\Theta_{i_{j}}^{i_{1}}P^{q-j}(\Theta)\hspace{1.1cm}\\
-\sum_{j=2}^{q}\sum_{i_{1},\cdots,i_{j}=2}^{n}(-1)^{j}
X^{i_{2}}\Theta_{i_{2}}^{i_{3}}\cdots\big(\nabla_{i_{1}}\Theta_{i_{j}}^{i_{1}}\big)P^{q-j}(\Theta)\hspace{4cm}\\-
\sum_{j=2}^{q}\sum_{i_{1},\cdots,i_{j+\ell}=2}^{n}(-1)^{j}X^{i_{2}}\Theta_{i_{2}}^{i_{3}}\cdots\Theta_{i_{j}}^{i_{1}}
\sum_{\ell=1}^{q-j}(-1)^{\ell+1}\big(\nabla_{i_{1}}\Theta_{i_{j+1}}^{i_{j+2}}\big)\cdots
\Theta_{i_{j+\ell}}^{i_{j+1}}P^{q-j}(\Theta).\label{5.12}
\end{multline}
By using \eqref{5.4} on the third term of \eqref{5.12}, the first
three terms of \eqref{5.12} are cancelled out. Change the indices
in the last term of \eqref{5.12}
\begin{eqnarray}
\lefteqn{-\sum_{j=2}^{q}\sum_{i_{1},\cdots,i_{j+\ell}=2}^{n}(-1)^{j}X^{i_{2}}\Theta_{i_{2}}^{i_{3}}\cdots\Theta_{i_{j}}^{i_{1}}
\sum_{\ell=1}^{q-j}(-1)^{\ell+1}\big(\nabla_{i_{1}}\Theta_{i_{j+1}}^{i_{j+2}}\big)\cdots
\Theta_{i_{j+\ell}}^{i_{j+1}}P^{q-j}(\Theta)}\nonumber\\&&=
-\sum_{j=2}^{q}\sum_{i_{1},\cdots,i_{j}=2}^{n}\sum_{k=3}^{j-1}(-1)^{j-1}
X^{i_{2}}\Theta_{i_{2}}^{i_{3}}\cdots\Theta_{i_{k-1}}^{i_{1}}\big(\nabla_{i_{1}}\Theta_{i_{k}}^{i_{k+1}}\big)\cdots
\Theta_{i_{j}}^{i_{k}}P^{q-j}(\Theta)\nonumber\\&&\quad
-\sum_{j=2}^{q}\sum_{i_{1},\cdots,i_{j}=2}^{n}(-1)^{j-1}
X^{i_{2}}\Theta_{i_{2}}^{i_{3}}\cdots\Theta_{i_{j-1}}^{i_{1}}\big(\nabla_{i_{1}}\Theta_{i_{j}}^{i_{j}}\big)
P^{q-j}(\Theta).\label{5.13}
\end{eqnarray}
By using \eqref{5.4} in \eqref{5.13}, we obtain
\begin{eqnarray}
\lefteqn{-\sum_{j=2}^{q}\sum_{i_{1},\cdots,i_{j+\ell}=2}^{n}(-1)^{j}X^{i_{2}}\Theta_{i_{2}}^{i_{3}}\cdots\Theta_{i_{j}}^{i_{1}}
\sum_{\ell=1}^{q-j}(-1)^{\ell+1}\big(\nabla_{i_{1}}\Theta_{i_{j+1}}^{i_{j+2}}\big)\cdots
\Theta_{i_{j+\ell}}^{i_{j+1}}P^{q-j}(\Theta)}\nonumber\\&&=
\sum_{j=2}^{q}\sum_{i_{1},\cdots,i_{j}=2}^{n}\sum_{k=3}^{j}(-1)^{j-1}
X^{i_{2}}\Theta_{i_{2}}^{i_{3}}\cdots\Theta_{i_{k-1}}^{i_{1}}\big(\nabla_{i_{k}}\Theta_{i_{1}}^{i_{k+1}}\big)\cdots
\Theta_{i_{j}}^{i_{k}}P^{q-j}(\Theta)\nonumber\\&&\quad
+\sum_{j=2}^{q}\sum_{i_{1},\cdots,i_{j}=2}^{n}(-1)^{j-1}
X^{i_{2}}\Theta_{i_{2}}^{i_{3}}\cdots\Theta_{i_{j}}^{i_{1}}\big(\nabla_{i_{j}}\Theta_{i_{1}}^{i_{j}}\big)P^{q-j}(\Theta).
\label{5.14}
\end{eqnarray}
By substituting \eqref{5.14} back to \eqref{5.12}, the last three
terms on of \eqref{5.12} are cancelled out. Therefore,
\eqref{5.12} is equal to $0$.
\end{proof}
\end{proof}
By Lemma \ref{lem5.1},we have
\[\int_{M}q\tilde{P}^{q}(\nabla
X,\Theta,\cdots,\Theta)\wedge\omega_{M}^{n-q}=\int_{M}\rm{div}(Y)\wedge\omega_{M}^{n-q}.\]
\begin{theorem}[General divergence theorem]
Given a holomorphic vector field with $(q-1,q-1)$-valued form $Y$
defined in Lemma \ref{lem5.1}, we have
\[\int_{M}\rm{div}(Y)\wedge\omega_{M}^{n-q}=\int_{M}L_{Y}(\omega^{n-q})=0.\]
\end{theorem}
\begin{proof}
Let $\tilde{Y}=
\sum_{i_{1}=2}^{n}\sum_{j=2}^{q}\sum_{i_{2},\cdots,i_{j}=2}^{n}
(-1)^{j-1}X^{i_{2}}\Theta_{i_{2}}^{i_{3}}\cdots\Theta_{i_{j}}^{i_{1}}P^{q-j}(\Theta)\otimes\frac{\partial}{\partial
z_{i_{1}}}$. We have $Y=
\sum_{i_{1}=2}^{n}X^{i_{1}}P^{q-1}(\Theta)\otimes\frac{\partial}{\partial
z_{i_{1}}}+\tilde{Y}$. Compute
\begin{multline*}
\int_{M}L_{\tilde{Y}}(\omega^{n-q})\\=(n-q)\int_{M}\sum_{i_{2},\cdots,i_{j}=2}^{n}
(-1)^{j-1}X^{i_{2}}\Theta_{i_{2}}^{i_{3}}\cdots
\Theta_{i_{j}}^{i_{1}}\wedge(\tilde{g}_{i_{1}\bar{\lambda}}d\bar{z}_{\lambda})\wedge
P^{q-j}(\Theta)\wedge\omega^{n-1-q}.\end{multline*} On K\"ahler
manifolds, we have
\begin{eqnarray*}\sum_{i_{1},\lambda=2}^{n}\Theta_{i_{j}}^{i_{1}}\wedge(\tilde{g}_{i_{1}\bar{\lambda}}d\bar{z}_{\lambda})
&=&\frac{\sqrt{-1}}{2\pi}\sum_{i_{1},\alpha,\beta,\lambda=2}^{n}R_{i_{j}\bar{\lambda}\alpha\bar{\beta}}dz_{\alpha}\wedge
d\bar{z}_{\beta}\wedge
d\bar{z}_{\lambda}\\&=&\frac{\sqrt{-1}}{2\pi}\sum_{i_{1},\alpha,\beta,\lambda=2}^{n}R_{i_{j}\bar{\beta}\alpha\bar{\lambda}}dz_{\alpha}\wedge
d\bar{z}_{\beta}\wedge d\bar{z}_{\lambda}=0.\end{eqnarray*}
Therefore, we obtain $\int_{M}L_{\tilde{Y}}(\omega^{n-q})=0$. By
the linearity of vector addition and \eqref{2.63}, we have
\begin{eqnarray}\int_{M}L_{Y}(\omega^{n-q})&=&\int_{M}L_{P^{q-1}(\Theta)\otimes X}(\omega^{n-q})
\nonumber\\&=&\int_{M}P^{q-1}(\Theta)\wedge
L_{X}(\omega^{n-q})\nonumber\\&=&
\int_{M}L_{X}(P^{q-1}(\Theta)\wedge
\omega^{n-q}).\label{5.15}\end{eqnarray} By \eqref{5.15}, Lemma
\ref{lem5.1} and Lemma \ref{lem2.7}, we obtain
\begin{eqnarray*}
\int_{M}L_{Y}(\omega^{n-q})&=&\int_{M}L_{X}(P^{q-1}(\Theta)\wedge
\omega^{n-q})\\&=&q\int_{M}\tilde{P}^{q}(\nabla
X,\Theta,\cdots,\Theta)\wedge\omega^{n-q}\\&=&\int_{M}{\rm
div}(Y)\wedge\omega^{n-q}\\&=&0.
\end{eqnarray*}
\end{proof}
\bibliographystyle{abbrv}
\bibliography{bando-futaki}

\begin{thebibliography}{10}

\bibitem{BS}
S.~Bando.
\newblock An obstruction for chern class forms to be harmonic.
\newblock unpublished, 1983.

\bibitem{B-M}
S.~Bando and T.~Mabuchi.
\newblock On some integral invariants on complex manifolds. {I}.
\newblock {\em Proc. Japan Acad. Ser. A Math. Sci.}, 62(5):197--200, 1986.

\bibitem{CE}
E.~Calabi.
\newblock Extremal {K}\"ahler metrics. {II}.
\newblock In {\em Differential geometry and complex analysis}, pages 95--114.
  Springer, Berlin, 1985.

\bibitem{C-T}
X.~X. Chen and G.~Tian.
\newblock Ricci flow on {K}\"ahler-{E}instein surfaces.
\newblock {\em Invent. Math.}, 147(3):487--544, 2002.

\bibitem{D-T}
W.~Y. Ding and G.~Tian.
\newblock K\"ahler-{E}instein metrics and the generalized {F}utaki invariant.
\newblock {\em Invent. Math.}, 110(2):315--335, 1992.

\bibitem{FA1}
A.~Futaki.
\newblock An obstruction to the existence of {E}instein {K}\"ahler metrics.
\newblock {\em Invent. Math.}, 73(3):437--443, 1983.

\bibitem{FA2}
A.~Futaki.
\newblock {\em K\"ahler-{E}instein metrics and integral invariants}, volume
  1314 of {\em Lecture Notes in Mathematics}.
\newblock Springer-Verlag, Berlin, 1988.

\bibitem{FA3}
A.~Futaki.
\newblock Asymptotic {C}how semi-stability and integral invariants.
\newblock {\em Internat. J. Math.}, 15(9):967--979, 2004.

\bibitem{F-M}
A.~Futaki and S.~Morita.
\newblock Invariant polynomials of the automorphism group of a compact complex
  manifold.
\newblock {\em J. Differential Geom.}, 21(1):135--142, 1985.

\bibitem{G-H}
P.~Griffiths and J.~Harris.
\newblock {\em Principles of algebraic geometry}.
\newblock Wiley-Interscience [John Wiley \& Sons], New York, 1978.
\newblock Pure and Applied Mathematics.

\bibitem{KS1}
S.~Kobayashi.
\newblock On compact {K}\"ahler manifolds with positive definite {R}icci
  tensor.
\newblock {\em Ann. of Math. (2)}, 74:570--574, 1961.

\bibitem{LNC}
N.~C. Leung.
\newblock Bando {F}utaki invariants and {K}\"ahler {E}instein metric.
\newblock {\em Comm. Anal. Geom.}, 6(4):799--808, 1998.

\bibitem{LZ1}
Z.~Lu.
\newblock On the {F}utaki invariants of complete intersections.
\newblock {\em Duke Math. J.}, 100(2):359--372, 1999.

\bibitem{LZ2}
Z.~Lu.
\newblock {$K$} energy and {$K$} stability on hypersurfaces.
\newblock {\em Comm. Anal. Geom.}, 12(3):601--630, 2004.

\bibitem{MT}
T.~Mabuchi.
\newblock {$K$}-energy maps integrating {F}utaki invariants.
\newblock {\em Tohoku Math. J. (2)}, 38(4):575--593, 1986.

\bibitem{P-S1}
D.~H. Phong and J.~Sturm.
\newblock Stability, energy functionals, and {K}\"ahler-{E}instein metrics.
\newblock {\em Comm. Anal. Geom.}, 11(3):565--597, 2003.

\bibitem{P-S2}
D.~H. Phong and J.~Sturm.
\newblock The {F}utaki invariant and the {M}abuchi energy of a complete
  intersection.
\newblock {\em Comm. Anal. Geom.}, 12(1-2):321--343, 2004.

\bibitem{TG1}
G.~Tian.
\newblock The {$K$}-energy on hypersurfaces and stability.
\newblock {\em Comm. Anal. Geom.}, 2(2):239--265, 1994.

\bibitem{TG2}
G.~Tian.
\newblock K\"ahler-{E}instein metrics with positive scalar curvature.
\newblock {\em Invent. Math.}, 130(1):1--37, 1997.

\bibitem{TG3}
G.~Tian.
\newblock Bott-{C}hern forms and geometric stability.
\newblock {\em Discrete Contin. Dynam. Systems}, 6(1):211--220, 2000.

\bibitem{T-Z}
G.~Tian and X.~Zhu.
\newblock A new holomorphic invariant and uniqueness of {K}\"ahler-{R}icci
  solitons.
\newblock {\em Comment. Math. Helv.}, 77(2):297--325, 2002.

\bibitem{WB}
B.~Weinkove.
\newblock Higher {K}-energy functionals and highier futaki invariants.
\newblock Preprint, 2002.

\bibitem{YM}
M.~Yotov.
\newblock Generalized futaki invariant of almost fano toric varieties,
  examples.
\newblock Preprint, 1999.

\end{thebibliography}
\end{document}